\documentclass[a4paper,fleqn]{cas-sc}

\usepackage[authoryear,longnamesfirst]{natbib}
\usepackage{amsmath,amssymb,amsthm,mathtools}
\usepackage{bm}
\usepackage{graphicx}
\usepackage{booktabs}
\usepackage{siunitx}
\usepackage{subcaption}
\usepackage{float}
\usepackage{enumitem}
\usepackage{xcolor}
\usepackage{array}
\usepackage[capitalize,nameinlink]{cleveref}
\usepackage{placeins}
\usepackage{setspace}
\usepackage{lineno}

\newcolumntype{L}[1]{>{\raggedright\arraybackslash}p{#1}}

\hypersetup{
    colorlinks=true,
    linkcolor=blue,
    citecolor=blue,
    urlcolor=blue
}
\usepackage{xcolor}

\newif\ifmarked
\markedfalse    

\newcommand{\revone}[1]{%
  \ifmarked\textcolor{blue}{#1}\else#1\fi}

\newcommand{\revtwo}[1]{%
  \ifmarked\textcolor{magenta}{#1}\else#1\fi}

\newcommand{\revboth}[1]{%
  \ifmarked\textcolor{teal}{#1}\else#1\fi}

\begin{document}
\let\WriteBookmarks\relax
\def\floatpagepagefraction{1}
\def\textpagefraction{.001}
\doublespacing

\shorttitle{Orbital Networks in the Three-Body Problem}
\shortauthors{A. Braik \& S.D. Ross}

\title [mode = title]{Orbital Networks in the Three-Body Problem}

\author[1]{Abdullah Braik}[orcid=0009-0007-8207-2997]
\cormark[1]
\fnmark[1]
\ead{braik@vt.edu}
\credit{Conceptualization; Methodology; Software; Validation; Formal analysis; Investigation; Data curation; Visualization; Writing -- original draft; Writing -- review and editing}

\author[1]{Shane D. Ross}[orcid=0000-0001-5523-2376]
\fnmark[2]
\ead{sdross@vt.edu}
\credit{Conceptualization; Methodology; Resources; Supervision; Project administration; Funding acquisition; Writing -- review and editing}

\affiliation[1]{organization={Aerospace and Ocean Engineering, Virginia Tech},
            addressline={495 Old Turner St},
            city={Blacksburg},
            state={VA},
            country={USA}
            }

\fntext[1]{PhD Student, Aerospace and Ocean Engineering, Virginia Tech, Blacksburg, VA, USA, braik@vt.edu}
\fntext[2]{Professor, Aerospace and Ocean Engineering, Virginia Tech, Blacksburg, VA, USA, sdross@vt.edu}

\begin{abstract}
Orbital transfers in multi-body systems are often studied as isolated trajectory design problems, making it difficult to identify the larger transport structure connecting families of periodic orbits, including which families act as hubs, gateways, relays, or persistently difficult-to-access regions. This work introduces a reachable-set-based framework for constructing orbital networks in the circular restricted three-body problem. Finite-$\Delta V$ and finite-time-of-flight reachable-set overlaps are used to infer accessibility relationships between representative periodic orbit families on a common Jacobi energy manifold and to assemble these relationships into a weighted orbital network. Applied to the Earth–Moon system, the resulting network reveals distinct accessibility regimes in which direct reachability, graph connectedness, and feasible multileg closure emerge separately. The analysis identifies multi-orbiter cycler orbits as the dominant hub, gateway, and relay families, with the (3,2)-cycler dominating across much of the sampled budget plane and the short-period (1,1)-cycler dominating in the low-time-of-flight regime, while the stable 2:1 resonant orbit remains persistently difficult to access. Although the maximum-budget network is nearly complete in a binary sense, its weighted accessibility remains strongly non-uniform. Selected proxy-supported connections are refined into concrete trajectories through differential correction, with corrected transfer costs remaining below the proxy estimates in all tested cases. Together, the results demonstrate how reachable-set overlap geometry can expose large-scale transport structure in nonlinear gravitational systems without requiring exhaustive pairwise trajectory optimization.
\end{abstract}




   \nonumnote{\copyright\ 2026. This manuscript version is made available under the CC-BY-NC-ND 4.0 license, \url{https://creativecommons.org/licenses/by-nc-nd/4.0/}. Published version: \url{https://doi.org/10.1016/j.asr.2026.10.100}}
   
\maketitle
\newpage
\tableofcontents
\newpage

\section{Introduction}

The cislunar region is evolving from a domain of primarily exploratory activity into an increasingly sustained operational environment.
Planned infrastructure, including NASA's lunar base, orbital platforms, scientific assets, and logistics support elements, will require repeated movement between distinct orbit regimes in the Earth--Moon system under explicit maneuver and time constraints \citep{whitley2016,cheetham2022}.
However, existing approaches do not yet provide a unified family-level view of how accessibility is organized across many coexisting cislunar orbit families under jointly varying maneuver and time budgets.
In this setting, the central question extends beyond designing transfers between preselected orbits to understanding the accessibility structure of  the broader cislunar orbit family landscape.
This shift motivates a transport architecture perspective that complements pairwise low-energy transfer design by exposing structural features, such as dominant access families, persistent hard-access regions, and the regimes in which multileg routing becomes operationally admissible, that are not visible from individual transfer studies alone.

The circular restricted three-body problem (CR3BP) offers a natural framework for studying Earth--Moon motion \citep{szebehely1967,koon2011}, and at a fixed Jacobi energy the region contains many coexisting periodic-orbit families, including Lyapunov, halo, resonant, cycler-type, and distant prograde orbits \citep{howell1984,broucke1968,henon1969,henon1997,rawat2026,ross2025cyclers,ross2026cyclers}.
Their geometry, stability, and bifurcation structure have been studied extensively \citep{doedel2007,henon1997}, with modern catalogs providing increasingly systematic coverage of the Earth--Moon family landscape \citep{leiva2006control,guzzetti2016,zimovan2017}.
Taken together, these families define the candidate set of hub and destination orbits for cislunar transport and constellation analysis.

Transfer design between selected periodic orbit families in the CR3BP is well established.
Invariant manifold methods have produced low-cost heteroclinic and near-heteroclinic connections between Lagrange point orbits \citep{koon2000,gomez2004}, with later extensions to inter-family transfers, resonant-orbit-assisted design, low-energy lunar trajectories, and operationally relevant structures such as near-rectilinear halo orbits (NRHOs) and butterfly orbits \citep{haapala2016,vaquero2014,parker2014,mccarthy2023,davis2017,zimovanspreen2022,braik2025}.
A closely related study by Capdevila and Howell organized Earth–Moon transfer options into a network \citep{capdevila2018}; 
the present work differs in that finite-$\Delta V$ reachable-set overlap on a common Jacobi manifold is used to characterize accessibility across many families simultaneously, and the resulting network is analyzed structurally to identify family-level roles rather than constructing selected transfer pathways.
Thus, while prior work provides mature tools for constructing selected transfers, it remains primarily pairwise and does not directly reveal the broader accessibility structure of the cislunar family set.

Graph-based and set-oriented methods provide a complementary way to organize transport information.
Prior studies have organized CR3BP transfer options as weighted graphs searched with shortest-path methods \citep{tsirogiannis2012}, with recent extensions using motion primitives, sampling-based kinodynamic planning, hierarchical tree search, and lobe-dynamics sequencing to improve automation and scalability \citep{smith2022,smith2023,bruchko2025,spear2026,hiraiwa2026}.
These methods demonstrate the value of graph abstractions for candidate transfer generation, but the graph is often used primarily as a computational tool for point-to-point design.
Set-oriented transport theory offers a complementary view, inferring transport structure from the propagation of sets on a discretized phase space \citep{dellnitz2005,DeJuLoMaPaPrRoTh2005,dellnitz2006target,JeJuRo2009}, and recent cislunar reachable-set work highlights the value of set-based descriptions under control constraints \citep{bowerfind2024}.
Together, these ideas motivate using local reachable sets to populate the edges of a family-to-family cislunar transport network.

Three issues remain insufficiently addressed in this body of work: family-to-family accessibility has not been characterized globally across a representative family set under jointly varying maneuver and time budgets; the cislunar family set has not been interpreted as a network in the stronger structural sense needed to identify hub, gateway, and relay roles; and the distinction between local pairwise reachability, graph connectedness, and budget-feasible multileg closure has not been systematically examined despite corresponding to operationally distinct regimes.

This paper makes three contributions.
First, it introduces a reachable-set based proxy 
accessibility cost for pairs of representative cislunar orbit families: 
a scalar derived from finite-$\Delta V$, finite-time-of-flight reachable-set overlap on a common Jacobi manifold, without constructing a continuous trajectory between the two families. The resulting family-pair matrix functions as a $\Delta V$ map across the representative cislunar set, 
extended from the conventional source-to-destination form to span a full family network and carry network-level structure such as centrality measures and budget-feasible multileg routes.
Second, it identifies the (3,2)-cycler family as a dominant direct-access, gateway, and relay family in the resulting cislunar transport network, with the (1,1)a-cycler dominant at low time-of-flight budgets, and the 2:1 stable resonant orbit as the persistent hard-access family.
Third, it separates direct accessibility, graph connectedness, and budget-feasible multileg closure into three distinct operational regimes (corresponding respectively to time-critical retasking, sustained coverage, and time-flexible redistribution) that pairwise transfer design alone does not expose.
The resulting network structure is interpretable in dynamical-systems terms: low-cost direct and multileg accessibility emerge from unstable representatives sharing a common chaotic region on the Jacobi manifold, while the persistent hard-access edge corresponds to representatives embedded in stable tori that resist low-energy invasion from that chaotic sea.

To support these contributions, this paper develops a family-level cislunar transport framework in the planar Earth--Moon CR3BP.
Representative periodic orbits from multiple coexisting families are used to construct local finite-$\Delta V$, finite-time-of-flight reachable-set atlases in a reduced three-dimensional 
phase-space description  along the Jacobi energy manifold.
Forward reachable sets are generated by energy-preserving heading-change maneuvers applied at points distributed along each representative orbit by arc length, while backward reachable information is obtained efficiently through time-reversal symmetry.
Pairwise family-to-family accessibility is then inferred from overlap between the forward reachable set of one family and the backward reachable set of another on a common discretized phase-space grid.
These pairwise proxy measures are assembled into a cislunar family network, which is used to identify hub, gateway, relay, and bottleneck roles and to study how those roles reorganize over varying budgets.
In this way, the reduced model functions as a screening layer that surfaces structural insight from proxy accessibility trends and identifies candidate connections for subsequent detailed correction,  rather than to solve individual transfer optimization problems directly.
Finally, selected proxy-supported family connections are realized through differential correction, providing corrected direct and relay trajectory examples that ground the network interpretation in concrete transfer solutions.

The remainder of the paper develops the CR3BP formulation, representative family set, reachable-set atlas construction, and overlap-based accessibility metrics.
These tools are then used to build and analyze the resulting cislunar family network across maneuver and time budgets, followed by selected corrected trajectory realizations and concluding remarks.

\section{Dynamical Framework, Symmetry, and Representative Periodic Orbits}\label{sec:framework}

\subsection{Planar Earth--Moon CR3BP Framework}

The analysis is formulated in the planar Earth--Moon CR3BP using the uniformly rotating synodic frame and standard nondimensional normalization; standard derivations of the governing equations are given in \cite{szebehely1967} and \cite{koon2011}. 
The primaries are located at \((-\mu,0)\) and \((1-\mu,0)\), and the planar state is \((x,y,\dot{x},\dot{y})\), as shown schematically in Figure~\ref{fig:cr3bp_schematic}. 
The distances from the particle to the Earth and Moon are, respectively,
\begin{equation}
r_1 = \sqrt{(x+\mu)^2 + y^2},
\qquad
r_2 = \sqrt{(x-1+\mu)^2 + y^2}.
\label{eq:r1r2}
\end{equation}

The equations of motion in the rotating frame are,
\begin{equation}
\ddot{x} - 2\dot{y} = -\bar U_x(x,y),
\qquad
\ddot{y} + 2\dot{x} = -\bar U_y(x,y),
\label{eq:cr3bp_planar_eom}
\end{equation}
where the effective potential is,
\begin{equation}
\bar U(x,y)
=
-\frac{1}{2}\left(x^2+y^2\right)
-
\frac{1-\mu}{r_1}
-
\frac{\mu}{r_2}.
\label{eq:cr3bp_potential}
\end{equation}
Its first derivatives are,
\begin{equation}
\bar U_x
=
-x + \frac{(1-\mu)(x+\mu)}{r_1^3} + \frac{\mu(x-1+\mu)}{r_2^3},
\qquad
\bar U_y
=
-y + \frac{(1-\mu)y}{r_1^3} + \frac{\mu y}{r_2^3}.
\label{eq:cr3bp_potential_derivs}
\end{equation}

The planar CR3BP admits the Jacobi integral,
\begin{equation}
C_J = -2\bar U(x,y) - \left(\dot{x}^2+\dot{y}^2\right),
\label{eq:jacobi_integral}
\end{equation}
which defines the fixed three-dimensional energy level used throughout the later accessibility analysis.

The physical and normalization quantities used in the study are summarized in Table~\ref{tab:em_cr3bp_constants}; the nondimensional Earth and Moon radii are used later in the admissible-domain definition.

\par\medskip
\begin{center}
\captionsetup{type=table}
\centering
\caption{\footnotesize Physical and normalization quantities used in the Earth--Moon CR3BP formulation.}
\label{tab:em_cr3bp_constants}
\begin{tabular}{llll}
\toprule
Quantity & Symbol & Physical value & Nondimensional value \\
\midrule
Earth--Moon mass parameter & $\mu$ & --- & 0.012150584270572 \\
Distance unit & $LU$ & \SI{384400}{km} & 1 \\
Lunar sidereal period & $T_{\mathrm{EM}}$ & \SI{27.321661}{day} & $2\pi$ \\
Time unit & $TU$ & \SI{4.34837740}{day} & 1 \\
Earth radius & $R_E$ & \SI{6378}{km} & 0.0165921 \\
Moon radius & $R_M$ & \SI{1737}{km} & 0.00451873 \\
\bottomrule
\end{tabular}
\end{center}
\medskip

\par\medskip
\begin{center}
\captionsetup{type=figure}
    \centering
    \includegraphics[width=0.6\textwidth]{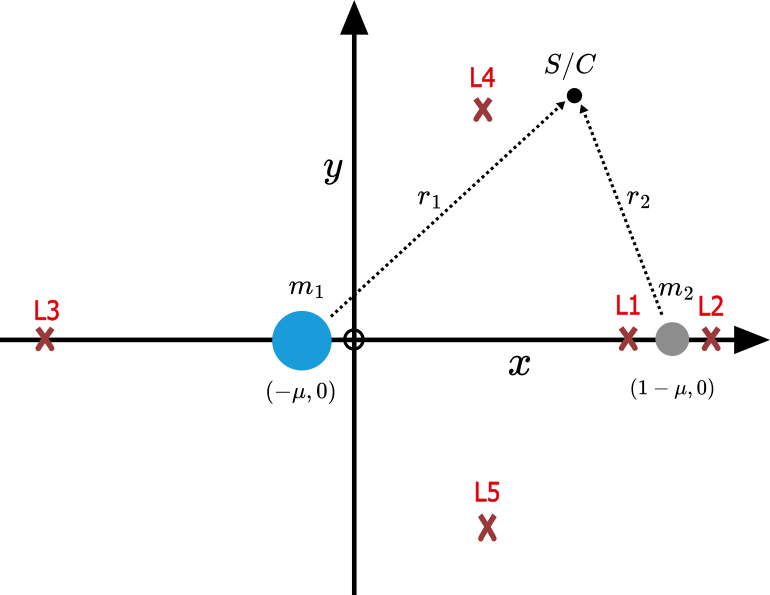}
    \caption{\footnotesize Planar Earth--Moon CR3BP geometry in the rotating synodic frame.}
    \label{fig:cr3bp_schematic}
\end{center}
\medskip


\subsection{Reduced \texorpdfstring{$(x,y,\theta)$}{(x,y,theta)} Formulation on a Fixed Jacobi Constant}\label{subsec:reduced}

For a fixed Jacobi constant \(C_J\), the planar CR3BP state is restricted to a three-dimensional energy manifold within the four-dimensional state space.
The present analysis exploits this constraint by representing the state using planar position and the local velocity heading.
Let the  velocity in the rotating frame be parameterized by its magnitude \(v\) and heading angle \(\theta\) as,

\begin{equation}
\dot{x} = v\cos\theta,
\qquad
\dot{y} = v\sin\theta.
\label{eq:vel_heading}
\end{equation}

The Jacobi integral \eqref{eq:jacobi_integral} then determines the speed from position:
\begin{equation}
v^2 = -2\bar U(x,y) - C_J.
\label{eq:v2_reduced}
\end{equation}
Therefore, \(v\) is no longer an independent state variable, and the planar CR3BP can be written in reduced coordinates \((x,y,\theta)\) as,
\begin{equation}
\dot{x} = v\cos\theta,
\qquad
\dot{y} = v\sin\theta,
\qquad
\dot{\theta}
=
-2 + \frac{\bar U_x(x,y)\sin\theta-\bar U_y(x,y)\cos\theta}{v},
\label{eq:reduced_model}
\end{equation}
with
\begin{equation}
v = \sqrt{-2\bar U(x,y)-C_J}.
\label{eq:reduced_speed}
\end{equation}

The reduced formulation is regular only where,
\begin{equation}
-2\bar U(x,y)-C_J > 0,
\label{eq:energy_admissible}
\end{equation}
so that \(v>0\).
The boundary,
\begin{equation}
2\bar U(x,y)+C_J = 0,
\label{eq:zvc_boundary}
\end{equation}
is the zero-velocity curve.
Although this boundary is part of the energetic boundary of the full planar CR3BP, it is singular in the reduced \((x,y,\theta)\) description because the heading  of motion in \eqref{eq:reduced_model} contains division by \(v\).
Thus, trajectories approaching the zero-velocity curve must be treated with the full planar state model rather than the reduced formulation.

This reduced representation is the state space used for the local reachable-set construction and the subsequent family-to-family accessibility proxy.


\subsection{Time-Reversal Symmetry About the \texorpdfstring{$x$}{x}-Axis}\label{subsec:time-reversal}

The planar CR3BP is reversible under the standard time-reversal symmetry about the \(x\)-axis \citep{szebehely1967,koon2011}.
\revtwo{
Let the full rotating-frame state be
\(\mathbf{X}=(x,y,\dot{x},\dot{y})\).
The corresponding state-space symmetry map is
}
\[
\mathcal{R}(\mathbf{X})
=
\mathcal{R}(x,y,\dot{x},\dot{y})
=
(x,-y,-\dot{x},\dot{y}).
\]

\revtwo{
Thus, if \(\mathbf{X}(t)\) is a solution, then the trajectory
\[
\widetilde{\mathbf{X}}(t)
=
\mathcal{R}\!\left(\mathbf{X}(-t)\right)
\]
is also a solution.
}
The reversed motion is obtained by reflecting the trajectory across the \(x\)-axis while reversing time.
The map preserves the Jacobi integral because \(\bar U(x,-y)=\bar U(x,y)\) and the speed magnitude is unchanged.
Therefore, the time-reversed image of a trajectory remains on the same fixed-\(C_J\) manifold.

For the reduced \((x,y,\theta)\) coordinates, the induced symmetry follows from \eqref{eq:vel_heading}.
Since \((\dot{x},\dot{y})=(v\cos\theta,v\sin\theta)\) maps to \((-\dot{x},\dot{y})\), the transformed heading satisfies,
\begin{equation}
\theta^*=\pi-\theta
\qquad
(\mathrm{mod}\ 2\pi),
\label{eq:theta_reverser_derivation}
\end{equation}
so the induced reduced-coordinate symmetry map is,
\begin{equation}
\mathcal{R}(x,y,\theta)
=
\left(
x,\,-y,\,\pi-\theta
\right),
\label{eq:reduced_reverser}
\end{equation}
with the angle understood modulo \(2\pi\).

This reduced symmetry is used later to convert forward-propagated reachable information into backward reachable information on the same fixed-\(C_J\) manifold.


\subsection{Admissible Cislunar Domain and Zero Velocity Curve}

The accessibility analysis is carried out at the common Jacobi constant \(C_J=3.1294\).
At this energy level, the admissible configuration-space region is determined by \eqref{eq:energy_admissible}, with boundary given by the zero-velocity curve in \eqref{eq:zvc_boundary}.
Figure~\ref{fig:admissible_domain} shows the corresponding admissible cislunar domain.

For the reachable-set construction, this admissible set is intersected with a finite barycenter-centered region.
The distance from the Earth--Moon barycenter is,
\begin{equation}
r = \sqrt{x^2+y^2},
\end{equation}
and let \(R_{\mathrm{dom}}\) denote the domain radius. In this study,
\begin{equation}
R_{\mathrm{dom}} = 1.2
\label{eq:rdom_value}
\end{equation}
is adopted in nondimensional units, so that the working region is restricted to,
\begin{equation}
r \le R_{\mathrm{dom}}.
\label{eq:rdom_condition}
\end{equation}
This restriction focuses the analysis on the interior cislunar region while retaining all representative orbits considered in this study, including the largest \(x\)-excursion of the representative L2 Lyapunov orbit.

The interiors of the Earth and Moon are also excluded:
\begin{equation}
r_1 > R_E,
\qquad
r_2 > R_M,
\label{eq:primary_exclusion}
\end{equation}
where \(R_E\) and \(R_M\) denote the nondimensional Earth and Moon radii listed in Table~\ref{tab:em_cr3bp_constants}.

Accordingly, the admissible working region may be written as,
\begin{equation}
\mathcal{D}(C_J)
=
\left\{
(x,y)\;:\;
r \le R_{\mathrm{dom}},\;
-2\bar U(x,y)-C_J \ge 0,\;
r_1 > R_E,\;
r_2 > R_M
\right\}.
\label{eq:working_domain}
\end{equation}
This set defines the working configuration-space domain used in the reachable-set computations.

As noted above, the zero-velocity curve is retained as the energetic boundary of the admissible domain, although reduced-state propagation is regular only for \(v>0\).

\par\medskip
\begin{center}
\captionsetup{type=figure}
    \centering
    \includegraphics[width=0.50\textwidth]{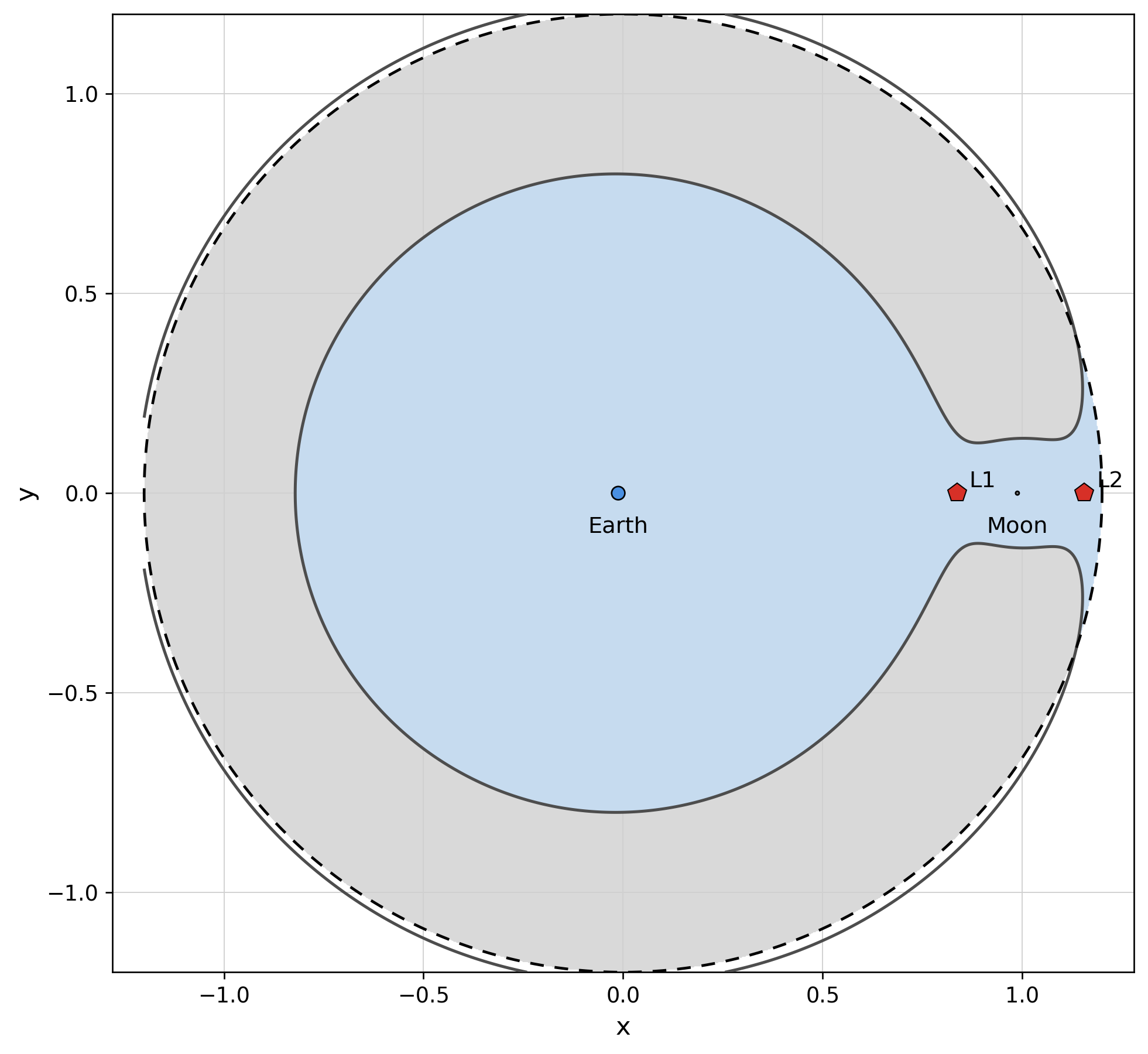}
    \caption{\footnotesize Admissible cislunar domain at the common Jacobi level \(C_J=3.1294\) used throughout the study.
    The retained region is defined by the intersection of the energy-admissible set \(-2\bar U-C_J \ge 0\), the finite barycenter-centered domain \(r \le R_{\mathrm{dom}}=1.2\), and the exclusion of Earth and Moon interiors.}
    \label{fig:admissible_domain}
\end{center}
\medskip


\subsection{Periodic Orbit Family Representatives on a Common Jacobi Energy Manifold}\label{subsec:families}

The accessibility framework is built on thirteen representative planar periodic orbits in the Earth--Moon CR3BP, each selected from a distinct cislunar family at the common Jacobi level \(C_J=3.1294\).
Using a shared Jacobi constant allows differences in accessibility to reflect orbit geometry and local transport structure rather than energy mismatch between representatives.
Unless otherwise noted, the term ``family'' is used below as shorthand for the representative periodic orbit selected from that family at $C_J=3.1294$; references to family-to-family accessibility, family roles, and the family network should therefore be understood in this representative-orbit sense.

\revboth{
The value \(C_J=3.1294\) is the largest common Jacobi constant for which all thirteen representatives considered in this study coexist.
Continuation of the (2,1)-cycler reaches an upper-\(C_J\) fold at approximately this value, so that this family does not provide the required representative at larger \(C_J\).
Selecting the upper end of the common-existence interval also retains a comparatively restrictive cislunar energy manifold: the \(L_1\) and \(L_2\) necks are open, permitting transport through the principal Earth--Moon gateways, while the \(L_3\) neck remains closed.
Thus, the selected value permits a common-energy comparison across the representative family set without moving to lower \(C_J\) levels associated with an additional change in the zero-velocity topology.
}

The selected representatives span Lagrange point, cycler, resonant, and distant prograde motions \citep{ross2025cyclers,ross2026cyclers,rawat2026,JPL_PeriodicOrbits}.
Their periods and stability parameters are listed in Table~\ref{tab:representative_nodes}, and their geometries are shown in Figure~\ref{fig:representative_orbits}.
\revboth{
The thirteen-family set is intended as a diverse, non-exhaustive cross-section of the periodic-orbit structures available on this common Jacobi manifold rather than as a prescribed coverage threshold.
Additional families may be incorporated into the same reachable-set and network construction when required for a particular application.
}

\par\medskip
\begin{center}
\captionsetup{type=table}
\centering
\caption{\footnotesize Representative periodic orbits used in the accessibility study at the common Jacobi level \(C_J=3.1294\).
Periods are reported in days, and \(\sigma\) denotes the instability rate.}
\label{tab:representative_nodes}
\begin{tabular}{lcccc}
\toprule
Orbit Family & Abbreviation & Period $T_{PO}$ & Instability rate $\sigma$ & Instability rate $\sigma_d$ \\
&  & [days]& [TU$^{-1}$]& [day$^{-1}$]                \\
\midrule
L1 Lyapunov             & LL1    & 12.811 & 2.4884 & 0.5722 \\
L2 Lyapunov            & LL2    & 15.117 & 1.9797 & 0.4552 \\
(1,1)a-cycler           & C11a  & 42.140 & 1.0482 & 0.2410 \\
(1,1)b-cycler           & C11b  & 55.995 & 0.9255 & 0.2128\\
(2,1)-cycler            & C21   & 84.533 & 0.1358 & 0.0312\\
(3,2)-cycler            & C32   & 78.613 & 0.6886 & 0.1583\\
2:1 stable resonant     & R21-S & 26.500 & 0 & 0 \\
2:1 unstable resonant   & R21-U & 31.039 & 0.8397 & 0.1931 \\
3:1 stable resonant     & R31-S & 27.252 & 0 & 0 \\
3:1 unstable resonant   & R31-U & 28.066 & 0.40124 & 0.0923 \\
5:2 stable resonant     & R52-S & 54.802 & 0  & 0\\
5:2 unstable resonant   & R52-U & 56.436 & 0.36547 & 0.0840\\
Distant prograde orbit  & DPO   & 11.184 & 1.5886 & 0.3653\\
\bottomrule
\end{tabular}
\end{center}
\medskip

Periodic orbits in the CR3BP differ in both their stability character and their period, so a meaningful comparison across families requires a period-normalized stability measure. For each representative orbit, the instability rate is defined as the 
Floquet exponent,
\begin{equation}
\sigma = \frac{\ln({\lambda_{\max}})}{T_{PO}},
\label{eq:stability_parameter}
\end{equation}
where $\lambda_{\max}$ is the largest magnitude nontrivial monodromy eigenvalue (or Floquet multiplier), and  $T_{PO}$ is the nondimensional period of the corresponding periodic orbit.
Stable orbits, for which $|\lambda_{\max}| = 1$, give $\sigma = 0$; 
larger values of $\sigma$
correspond to faster local divergence per unit nondimensional time.

\par\medskip
\begin{center}
\captionsetup{type=figure}
    \centering

    \begin{subfigure}[t]{0.23\textwidth}
        \centering
        \includegraphics[width=\textwidth]{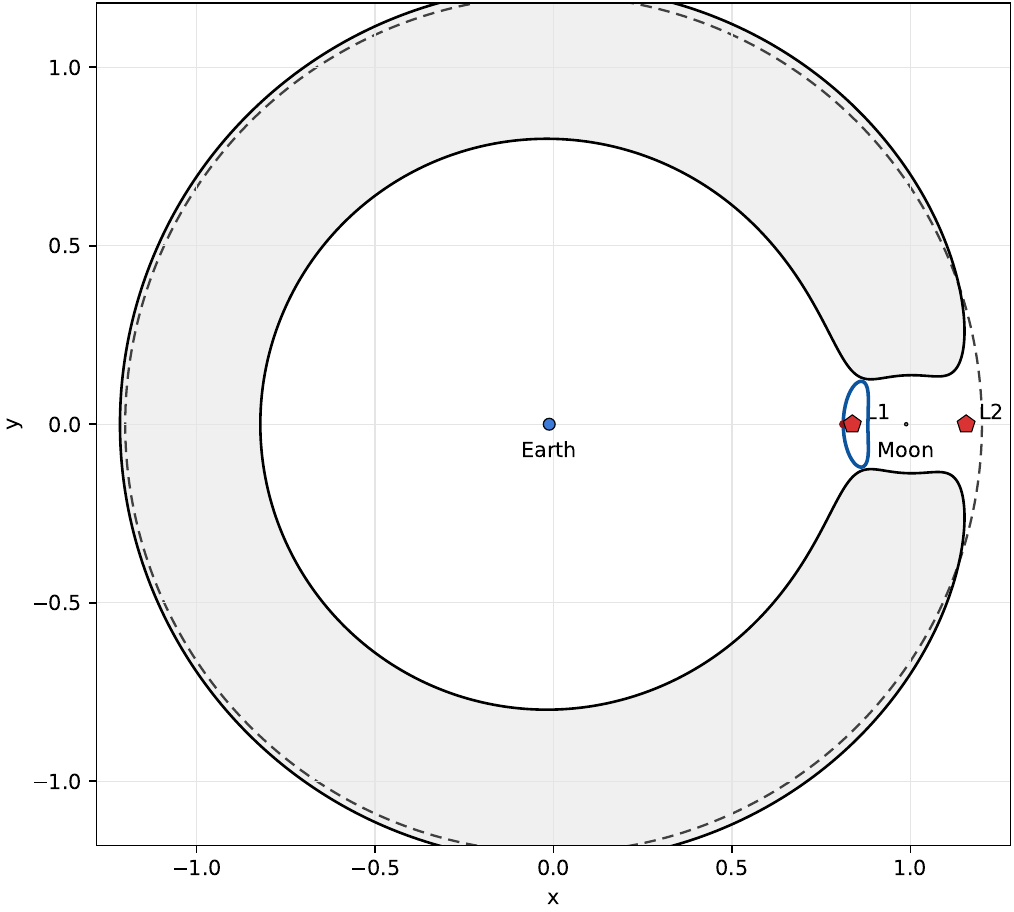}
        \caption{\footnotesize L1 Lyapunov}
    \end{subfigure}
    \hfill
    \begin{subfigure}[t]{0.23\textwidth}
        \centering
        \includegraphics[width=\textwidth]{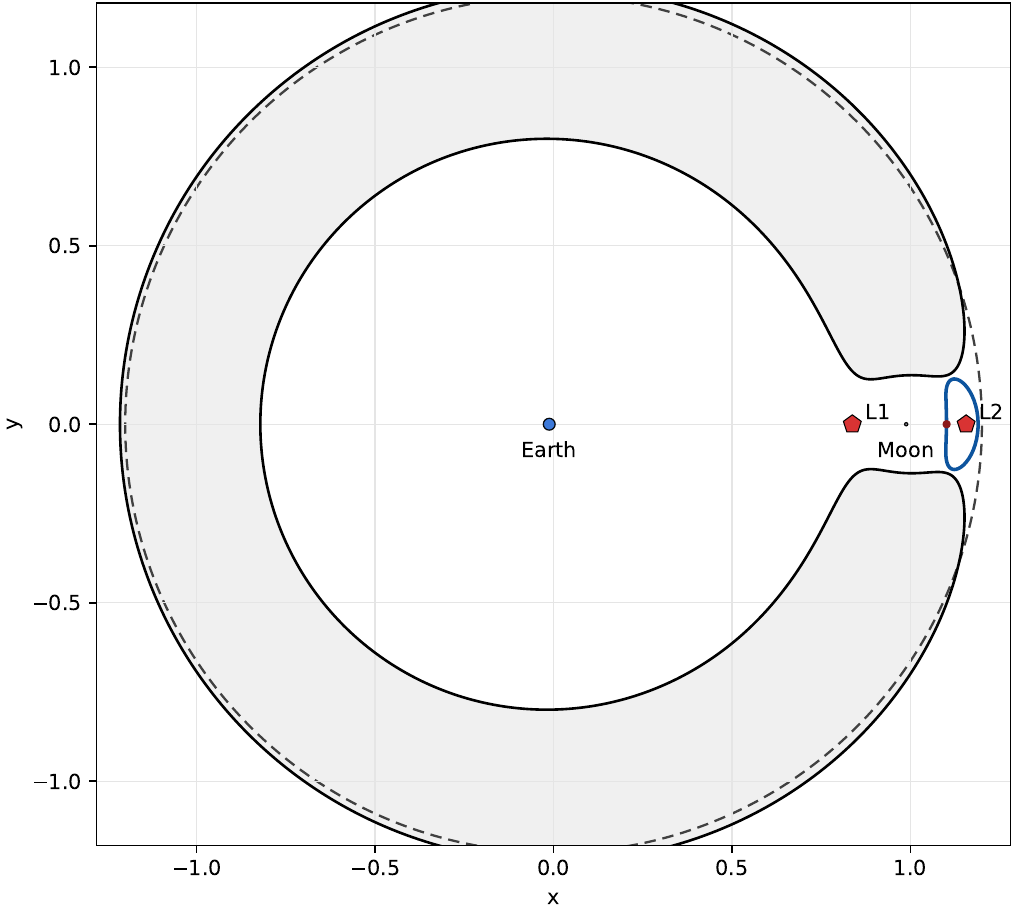}
        \caption{\footnotesize L2 Lyapunov}
    \end{subfigure}
    \hfill
    \begin{subfigure}[t]{0.23\textwidth}
        \centering
        \includegraphics[width=\textwidth]{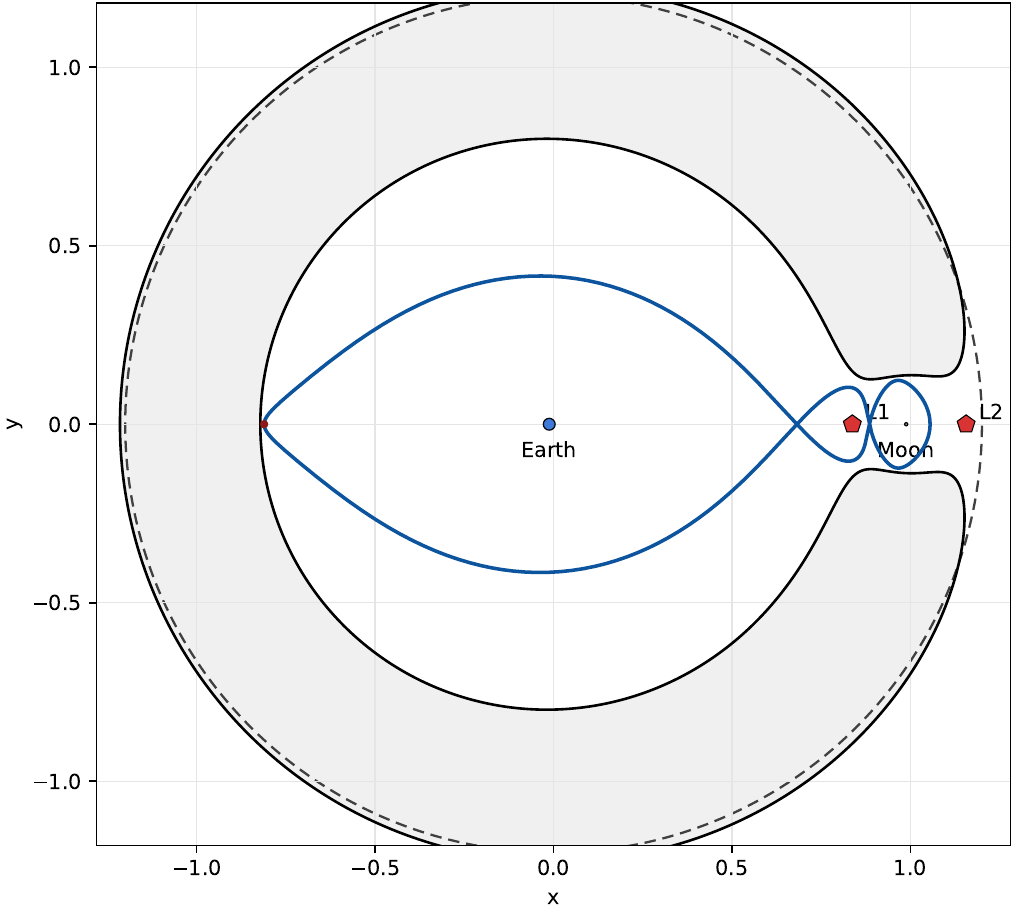}
        \caption{\footnotesize (1,1)a-cycler}
    \end{subfigure}
    \hfill
    \begin{subfigure}[t]{0.23\textwidth}
        \centering
        \includegraphics[width=\textwidth]{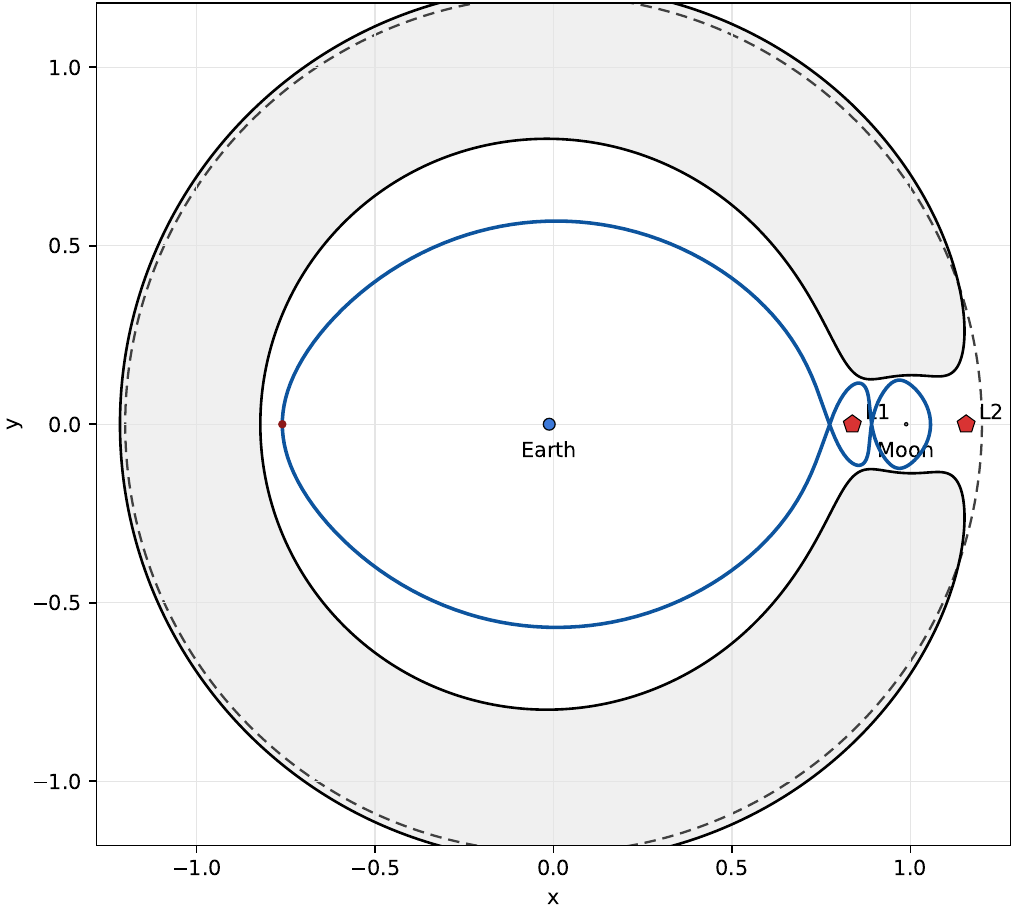}
        \caption{\footnotesize (1,1)b-cycler}
    \end{subfigure}
    
    \vspace{0.15em}

    \begin{subfigure}[t]{0.23\textwidth}
        \centering
        \includegraphics[width=\textwidth]{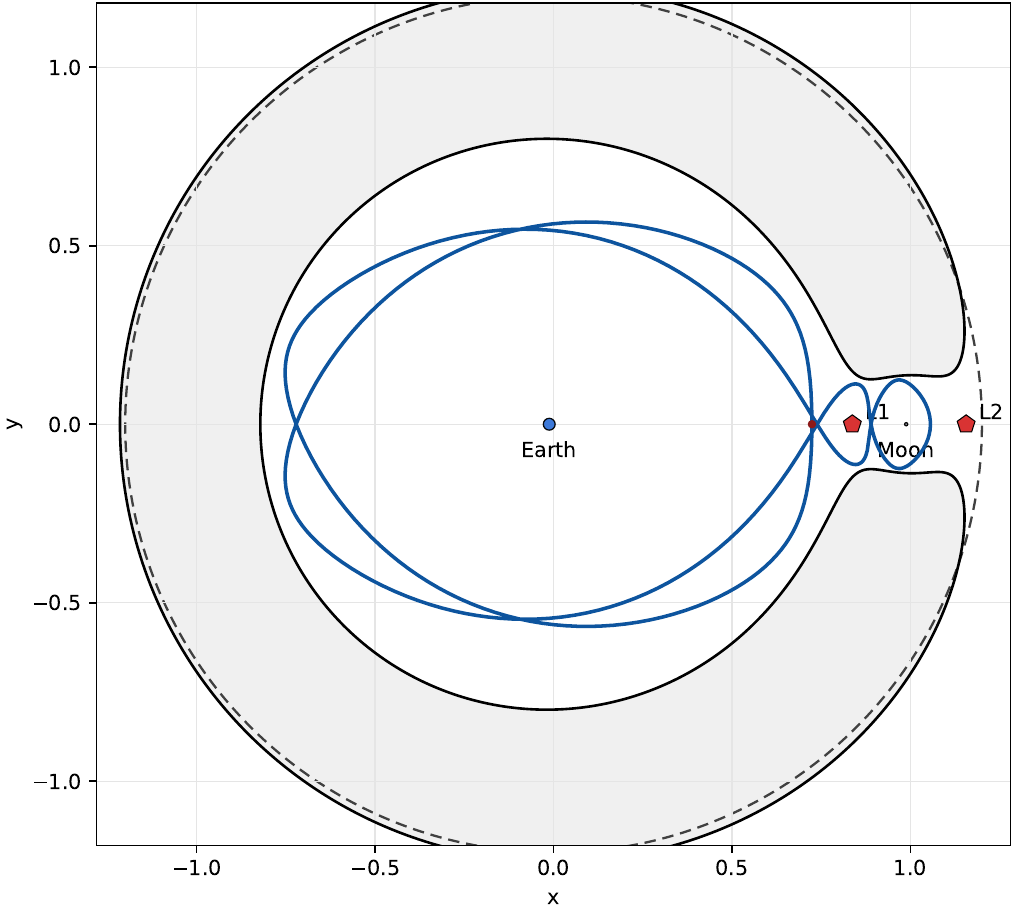}
        \caption{\footnotesize (2,1)-cycler}
    \end{subfigure}
    \hfill
    \begin{subfigure}[t]{0.23\textwidth}
        \centering
        \includegraphics[width=\textwidth]{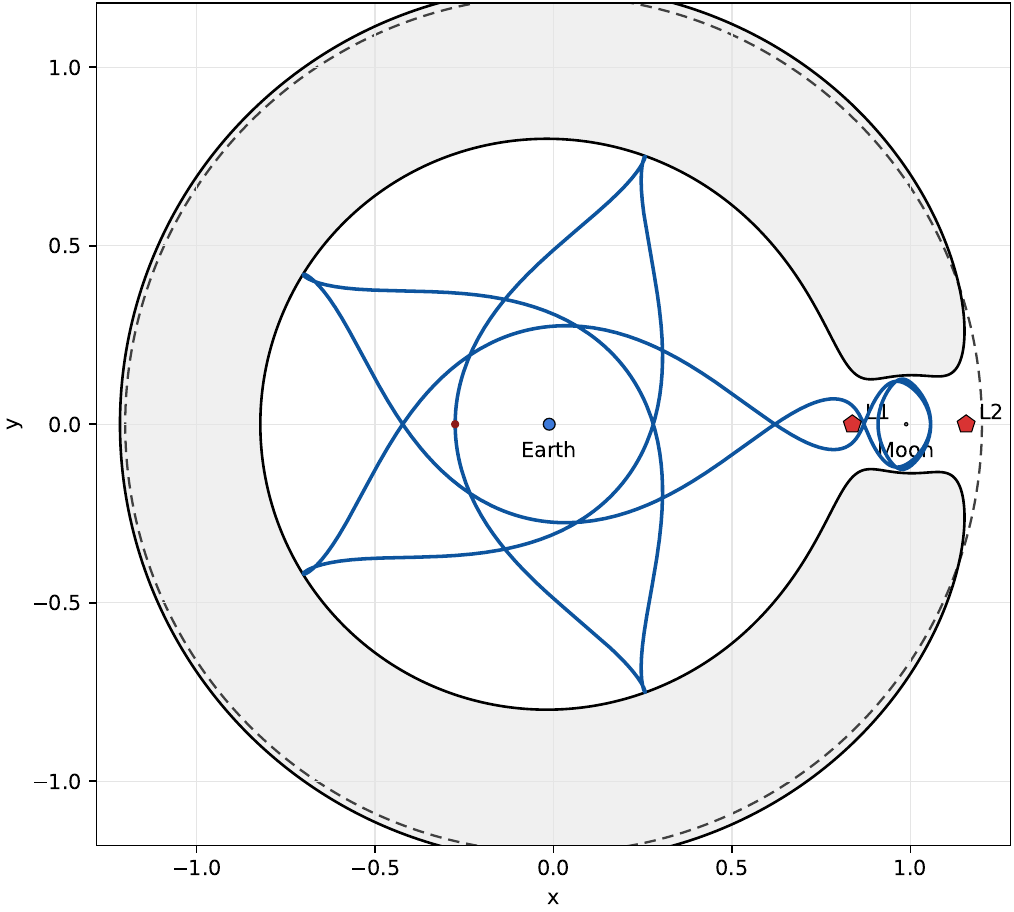}
        \caption{\footnotesize (3,2)-cycler}
    \end{subfigure}
    \hfill
    \begin{subfigure}[t]{0.23\textwidth}
        \centering
        \includegraphics[width=\textwidth]{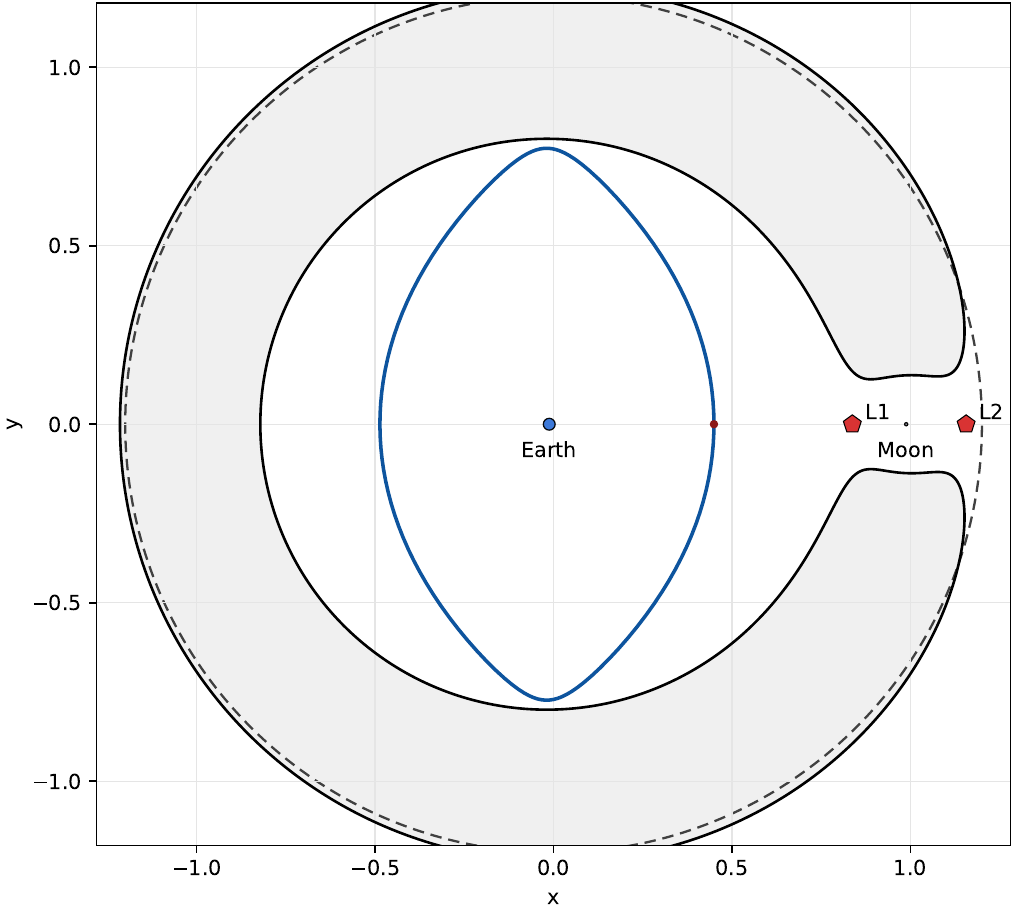}
        \caption{\footnotesize 2:1 stable resonant}
    \end{subfigure}
    \hfill
    \begin{subfigure}[t]{0.23\textwidth}
        \centering
        \includegraphics[width=\textwidth]{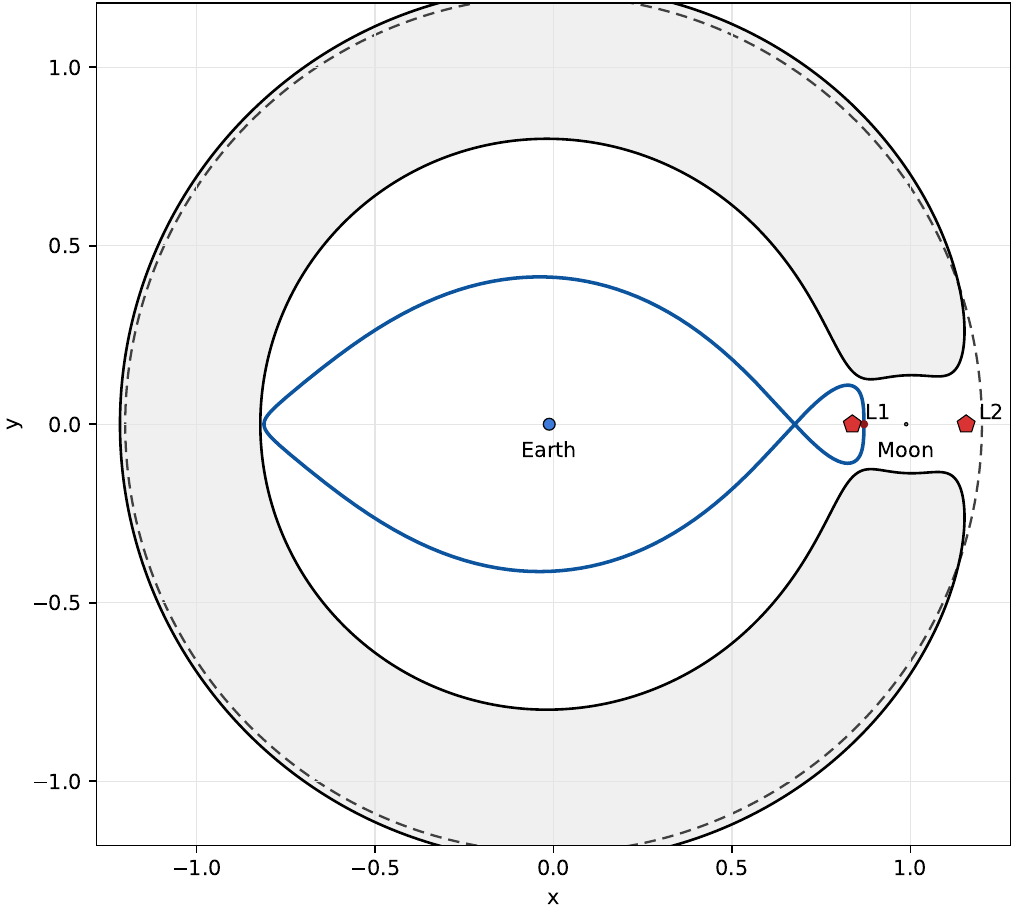}
        \caption{\footnotesize 2:1 unstable resonant}
    \end{subfigure}

    \vspace{0.15em}

    \begin{subfigure}[t]{0.23\textwidth}
        \centering
        \includegraphics[width=\textwidth]{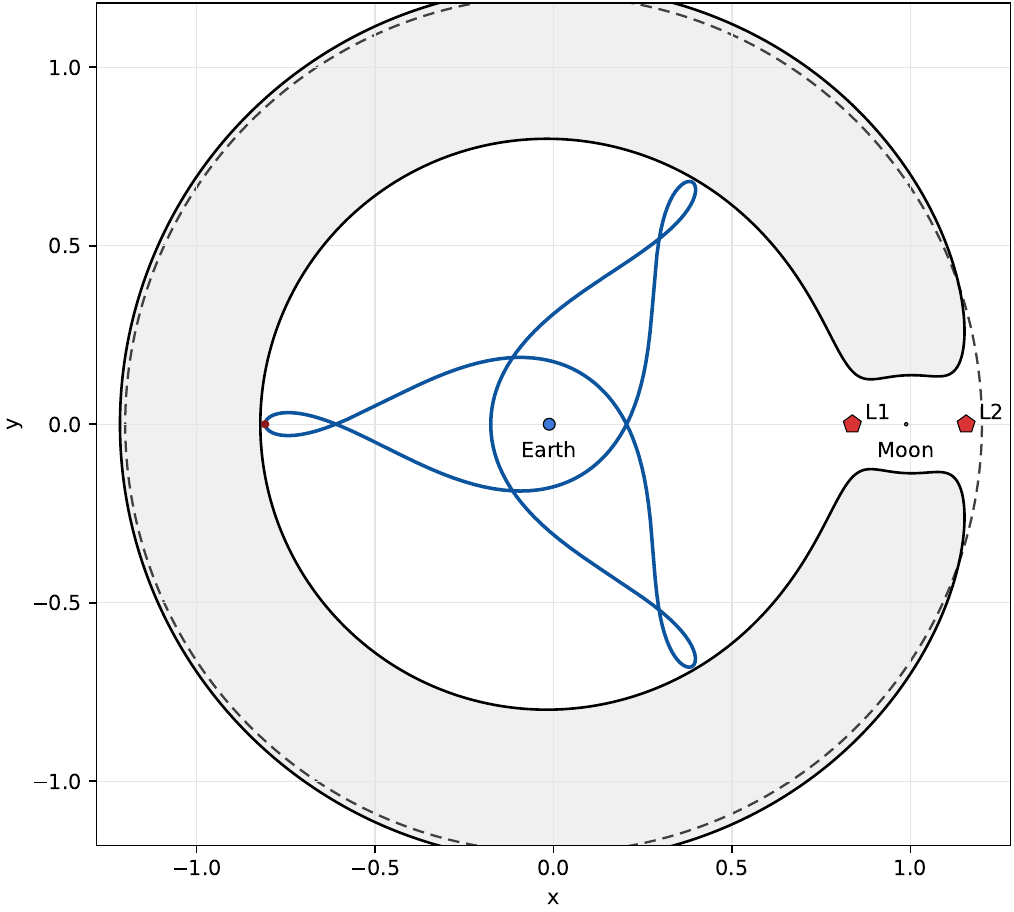}
        \caption{\footnotesize 3:1 stable resonant}
    \end{subfigure}
    \hfill
    \begin{subfigure}[t]{0.23\textwidth}
        \centering
        \includegraphics[width=\textwidth]{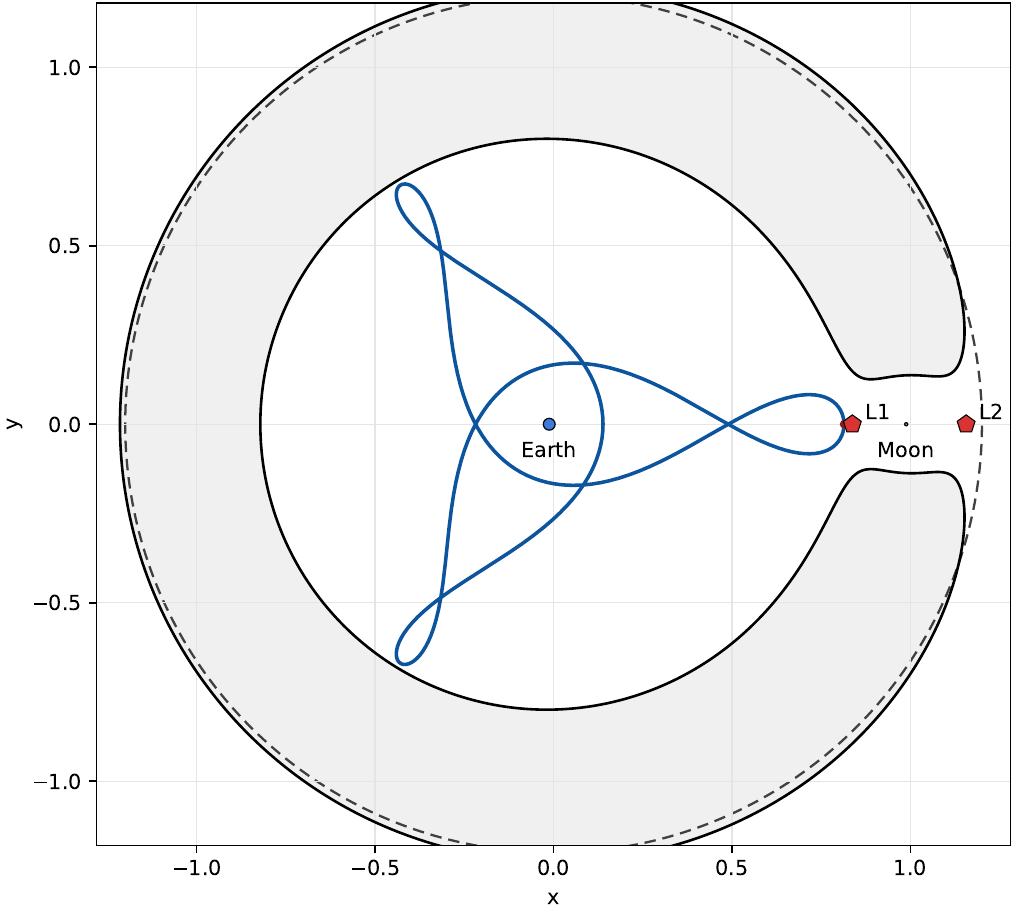}
        \caption{\footnotesize 3:1 unstable resonant}
    \end{subfigure}
    \hfill
    \begin{subfigure}[t]{0.23\textwidth}
        \centering
        \includegraphics[width=\textwidth]{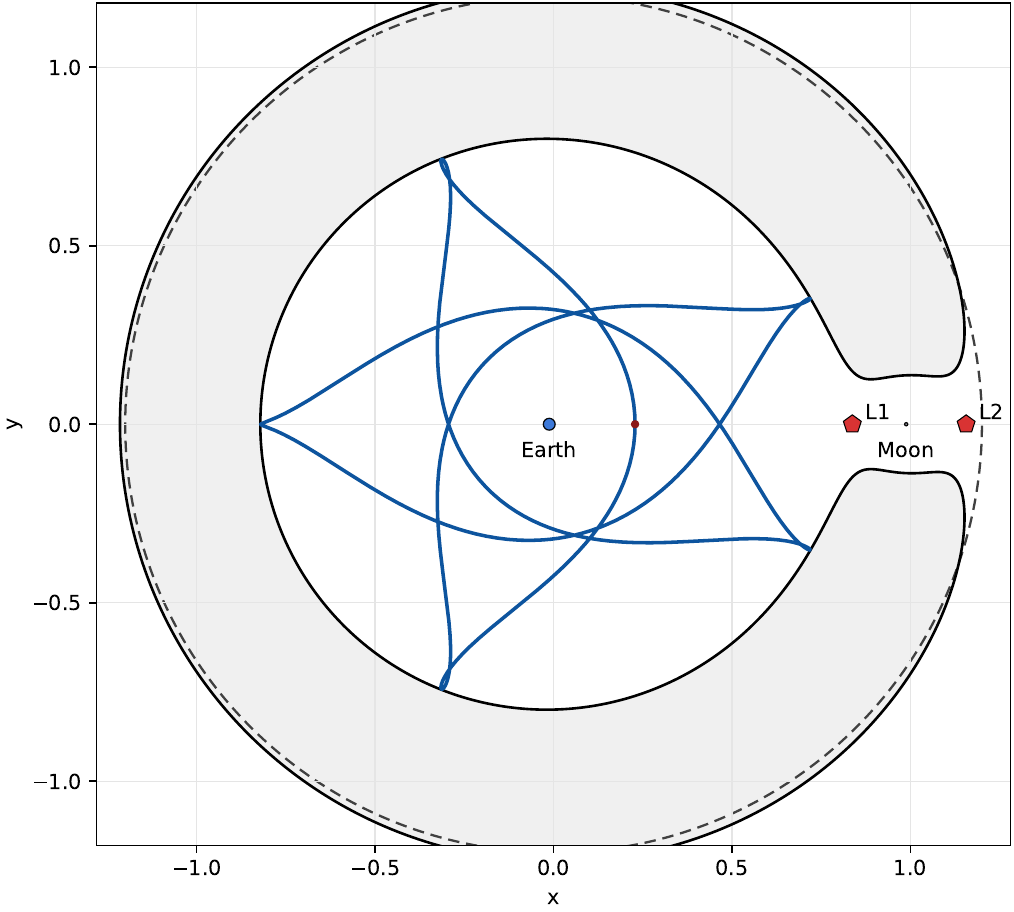}
        \caption{\footnotesize 5:2 stable resonant}
    \end{subfigure}
    \hfill
    \begin{subfigure}[t]{0.23\textwidth}
        \centering
        \includegraphics[width=\textwidth]{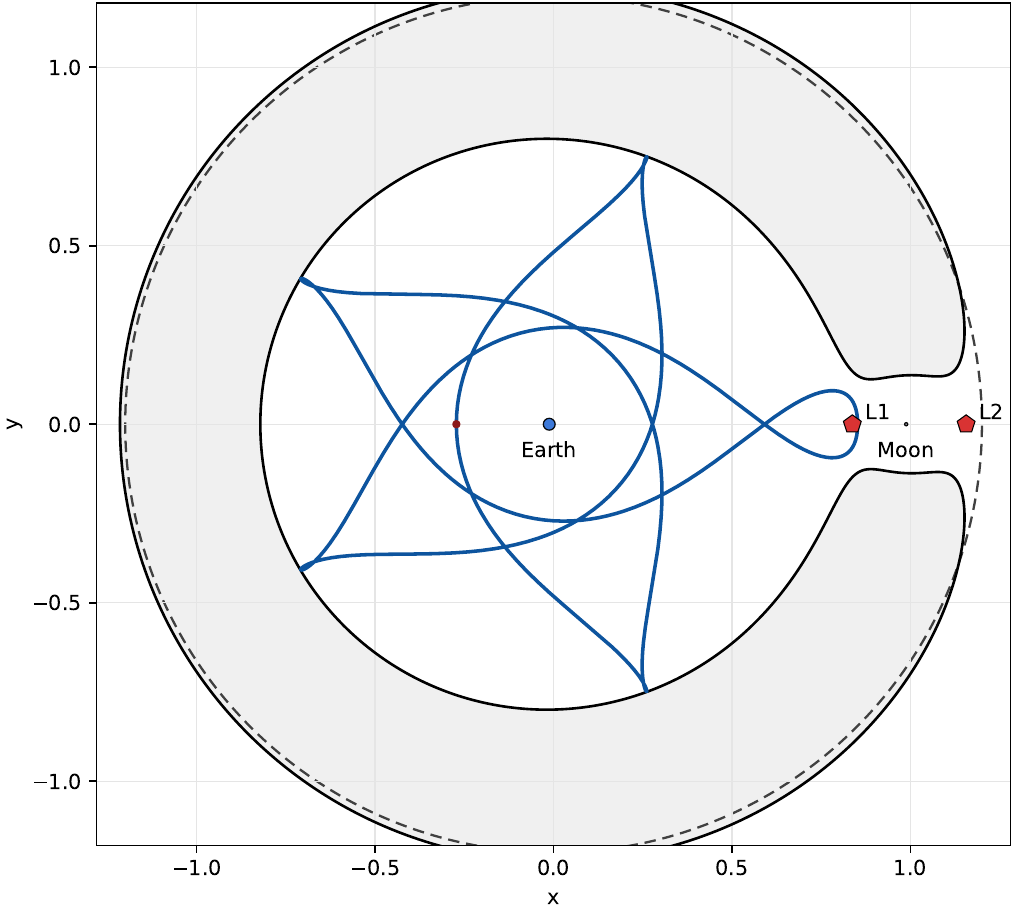}
        \caption{\footnotesize 5:2 unstable resonant}
    \end{subfigure}

    \vspace{0.15em}

    \begin{subfigure}[t]{0.23\textwidth}
        \centering
        \includegraphics[width=\textwidth]{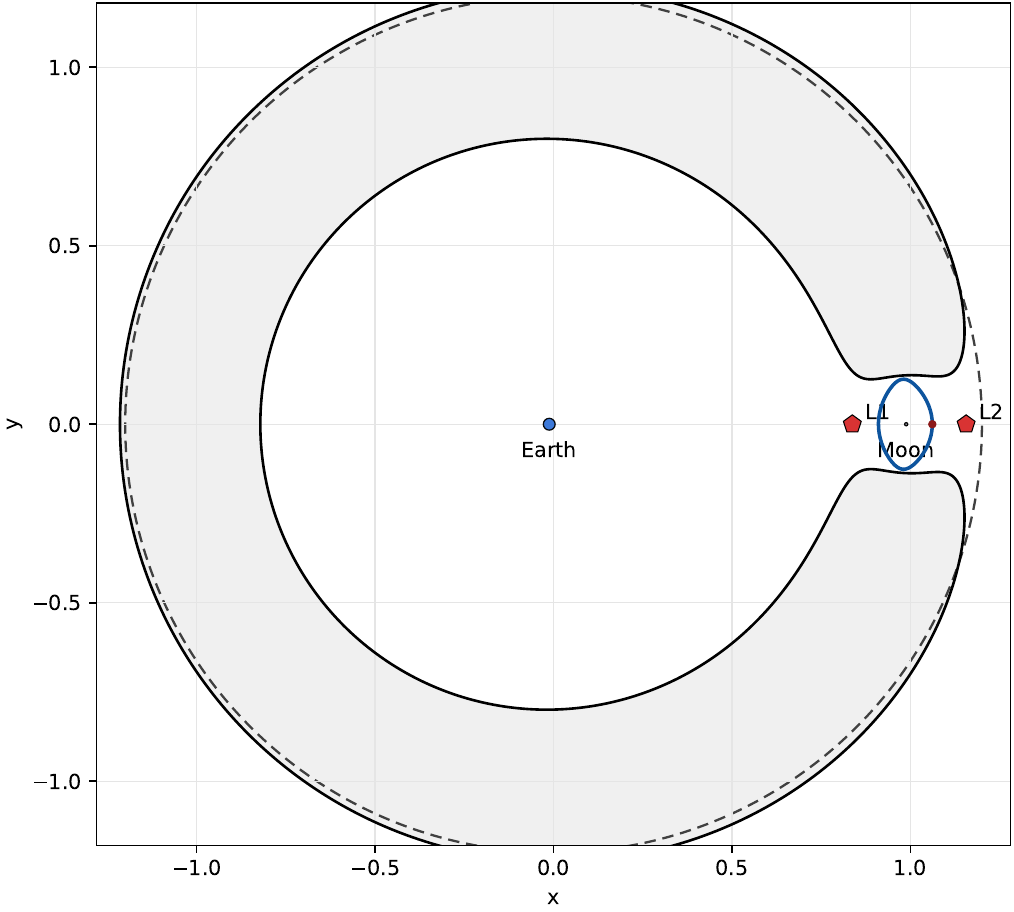}
        \caption{\footnotesize Distant prograde orbit}
    \end{subfigure}

    \caption{\footnotesize Representative periodic orbits used in the family-to-family cislunar accessibility study.
    Each panel shows one representative planar Earth--Moon CR3BP periodic orbit at the common Jacobi level \(C_J=3.1294\).}
    \label{fig:representative_orbits}
\end{center}
\medskip

Geometrically, the stability parameter encodes which region of the mixed phase space \citep{ScOt1997} of the Jacobi manifold each representative orbit occupies.
The stable resonant orbits, with $\sigma = 0$,  sit at the center of nested families of invariant tori on the manifold, while the ten unstable representatives, with $\sigma > 0$, are embedded in what we hypothesize is a single connected chaotic region of the manifold that fills the complement of those tori \citep{kumar2006}. 
Within such a connected chaotic region, trajectories starting near any unstable representative can, in principle, visit the neighborhood of any other unstable representative through the natural dynamics alone, though the natural transport time may be long \citep{koon2011}. 
Finite-$\Delta V$ heading-change  maneuvers shorten that natural transport, so the chaotic-sea geometry provides a low-cost connectivity medium for the unstable members of the family set.  
The stable resonant representatives, by contrast, are surrounded by resonant tori that resist easy access from the chaotic sea at the budgets considered here, corresponding to the high-cost edge of the accessibility hierarchy developed in Sections~\ref{sec:overlap}--\ref{sec:budget}.

The selected representatives span Lagrange-point, cycler, resonant, and distant prograde dynamics, sampling a range of geometric scales and stability characters on the common Jacobi manifold.
They provide the family set on which the reachable-set framework of Section \ref{sec:reachable} is built.

\section{Construction of Local Reachable Set Atlases}\label{sec:reachable}

For each  family representative, a local reachable-set atlas is constructed on the common Jacobi manifold. A forward reachable set is obtained by direct propagation of a heading-fan maneuver model from arc-length-spaced seeds along the representative orbit, and the corresponding backward reachable set is obtained from the time-reversal symmetry of the planar CR3BP without additional integration.
The remainder of this section develops the seed selection, the fixed-$C_J$ heading-change maneuver model, the finite-time propagation, and the voxel-level reachability logging used in this construction.


\subsection{Seed Selection Along Representative Periodic Orbits}

Each local atlas is constructed from seed states placed on the representative periodic orbit of a selected cislunar family.
For a given family, the representative orbit is first propagated over one full period in the reduced \((x,y,\theta)\) model introduced in Section~\ref{sec:framework}, yielding a dense discretization,
\begin{equation}
\mathbf{X}_{\mathrm{PO}}(t)=\bigl(x_{\mathrm{PO}}(t),\,y_{\mathrm{PO}}(t),\,\theta_{\mathrm{PO}}(t)\bigr),
\qquad
t\in[0,T_{\mathrm{PO}}].
\end{equation}
This dense trajectory is used only to parameterize the closed orbit and interpolate seed locations.

Seeds are distributed along the orbit by arc length rather than uniform time spacing.
Let \(s\) denote cumulative arc length along the trajectory and \(L_{\mathrm{PO}}\) the total orbit length.
Seeds are selected at uniform arc-length intervals around the loop.
With \(N_s\) seeds, the seed arc-length locations are,
\begin{equation}
s_j=j\,\Delta s_{\mathrm{seed}},
\qquad
\Delta s_{\mathrm{seed}}=\frac{L_{\mathrm{PO}}}{N_s},
\qquad
j=0,1,\dots,N_s-1.
\end{equation}
Arc-length sampling avoids the clustering that uniform time sampling would produce in low-speed regions and provides more even geometric coverage of the orbit.

The corresponding reduced states are obtained by interpolation along the dense periodic orbit.
Denoting the \(j\)-th seed state by,
\begin{equation}
\mathbf{X}_j=\bigl(x_j,\,y_j,\,\theta_j\bigr),
\end{equation}
the collection \(\{\mathbf{X}_j\}\) defines the initial states for the local forward reachable set.
Because the seeds lie on the representative periodic orbit, they inherit its Jacobi level by construction.
Arc-length seed placement is illustrated in Figure~\ref{fig:atlas_construction_overview}(a).

\par\medskip
\begin{center}
\captionsetup{type=figure}
    \centering

    \begin{subfigure}[t]{0.27\textwidth}
        \centering
        \includegraphics[width=0.50\linewidth]{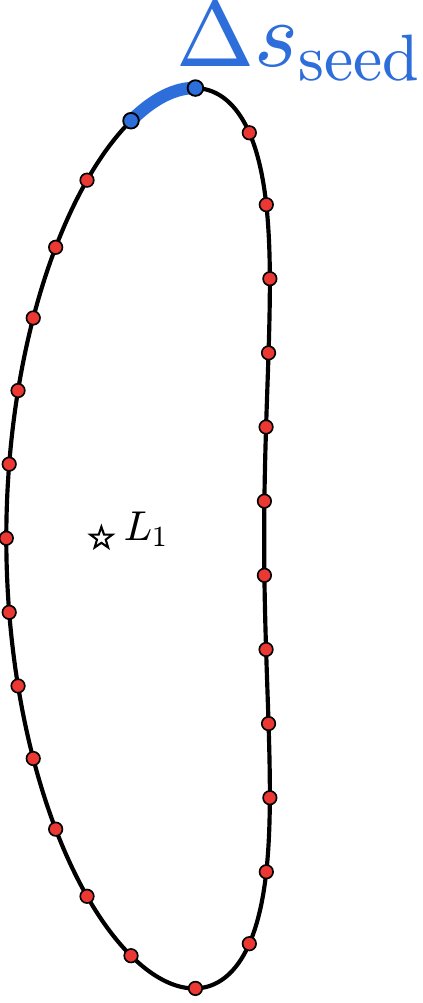}
        \caption{\footnotesize Arc-length seeds.}
    \end{subfigure}
    \hfill
    \begin{subfigure}[t]{0.29\textwidth}
        \centering
        \includegraphics[width=\textwidth]{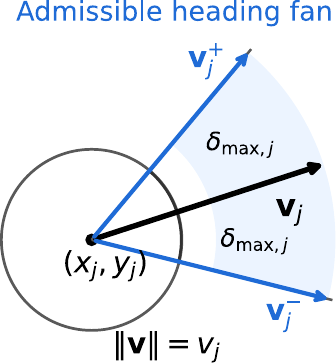}
        \caption{\footnotesize Heading-change fan.}
    \end{subfigure}
    \hfill
    \begin{subfigure}[t]{0.38\textwidth}
        \centering
        \includegraphics[width=\textwidth]{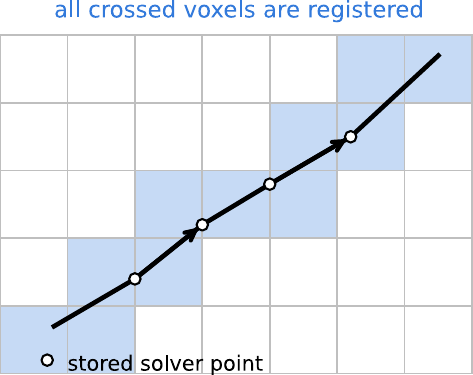}
        \caption{\footnotesize Segment-walk registration.}
    \end{subfigure}

    \caption{\footnotesize Overview of the local reachable-set atlas construction.
    Panel (a) shows arc-length seed placement along a representative periodic orbit.
    Panel (b) shows the fixed-$C_J$ heading-change maneuver fan applied at each seed.
    Panel (c) shows a schematic illustration of segment-walk voxel registration,
    which records crossed voxels between successive propagated states.
    \revone{The voxel grid shown in panel (c) is illustrative and is not drawn at the
    numerical resolution used in the analysis; the nominal voxel dimensions used in the computations are reported in Table~\ref{tab:atlas_numerics}.}}
    \label{fig:atlas_construction_overview}
\end{center}
\medskip


\subsection{Constant-\(C_J\) Heading Change Maneuver Model}

Local accessibility is probed by applying an instantaneous in-plane heading change at each seed location along an orbit, while holding position and speed fixed.
Because the speed is preserved, the maneuver also preserves \(C_J\), so the perturbed state remains on the same fixed-energy manifold.

At the seed \(\mathbf{X}_j\), the local speed is,
\begin{equation}
v_j=\sqrt{-2\bar U(x_j,y_j)-C_J}.
\end{equation}
A heading change is introduced by an angular offset \(\delta\), so that the post-maneuver heading becomes,
\begin{equation}
\theta_j^{+}=\theta_j+\delta.
\label{eq:theta_post_maneuver}
\end{equation}
The required impulse is the vector difference between two equal-magnitude velocity vectors separated by the turning angle \(\delta\):
\begin{equation}
\Delta V_{\mathrm{turn}} = 2v_j\sin\!\left(\frac{|\delta|}{2}\right).
\label{eq:dv_turn}
\end{equation}

The admissible heading changes are constrained by a one-sided maneuver cap \(\Delta V_{\mathrm{a}}\), so that \(\Delta V_{\mathrm{turn}} \le \Delta V_{\mathrm{a}}\).
Here, ``one-sided'' means that trajectories are propagated away from a single family.
For pairwise accessibility under total budgets \((\Delta V_{\mathrm{cap}},T_{\mathrm{cap}})\), the construction uses the symmetric split,
\begin{equation}
\Delta V_{\mathrm{a}} = \frac{\Delta V_{\mathrm{cap}}}{2},
\qquad
T_{\mathrm{a}} = \frac{T_{\mathrm{cap}}}{2},
\end{equation}
so that each family contributes one local leg under half of the total budget.
Hence,
\begin{equation}
|\delta|\le \delta_{\max,j},
\qquad
\delta_{\max,j}
=
2\arcsin\!\left(\frac{\Delta V_{\mathrm{a}}}{2v_j}\right),
\label{eq:deltamax}
\end{equation}
when \(0\le \Delta V_{\mathrm{a}} \le 2v_j\).
If \(\Delta V_{\mathrm{a}}>2v_j\), then the full angular range becomes admissible from that seed.

For each seed, the admissible interval \(-\delta_{\max,j}\le \delta \le \delta_{\max,j}\) is discretized into a symmetric heading fan about the nominal orbit direction using increment \(\Delta\delta_{\mathrm{fan}}\):
\begin{equation}
\delta_k = k\,\Delta\delta_{\mathrm{fan}},
\qquad
|\delta_k|\le \delta_{\max,j}.
\end{equation}
The fan increment \(\Delta\delta_{\mathrm{fan}}\) is distinct from the heading voxel spacing \(\Delta\theta\) used later to discretize the reduced phase space.
Each seed therefore launches a finite set of directionally perturbed trajectories on the same fixed-\(C_J\) manifold.
The local heading fan is illustrated in Figure~\ref{fig:atlas_construction_overview}(b).

This maneuver model is a deliberate control abstraction. It does not represent general impulsive transfer design, but isolates the effect of bounded velocity-direction changes on a fixed-\(C_J\) manifold:
by construction, every admissible maneuver is a pure rotation of the rotating-frame velocity vector at fixed rotating-frame speed, with no change in rotating-frame speed magnitude. This is a strong restriction relative to general impulsive maneuvers, and in many configurations it corresponds to a substantial normal component of the inertial impulse rather than the tangential burns that two-body intuition associates with efficient orbit transfer \citep{BaMuWh1971}. 
The results that follow show that this restricted maneuver model is nonetheless sufficient to achieve low-$\Delta V$ family-to-family accessibility on the chaotic region of the Jacobi manifold; the chaotic-region connectivity provides the geometric medium that purely directional, fixed-rotating-speed control can exploit. 
The angular sampling increment \(\Delta\delta_{\mathrm{fan}}\) is specified later with the discretization choices.


\subsection{Finite-Time Propagation and Forward Reachable Sets}

For each seed \(\mathbf{X}_j\) and each admissible heading offset \(\delta\), the post-maneuver state is propagated forward over a one-sided propagation horizon \(T_{\mathrm{a}}\).
The resulting trajectory is written as,
\begin{equation}
\mathbf{X}(t;\mathbf{X}_j,\delta)
=
\bigl(x(t;\mathbf{X}_j,\delta),\,y(t;\mathbf{X}_j,\delta),\,\theta(t;\mathbf{X}_j,\delta)\bigr),
\qquad
0\le t\le T_{\mathrm{a}},
\end{equation}
with initial condition \(\mathbf{X}(0;\mathbf{X}_j,\delta)=\bigl(x_j,\,y_j,\,\theta_j+\delta\bigr)\).
Propagation terminates when time \(T_{\mathrm{a}}\) is reached, when the trajectory impacts either primary, or when it exits \(\mathcal{D}(C_J)\).
During reduced-state integration, \(v^2=-2\bar U-C_J\) is monitored to avoid the singular heading dynamics near the zero-velocity curve.
If \(v^2\) falls below a prescribed threshold, propagation is continued in the full planar state model.

The propagated trajectories are recorded on a reduced phase-space voxel grid, consistent with the set-oriented view of transport as set propagation on a partitioned phase space \citep{dellnitz2005}.
Define the energy-reduced phase space as,
\begin{equation}    
\mathcal{Q}(C_J)=\mathcal{D}(C_J)\times S^1,
\end{equation}
where \(\mathcal{D}(C_J)\) is the admissible cislunar domain and \(S^1\) is the circle corresponding to the heading coordinate.
This space is partitioned into voxels (three-dimensional boxes),
\begin{equation}
\mathcal{Z}
=
\left\{
Z_{\alpha\beta\gamma}
=
I_\alpha \times J_\beta \times K_\gamma
\right\},
\end{equation}
where \(I_\alpha\), \(J_\beta\), and \(K_\gamma\) are the intervals associated with the \(x\)-, \(y\)-, and \(\theta\)-discretizations, respectively.

For a given representative orbit, the forward reachable set is the collection of voxels visited by at least one admissible trajectory launched from any seed.
\revone{Denote this family-level set by
\(\mathcal{F}_{\mathrm{fam}}(\Delta V_a,T_a)\), where,}
\begin{equation}
\mathcal{F}_{\mathrm{fam}}(\Delta V_{\mathrm{a}},T_{\mathrm{a}})
=
\left\{
Z\in\mathcal{Z}
\;:\;
\exists\, j,\;
\exists\,\delta \text{ with } |\delta|\le \delta_{\max,j},\;
\exists\, t\in[0,T_{\mathrm{a}}]
\text{ such that }
\mathbf{X}(t;\mathbf{X}_j,\delta)\in Z
\right\}.
\label{eq:frs_family}
\end{equation}
Computing this set from stored propagation output alone is unreliable.
With adaptive integration, successive stored states may cross intermediate voxels without placing an output point inside them, and stored states are in general not aligned with voxel boundaries.
To register every crossed voxel, the trajectory between successive propagated states is walked through the voxel grid, interpolated  position and heading along each segment. 
Segment-walk registration is illustrated in Figure~\ref{fig:atlas_construction_overview}(c).

The nominal discretization and numerical parameters are summarized in Table~\ref{tab:atlas_numerics}.
These settings define the voxel resolution, seed and heading-fan sampling, one-sided budgets, and propagation tolerances used for the primary results.
Sensitivity to alternative numerical choices is assessed separately in the validation study.

\par\medskip
\begin{center}
\captionsetup{type=table}
\centering
\caption{\footnotesize Nominal discretization and numerical parameters used in the reachable-set construction.}
\label{tab:atlas_numerics}
\begin{tabular}{lll}
\toprule
Quantity & Symbol & Value \\
\midrule
Finite cislunar domain radius & $R_{\mathrm{dom}}$ & 1.2 \\
Spatial voxel spacing & $\Delta x=\Delta y$ & 0.001 \; (\SI{384.4}{km}) \\
Heading voxel spacing & $\Delta\theta$ & $1^\circ$ \\
Seed arc-length spacing & $\Delta s_{\mathrm{seed}}$ & 0.01 \\
Heading-fan sampling increment & $\Delta\delta_{\mathrm{fan}}$ & $0.5^\circ$ \\
One-sided maneuver cap & $\Delta V_{\mathrm{a}}$ & 0.2 \; (\SI{204.6}{m.s^{-1}}) \\
One-sided propagation horizon & $T_{\mathrm{a}}$ & $\pi$ \; (\SI{13.66}{day}) \\
Relative tolerance & $\mathrm{RelTol}$ & $10^{-8}$ \\
Absolute tolerance & $\mathrm{AbsTol}$ & $10^{-8}$ \\
Reduced-model fallback threshold & $v^2_{\mathrm{tol}}$ & $10^{-8}$ \\
\bottomrule
\end{tabular}
\end{center}
\medskip

The discretization in Table~\ref{tab:atlas_numerics} is chosen so that the reachable-set construction is limited by the underlying dynamics rather than by the numerical grid, while remaining computationally tractable. The spatial voxel spacing $\Delta x = \Delta y = 0.001$
is fine enough that two reachable sets do not register as overlapping at the voxel scale unless their continuous images actually approach within roughly 384 km ($=0.001~LU$) on the common Jacobi manifold. The heading voxel spacing $\Delta\theta = 1^\circ$
and the maneuver-fan increment $\Delta\delta_{\mathrm{fan}} = 0.5^\circ$
resolve the heading direction finely and oversample the maneuver fan by a factor of two, so that the overlap proxy depends on the dynamics rather than on the angular bin structure. The seed arc-length spacing $\Delta s_{\mathrm{seed}} = 0.01$
distributes seeds uniformly along each representative orbit, giving even geometric coverage independent of the local speed. Appendix~\ref{app:baseline_validation} reports a sensitivity study confirming that the dominant accessibility trends are preserved under coarser and finer discretizations.

\revone{
For reference, a representative computational benchmark was performed on
Virginia Tech's Owl cluster using 45 MATLAB parallel workers. The dominant
forward-propagation stage required approximately 23--230~s per family,
depending on the number of seed--heading trajectories, corresponding to
approximately 21~min for the complete thirteen-family atlas set. Once the
atlases were available, the complete set of 78 pairwise voxel-overlap searches
required approximately 104~s. The cached atlas set occupied approximately
5.2~GB, whereas the complete pairwise-overlap output occupied approximately
196~kB. These values are representative and depend on the selected spatial,
heading, seed, maneuver, and propagation resolutions.
}

\subsection{Forward and Backward Reachable Sets via Time-Reversal Symmetry}

The time-reversal symmetry introduced in Section~\ref{subsec:time-reversal} is used to obtain backward-reachable information from the forward reachable set, avoiding separate backward integrations for each seed and admissible heading.

Under the reduced symmetry map \(\mathcal{R}\), a forward trajectory \(\mathbf{X}(t;\mathbf{X}_j,\delta)\) generates the reversed trajectory \(\mathcal{R}\mathbf{X}(-t;\mathbf{X}_j,\delta)\) on the same Jacobi manifold.
The heading offset transforms consistently: since \(\theta_j^{+}=\theta_j+\delta\), the mirrored heading is \(\pi-(\theta_j+\delta)=(\pi-\theta_j)-\delta\).
Thus, the reversed control offset is \(-\delta\), which remains admissible because the heading fan is symmetric about the nominal direction.

The forward reachable set is computed by propagation and voxel logging, and the backward reachable set is obtained by applying the symmetry map, \eqref{eq:reduced_reverser}.
For a family representative,
\begin{equation}
\mathcal{B}_{\mathrm{fam}}(\Delta V_{\mathrm{a}},T_{\mathrm{a}})
=
\mathcal{R}\!\left(\mathcal{F}_{\mathrm{fam}}(\Delta V_{\mathrm{a}},T_{\mathrm{a}})\right).
\label{eq:brs_from_frs}
\end{equation}
Thus, occupied forward voxels are mapped by \(y\mapsto -y\) and \(\theta\mapsto\pi-\theta\).

Together, \(\mathcal{F}_{\mathrm{fam}}\) and \(\mathcal{B}_{\mathrm{fam}}\) form the local reachable-set atlas for each  family representative, with dependence on \(\Delta V_{\mathrm{a}}\) and \(T_{\mathrm{a}}\) suppressed when unambiguous.
Figure~\ref{fig:frs_brs_symmetry} shows an example for the  \(2\!:\!1\) unstable (R21-U) orbit family representative.

\par\medskip
\begin{center}
\captionsetup{type=figure}
    \centering
    \begin{subfigure}[t]{0.47\textwidth}
        \centering
        \includegraphics[width=\textwidth]{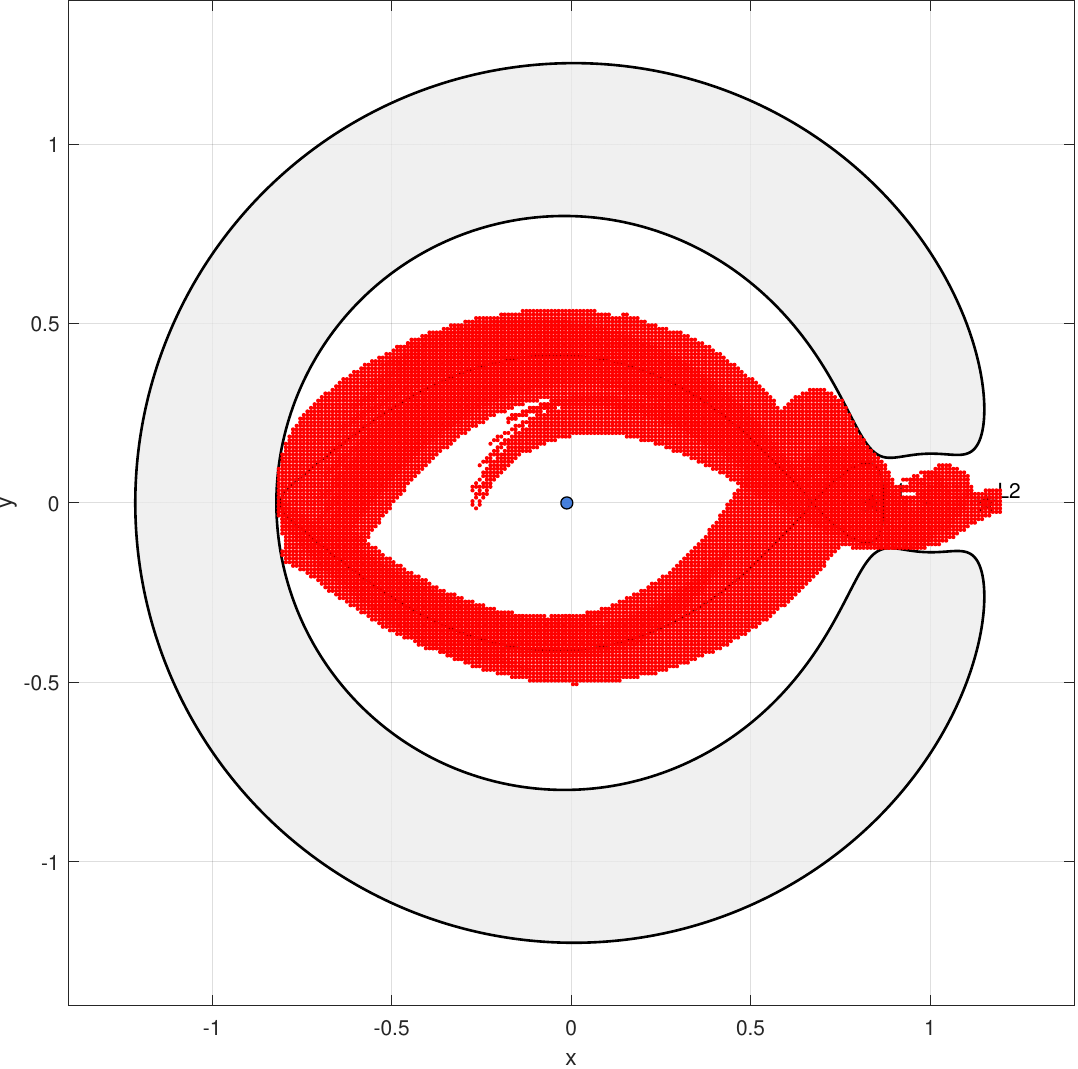}
        \caption{Forward reachable set \(\mathcal{F}_{\mathrm{fam}}\).}
    \end{subfigure}
    \hfill
    \begin{subfigure}[t]{0.47\textwidth}
        \centering
        \includegraphics[width=\textwidth]{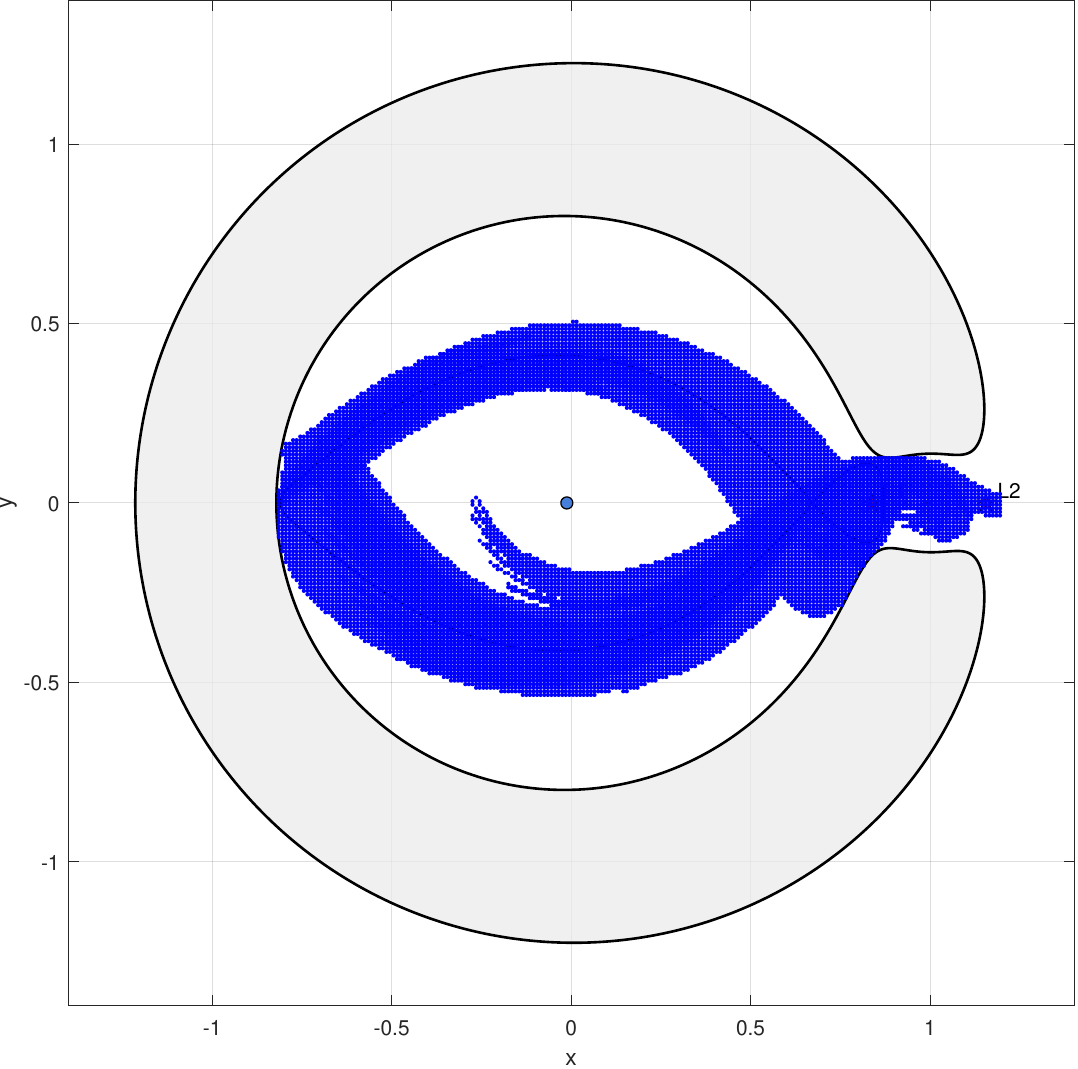}
        \caption{Backward reachable set \(\mathcal{B}_{\mathrm{fam}}=\mathcal{R}(\mathcal{F}_{\mathrm{fam}})\).}
    \end{subfigure}
    \caption{\footnotesize Example of the symmetry-based construction of forward and backward reachable sets for the  \(2\!:\!1\) unstable resonant  family representative.
    Panel (a) shows the forward reachable set obtained directly by propagation and voxel logging.
    Panel (b) shows the backward reachable set obtained by applying \(\mathcal{R}(x,y,\theta)=(x,-y,\pi-\theta)\) to the forward reachable set.}
    \label{fig:frs_brs_symmetry}
\end{center}
\medskip

\section{Inferring Pairwise Family Accessibility from Reachable Set Overlap}\label{sec:overlap}

The local reachable-set atlases of Section~\ref{sec:reachable} describe how individual families propagate forward and backward on the common Jacobi manifold, but they do not yet specify when two families should be considered accessible to each other. This section develops that link. Overlap between forward and backward reachable sets is used to identify reduced phase-space regions in which departing motion from one family is compatible with approaching motion toward another, and voxel-level proxy metrics on those regions reduce each family pair to a single accessibility cost and time.


\subsection{Pairwise Overlap Between Forward and Backward Reachable Sets}

Two families are considered locally accessible if a trajectory departing from one and a trajectory approaching the other can meet in the same region of reduced phase space within the prescribed maneuver and time budgets. This idea is made precise by intersecting the forward reachable set of the departure family with the backward reachable set of the arrival family on the common voxel grid \(\mathcal{Q}(C_J)=\mathcal{D}(C_J)\times S^1\) at the shared Jacobi level \(C_J=3.1294\).

The pairwise accessibility problem is posed at total maneuver and time budgets \((\Delta V_{\mathrm{cap}},T_{\mathrm{cap}})\).
In the present construction, the local atlases are used symmetrically, so the one-sided atlas parameters satisfy,
\begin{equation}
\Delta V_{\mathrm{a}}=\frac{\Delta V_{\mathrm{cap}}}{2},
\qquad
T_{\mathrm{a}}=\frac{T_{\mathrm{cap}}}{2}.
\end{equation}

Because all family atlases are constructed on the same voxelization of \(\mathcal{Q}(C_J)\), pairwise overlap is well defined at the voxel level. Let \(\mathcal{F}_A(\Delta V_{\mathrm{a}},T_{\mathrm{a}})\) denote the forward reachable set of family \(A\), and let \(\mathcal{B}_B(\Delta V_{\mathrm{a}},T_{\mathrm{a}})\) denote the backward reachable set of family \(B\). The pairwise overlap is then defined by,
\begin{equation}
\mathcal{O}_{A\to B}(\Delta V_{\mathrm{cap}},T_{\mathrm{cap}})
=
\mathcal{F}_A\!\left(
\frac{\Delta V_{\mathrm{cap}}}{2},
\frac{T_{\mathrm{cap}}}{2}
\right)
\cap
\mathcal{B}_B\!\left(
\frac{\Delta V_{\mathrm{cap}}}{2},
\frac{T_{\mathrm{cap}}}{2}
\right).
\label{eq:pairwise_overlap}
\end{equation}
It should be interpreted as a necessary, but not sufficient, condition for a direct family-to-family connection in the present framework: if
\begin{equation}
\mathcal{O}_{A\to B}(\Delta V_{\mathrm{cap}},T_{\mathrm{cap}})=\emptyset,
\label{eq:empty_overlap}
\end{equation}
then no direct connection from family \(A\) to family \(B\) is inferred within the present reachable set framework.

An example of this construction is shown in Figure~\ref{fig:pairwise_overlap_example} for the overlap between the forward reachable set of L1 Lyapunov (LL1) and the backward reachable set of the R21-U family.

\par\medskip
\begin{center}
\captionsetup{type=figure}
    \centering
    \includegraphics[width=0.62\textwidth]{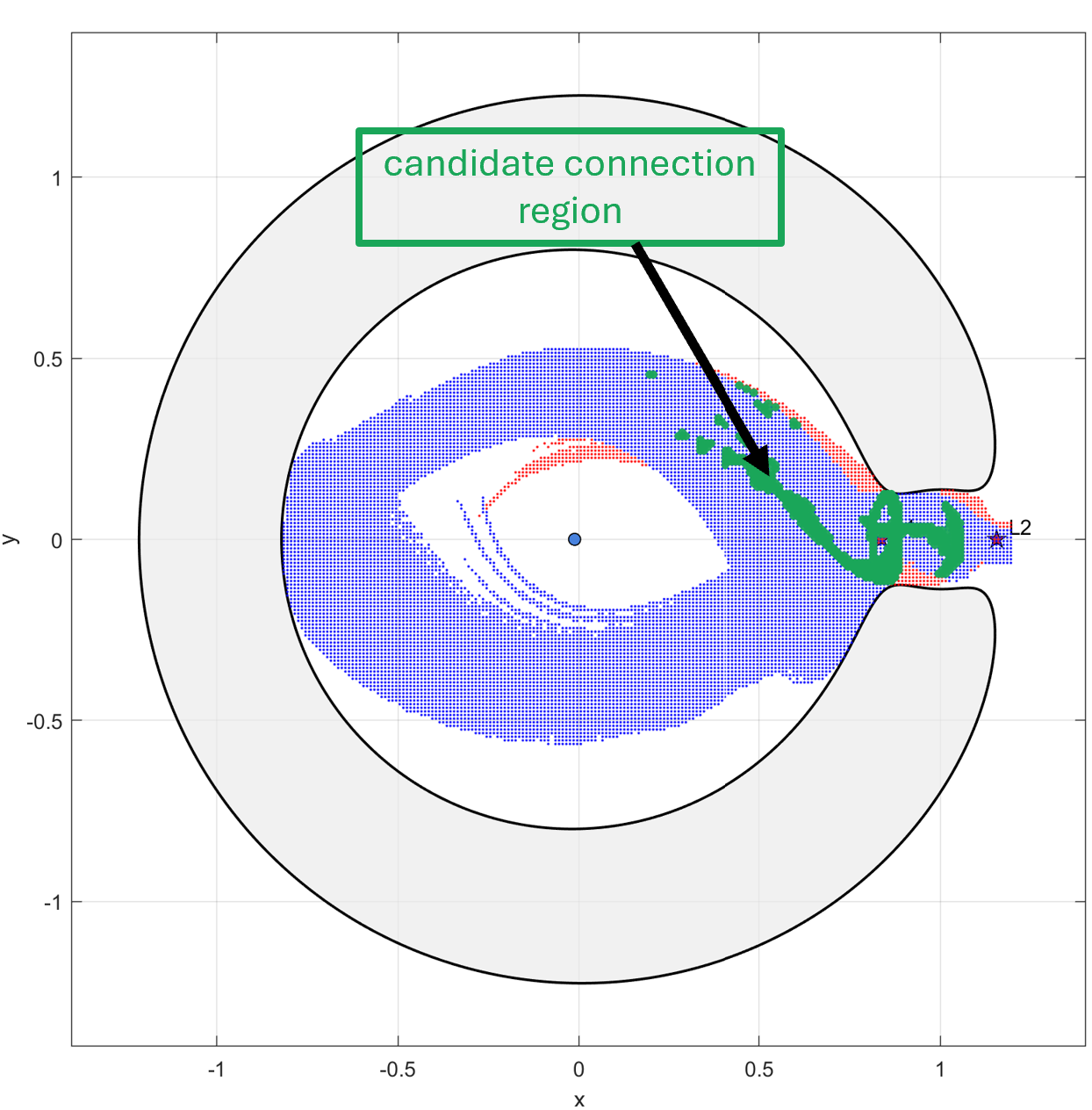}
    \caption{\footnotesize Example of pairwise accessibility from reachable-set overlap for L1 Lyapunov to \(2\!:\!1\) unstable resonant.
    The annotated green co-occupied voxels indicate candidate connection regions, which are ranked in Section~\ref{subsec:ranking} using proxy cost and time metrics.
    All reachable sets are shown as projections onto the \((x,y)\) plane.}
    \label{fig:pairwise_overlap_example}
\end{center}
\medskip


\subsection{Overlap Voxels as Candidate Local Connection Regions}

The pairwise overlap defined by \eqref{eq:pairwise_overlap} is not only a binary indicator of accessibility.
It also localizes reduced phase-space regions where motion departing from family \(A\) is compatible with motion approaching family \(B\).

If \(Z\in\mathcal{O}_{A\rightarrow B}\), then by construction that voxel is reached by departing motion from family \(A\) and is also compatible with approaching motion toward family \(B\).
Thus, the voxel lies in a local region of the common Jacobi energy manifold where the two one-sided atlases are simultaneously compatible under the prescribed pairwise budgets.
Each overlap voxel is therefore interpreted as a candidate local transfer region between the two families.

Accordingly, the overlap set may be resolved into a finite collection of candidate voxels,

\begin{equation}
\mathcal{O}_{A\rightarrow B}
=
\left\{
Z_k
\right\}_{k=1}^{N_{\mathrm{ov}}},
\end{equation}
where \(N_{\mathrm{ov}}\) denotes the number of overlap voxels.
Each voxel \(Z_k\) is identified by a local position and heading region \((x,y,\theta)\) compatible with both families.

The overlap should be interpreted geometrically rather than as a transfer solution.
An overlap voxel does not provide a corrected trajectory, a unique maneuver history, or an optimized transfer cost; it identifies a local compatibility region that can be ranked through voxel-level proxy metrics.

Although the number and distribution of overlap voxels describe where compatibility occurs, the accessibility matrix is not based on raw overlap volume.
Instead, each pair \((A,B)\) is ranked using voxel-level proxy metrics and then reduced to a single pairwise value.
Those metrics are introduced next.


\subsection{Voxel-Level Proxy Connection Metrics and Pairwise Aggregation}\label{subsec:ranking}

Each overlap voxel \(Z_k\in\mathcal{O}_{A\rightarrow B}\) is assigned a proxy maneuver cost and a proxy transfer time of flight that rank the candidate local connection regions. These quantities are not transfer-design outputs: they are computed entirely from overlap geometry and the contributing atlas histories, without constructing or correcting a continuous trajectory between families $A$ and $B$. Each is a screening metric that approximates what a corrected family-to-family transfer would cost (or take) at the local overlap region, a proxy in the sense that it stands in for an underlying transfer quantity that has not yet been computed.

For a given overlap voxel \(Z_k\), family \(A\) contributes departing trajectories that reach the voxel, while family \(B\) contributes approaching trajectories that are compatible with it.
Each contributing trajectory already carries the turning cost associated with the heading-change maneuver applied at its seed.

Let,
\begin{equation}
\Delta V_{\mathrm{turn},A}^{\min}(Z_k)
\end{equation}
denote the minimum departure-side turning cost represented in \(Z_k\), and let,
\begin{equation}
\Delta V_{\mathrm{turn},B}^{\min}(Z_k)
\end{equation}
denote the corresponding minimum approach-side turning cost.
These quantities represent the least turning effort carried into the voxel from the two families.

Because the overlap condition is imposed at finite voxel resolution, it does not imply exact state coincidence within \(Z_k\).
A local patching term is therefore introduced to represent the unresolved directional mismatch at the voxel scale.
If \((x_k,y_k)\) denotes the voxel center and \(\Delta\theta\) is the heading bin width, then,
\begin{equation}
\Delta V_{\mathrm{patch}}(Z_k)
=
2\,v_k\,\sin\!\left(\frac{\Delta\theta}{2}\right),
\label{eq:dv_patch}
\end{equation}
where
\begin{equation}
v_k
=
\sqrt{-2\bar U(x_k,y_k)-C_J},
\qquad
C_J=3.1294.
\label{eq:vk_common}
\end{equation}
To confirm that the patching term acts as a voxel-scale correction rather than the dominant cost contribution, at the maximum-budget reference case, \(\Delta V_{\mathrm{patch}}\) over the finite direct pairs ranges from \(0.46\) to \(35.92~\mathrm{m/s}\), with median \(2.83~\mathrm{m/s}\), compared with a total proxy-cost range of \(0.62\) to \(394.82~\mathrm{m/s}\) and median \(42.76~\mathrm{m/s}\).
The largest patching term occurs for the \(3{:}1\) stable resonant (R31-S) to \(3{:}1\) unstable resonant (R31-U) pair, whose total proxy cost is \(236.56~\mathrm{m/s}\).

The base turning contribution associated with \(Z_k\) is,
\begin{equation}
\Delta V_{\mathrm{base}}(Z_k)
=
\Delta V_{\mathrm{turn},A}^{\min}(Z_k)
+
\Delta V_{\mathrm{turn},B}^{\min}(Z_k),
\label{eq:dv_base}
\end{equation}
and the voxel-level proxy cost is defined by
\begin{equation}
\Delta V_{\mathrm{proxy}}(Z_k)
=
\Delta V_{\mathrm{base}}(Z_k)
+
\Delta V_{\mathrm{patch}}(Z_k).
\label{eq:dv_proxy}
\end{equation}

A corresponding proxy time is defined from the trajectories contributing to the same voxel.
Let,
\begin{equation}
T_{A}^{\mathrm{mean}}(Z_k),
\end{equation}
denote the average time taken by the departing trajectories from family \(A\) that reach \(Z_k\), and let,
\begin{equation}
T_{B}^{\mathrm{mean}}(Z_k),
\end{equation}
denote the average time associated with the approaching trajectories toward family \(B\) that pass through \(Z_k\).
Then, the voxel-level proxy time is,
\begin{equation}
T_{\mathrm{proxy}}(Z_k)
=
T_{A}^{\mathrm{mean}}(Z_k)
+
T_{B}^{\mathrm{mean}}(Z_k).
\label{eq:t_proxy}
\end{equation}

The pairwise aggregation is obtained by selecting the overlap voxel with minimum proxy cost,
\begin{equation}
Z_{A,B}^{\star}
=
\arg\min_{Z_k\in\mathcal{O}_{A\rightarrow B}}
\Delta V_{\mathrm{proxy}}(Z_k).
\label{eq:best_voxel}
\end{equation}
The direct pairwise proxy cost is then defined by,
\begin{equation}
\Delta V_{A,B}
=
\Delta V_{\mathrm{proxy}}(Z_{A,B}^{\star}),
\label{eq:pairwise_proxy_metric}
\end{equation}
and the associated direct pairwise proxy time is taken from the same selected voxel,
\begin{equation}
T_{A,B}
=
T_{\mathrm{proxy}}(Z_{A,B}^{\star}).
\label{eq:pairwise_proxy_time}
\end{equation}

Thus, each family pair is reduced from a set-valued overlap to a single representative direct candidate in the reduced model. 
If \(\mathcal{O}_{A\rightarrow B}=\varnothing\), then no candidate direct connection is assigned for that pair under the prescribed budgets. 
Otherwise, the aggregated pairwise quantities \(\Delta V_{A,B}\) and \(T_{A,B}\) provide the direct family-to-family screening metrics used in the matrix construction of Section~\ref{subsec:matrix}.


\subsection{Family-to-Family Direct Accessibility Matrix}\label{subsec:matrix}

For each family pair \(A\) and \(B\), the voxel-level aggregation of Section~\ref{subsec:ranking} yields a direct proxy cost \(\Delta V_{A,B}\) and an associated proxy time \(T_{A,B}\) whenever the overlap set \(\mathcal{O}_{A\rightarrow B}\) is nonempty.
Collecting these quantities over the representative orbit set defines the weighted family-to-family accessibility matrices,
\begin{equation}
\bigl(\mathbf{W}^{(\Delta V)}\bigr)_{A,B} = \Delta V_{A,B},
\qquad
\bigl(\mathbf{W}^{(T)}\bigr)_{A,B} = T_{A,B}.
\end{equation}
Because the planar CR3BP is treated together with time-reversal symmetry about the \(x\)-axis (Section~\ref{subsec:time-reversal}), the family-level accessibility relation is taken to be undirected, so that,
\begin{equation}
\Delta V_{A,B}
=
\Delta V_{B,A},
\qquad
T_{A,B}
=
T_{B,A},
\label{eq:family_direct_symmetry}
\end{equation}
and one symmetric value is retained for each family pair.

Figure~\ref{fig:baseline_direct_proxy_matrix} shows the direct proxy-cost matrix at the maximum-budget reference case.
At the family level, this matrix functions as a $\Delta V$ map across the representative cislunar family set: each off-diagonal entry quantifies the proxy maneuver cost required to access one family directly from another on the common Jacobi manifold. 
Families are ordered by increasing median direct proxy cost to clarify the weighted structure of the matrix.
The matrix has $\binom{13}{2} = 13 \cdot 12 / 2 = 78$ distinct off-diagonal entries, one per unordered family pair, and families are ordered along both axes by increasing median direct proxy cost to clarify the weighted structure. 
At the maximum-budget reference case, 75 of the 78 family-to-family transfers were found to be within the maximum propulsive and time-of-flight budgets.

\par\medskip
\begin{center}
\captionsetup{type=figure}
    \centering
    \includegraphics[width=1\textwidth]{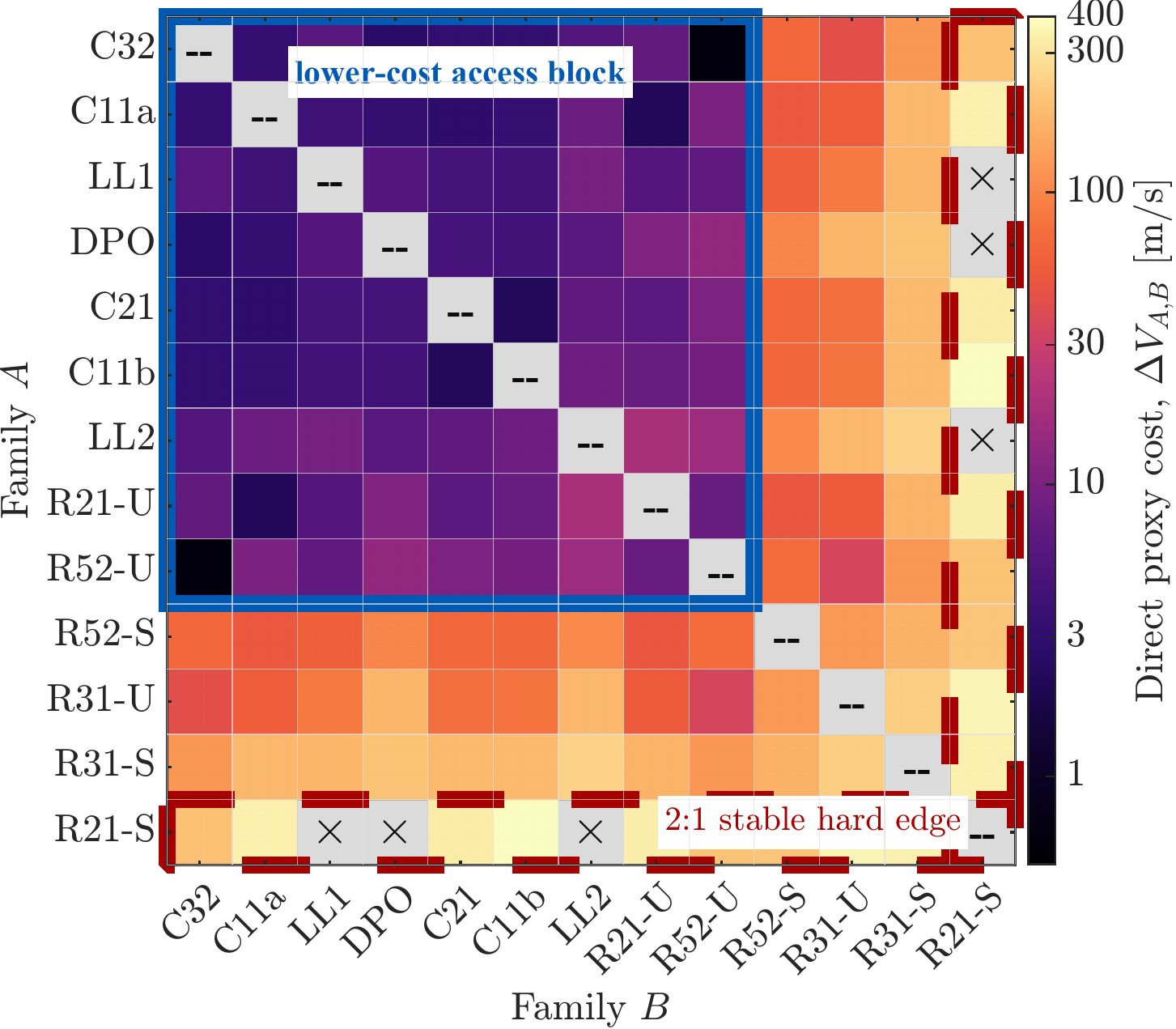}
    \caption{\footnotesize Family-to-family direct proxy-cost matrix \([\Delta V_{A,B}]\) at the maximum-budget reference case, shown on a logarithmic color scale with families ordered by increasing median finite direct proxy cost.
    Annotations highlight the lower-cost access block extending from (3,2)-cycler (C32) through 5:2 unstable resonant (R52-U) and the hard 2:1 stable resonant (R21-S) edge.
    The `\(\times\)' markers denote pairs with no direct overlap under the maximum-budget reference case constraints; diagonal entries are omitted.}
    \label{fig:baseline_direct_proxy_matrix}
\end{center}
\medskip

The main feature of the maximum-budget matrix is not the near-completeness of the matrix, but the substantial spread among the finite values of \(\Delta V_{A,B}\).
Direct family accessibility is therefore already highly structured before any network description is introduced, and the weighted matrix carries substantially more information than the binary question of whether a direct pair is present.

A clear pattern in the maximum-budget matrix is that the cycler families, together with the Lyapunov families and the distant prograde orbit, occupy much of the lower-cost portion of the direct accessibility structure, whereas the resonant families populate much of the more expensive portion.
This separation is only partial, however, and the resonant subset exhibits a nontrivial internal ordering rather than a uniform behavior.

\revone{
The low-cost access block is also consistent with the underlying dynamical
transport structure of the unstable representatives. As discussed in
Section~\ref{subsec:families}, these representatives occupy a common chaotic
region of the Jacobi manifold in which natural transport pathways can bring
their neighborhoods into close dynamical proximity. The LL1--LL2 pair provides
a familiar example: heteroclinic and near-heteroclinic connections between
the $L_1$ and $L_2$ Lyapunov families are well known
\citep{koon2000,gomez2004}, and the low proxy cost recovered here
($\Delta V_{\mathrm{LL1,LL2}}\approx 9.36~\mathrm{m/s}$) is consistent with
that manifold-mediated transport structure. More generally, the low-cost
block shows that the reachable-set overlap construction is able to detect
such dynamical proximity across several distinct unstable family types rather
than only for this individual Lyapunov pair.
This reading has a sharper form in the limit of a large time-of-flight budget, where, if their stable and unstable manifolds belong to the same connected chaotic region, the
base turning cost of an unstable--unstable pair is expected to vanish while stable
resonant targets retain a nonzero floor. That limit is examined in
Appendix~\ref{app:longtof}.}

In particular, the matrix does not support a clean division based only on stable versus unstable resonant type.
At the aggregate level, the unstable resonant families do tend to lie on the less expensive side of the stable resonant families, but the trend is not absolute.
Within the stable resonant subset, the R52-S family is markedly less expensive than the R31-S and R21-S families.
At the most difficult end of the maximum-budget matrix, the R21-S family is the most isolated case, with median finite direct proxy cost about \(340~\mathrm{m/s}\), compared with about \(182~\mathrm{m/s}\) for the R31-S family and about \(68~\mathrm{m/s}\) for the R52-S family.
Thus, the direct accessibility hierarchy resolved by the matrix is more detailed than a simple stability-based classification.

The three missing direct pairs are LL1  to the R21-S family, LL2 to the R21-S family, and the distant prograde orbit (DPO) to the R21-S family.
Within the present framework, these absences mean only that no direct overlap was found for those pairs at the maximum-budget reference case; they do not imply global impossibility.
Their placement is nevertheless informative because they are localized at the hardest-to-reach edge of an otherwise nearly complete weighted direct accessibility matrix.

The matrices \(\mathbf{W}^{(\Delta V)}\) and \(\mathbf{W}^{(T)}\) therefore provide the direct family-to-family accessibility description at the maximum-budget reference case.
Section~\ref{sec:network} takes \(\mathbf{W}^{(\Delta V)}\) as the starting point for the cislunar network interpretation.

\section{Maximum-Budget Reference Cislunar Orbital Network}\label{sec:network}

The direct family-to-family accessibility matrix of Section~\ref{sec:overlap} carries more information than its near-completeness at the maximum-budget reference case suggests. The substantial spread in finite proxy costs across family pairs implies that some families are inexpensively accessible from many others, some serve as efficient global staging points, and some lie on admissible minimum-cost multileg routes between still other families. To resolve these distinctions, the maximum-budget direct accessibility matrix is interpreted as an undirected weighted graph: the \emph{cislunar orbital network} on the common Jacobi manifold at $C_J = 3.1294$.

Three mission-design roles are read from this network.
\emph{Hub} families have many inexpensive direct links to other families.
\emph{Gateway} families provide comparatively low-cost access to a broad portion of the family set.
\emph{Relay} families frequently lie on admissible minimum-cost multileg paths between other families.
Each role is an operational interpretation of a corresponding cost-aware centrality measure introduced in Section~\ref{subsec:centrality}

As in Section~\ref{subsec:families}, ``family'' refers throughout to the representative periodic orbit at the common Jacobi level. The remainder of the section constructs the network from the direct accessibility matrix, defines the cost-aware centrality measures, and identifies the maximum-budget hub, gateway, and relay structure. Additional numerical validation is reported in Appendix~\ref{app:baseline_validation}.


\subsection{Network Construction from the Direct Accessibility Matrix}

Let \(\mathcal{S}\) denote the representative family set of Section~\ref{sec:framework}, with $N_f = |\mathcal{S}|=13$.
For prescribed total budgets $(\Delta V_{\mathrm{cap}},T_{\mathrm{cap}})$, the pairwise direct  proxy propulsion costs $\Delta V_{A,B}$ and times-of-flight $T_{A,B}$ from Section~\ref{sec:overlap} define an undirected weighted graph \citep{newman2010},
\begin{equation}
\mathcal{G}(\Delta V_{\mathrm{cap}},T_{\mathrm{cap}})
=
(\mathcal{V},\mathcal{E},w),
\end{equation}
in which each node in $\mathcal{V}$ corresponds to one representative family in $\mathcal{S}$.

An undirected edge \((A,B)\in\mathcal{E}\) is retained when the family pair admits a direct candidate transfer between $A$ and $B$ within the total maneuver budget, that is, when the direct proxy cost is finite and does not exceed $\Delta V_{\mathrm{cap}}$:
\begin{equation}
(A,B)\in\mathcal{E}
\quad\Longleftrightarrow\quad
\Delta V_{A,B}
\le
\Delta V_{\mathrm{cap}}.
\label{eq:network_edge_rule}
\end{equation}
Each retained edge carries the cost weight,
\begin{equation}
w(A,B)
=
\Delta V_{A,B}.
\label{eq:edge_weight_def}
\end{equation}
The corresponding proxy time $T_{A,B}$ is recorded  as an auxiliary edge attribute and is used later for multileg path interpretation, but does not enter the edge-retention rule or the centrality measures of Section~\ref{subsec:centrality}..


At the maximum-budget reference case, \(\Delta V_{\mathrm{cap}}=409.3~\mathrm{m/s}\) and \(T_{\mathrm{cap}}=27.32~\mathrm{days}\); the corresponding maximum-budget cislunar orbital network is denoted $\mathcal{G}_{\max}$
and is shown in
Figure~\ref{fig:baseline_reduced_family_network},
with families arranged clockwise by increasing median direct proxy cost and the edge color encoding \(w(A,B)=\Delta V_{A,B}\) on a logarithmic scale.
\(\mathcal{G}_{\max}\) contains the full 13-family representative set and 75 of the 78 distinct family pairs as direct edges;
the three missing direct pairs are those involving the R21-S family, as identified in Section~\ref{subsec:matrix}.
The maximum-budget reference network is therefore nearly complete in a binary sense but remains strongly nonuniform in its edge costs.

Multileg routes on $\mathcal{G}$ are defined and constrained as follows.
For a path \(P=(A_0,A_1,\dots,A_m)\) with \(A_0,\dots,A_m\in\mathcal{S}\) and \((A_\ell,A_{\ell+1})\in\mathcal{E}\) for \(\ell=0,\dots,m-1\), the cumulative proxy maneuver cost is,
\begin{equation}
\Delta V(P)
=
\sum_{\ell=0}^{m-1}
w(A_\ell,A_{\ell+1}),
\label{eq:path_cost_def}
\end{equation}
and $P$  is regarded as admissible  when,
\begin{equation}
\Delta V(P)\le \Delta V_{\mathrm{cap}}.
\label{eq:path_budget_rule}
\end{equation}
Multileg feasibility is enforced  through cumulative maneuver cost alone. Cumulative multileg time-of-flight is not imposed as a separate path constraint because the family-level abstraction aggregates over orbital phase, so the phase dependent coasting time at each intermediate family is not represented within the present network. Phase-resolved time accounting is a natural extension and is discussed as future work in Section~\ref{sec:conclusion}.

\par\medskip
\begin{center}
\captionsetup{type=figure}
    \centering
    \includegraphics[width=0.9\textwidth]{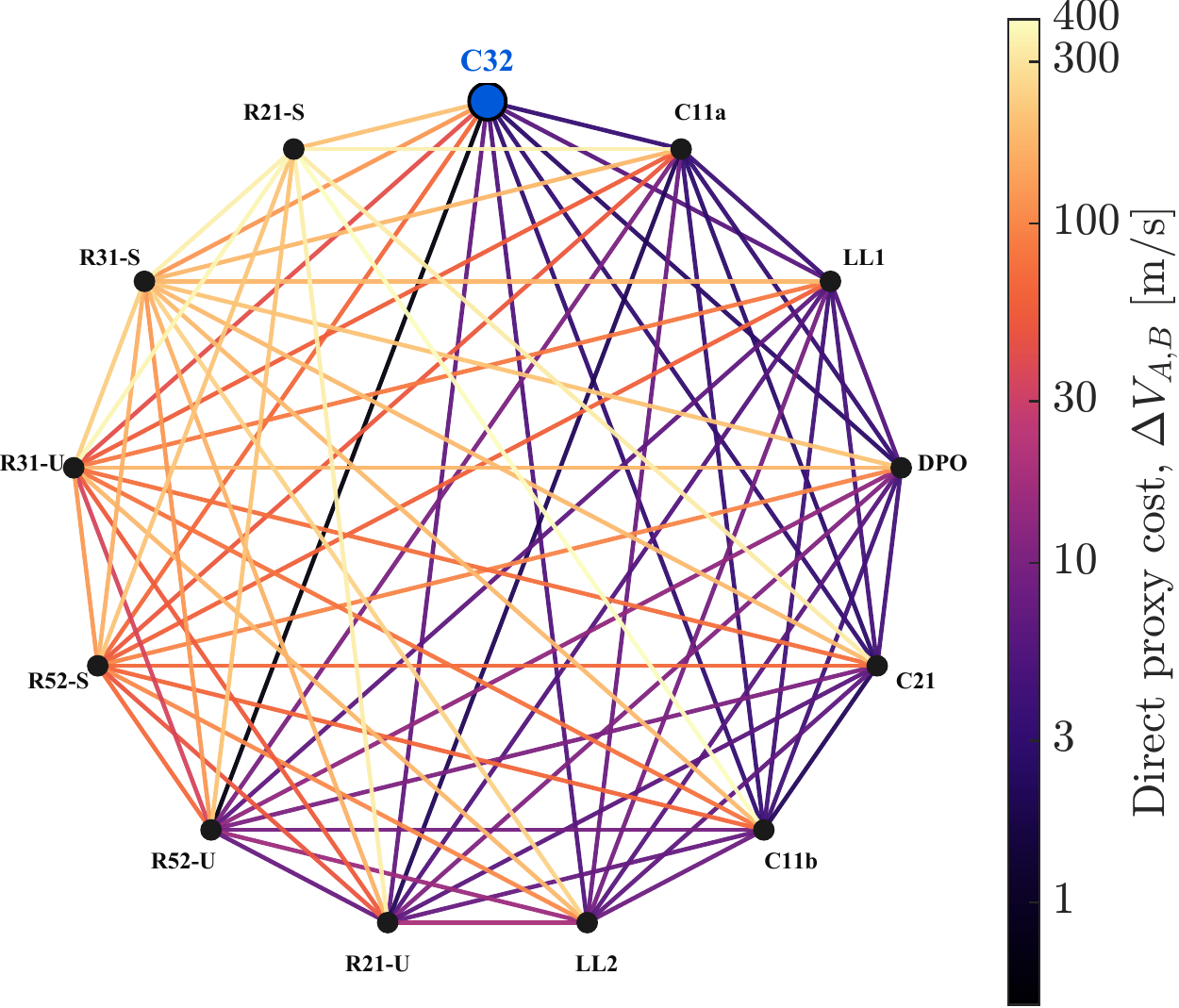}
    \caption{\footnotesize Maximum-budget cislunar orbital network.
    Families are arranged clockwise by increasing median direct proxy cost, and edge color indicates the direct proxy cost \(w(A,B)=\Delta V_{A,B}\) on a logarithmic scale.
   The (3,2)-cycler (C32) is highlighted because it is the dominant hub, gateway, and relay family at the maximum-budget reference case.
    Missing direct pairs appear as absent edges.}
    \label{fig:baseline_reduced_family_network}
\end{center}
\medskip


\subsection{Cost-Aware Centrality Measures and Their Physical Meaning}\label{subsec:centrality}

To interpret the cislunar orbital network, three cost-aware node measures are used \citep{SaxenaIyengar2020Survey}.
The three are chosen so that they resolve distinct aspects of accessibility within the same weighted family graph: direct-access prominence, global accessibility through both direct and multileg paths, and intermediary role on admissible minimum-cost routes.
Throughout this subsection, edge costs are the direct proxy costs \(w(A,B)=\Delta V_{A,B}\) introduced in Section~5.1.

For any two distinct families \(A,B\in\mathcal{S}\), let the budget-feasible shortest-path cost be,
\begin{equation}
d(A,B)
=
\min_{P:A\rightarrow B}\Delta V(P),
\label{eq:budget_feasible_distance}
\end{equation}
where the minimum is taken over all admissible multileg paths \(P\) from \(A\) to \(B\), with admissibility defined by \eqref{eq:path_budget_rule};
if no such path exists, then \(d(A,B)=\infty\). The quantity \(d(A,B)\) is the minimum cumulative proxy maneuver cost to reach family $B$ from family $A$  on the  network within the prescribed total maneuver budget \citep{newman2010,dijkstra1959}.

\paragraph{Strength.}
The first measure is the normalized strength \citep{barrat2004}. 
Because larger direct proxy cost corresponds to weaker direct accessibility, strength is computed using reciprocal-cost edge weights: 
\begin{equation}
\widehat{w}(A,B)
=
\frac{1}{w(A,B)}
=
\frac{1}{\Delta V_{A,B}},
\label{eq:reciprocal_cost_weight}
\end{equation}
and the normalized strength of family \(A\) is,
\begin{equation}
S(A)
=
\frac{1}{N_f-1}
\sum_{B\in\Gamma(A)}
\widehat{w}(A,B),
\label{eq:strength_def}
\end{equation}
where \(\Gamma(A)\) denotes the set of families directly connected to \(A\). 
A family with high normalized strength has many inexpensive direct links and is therefore interpreted as a \emph{hub}.

\paragraph{Harmonic closeness.}
The second measure is the normalized harmonic closeness,
\begin{equation}
H(A)
=
\frac{1}{N_f-1}
\sum_{\substack{B\in\mathcal{S}\\B\neq A}}
\frac{1}{d(A,B)}.
\label{eq:harmonic_closeness_def}
\end{equation}
The harmonic form is used in place of standard closeness  because $d(A,B)$ can be infinite under the cap-feasibility rule whenever no  admissible path between $A$ and $B$ exists; 
reciprocal terms remain well defined in that case and contribute zero, so the measure is meaningful across the full budget plane considered in Section~\ref{sec:budget} \citep{marchiori2000,latora2001,boldi2014}.
A family with high harmonic closeness reaches the rest of the family set, on average, through a low minimum path cost, whether by  direct or multileg routes, and is interpreted as a \emph{gateway} or \emph{staging} family.

\paragraph{Betweenness.}
The third measure is the normalized betweenness \citep{freeman1977},
\begin{equation}
B(A)
=
\frac{2}{(N_f-1)(N_f-2)}
\sum_{\substack{P<Q\\P,Q\in\mathcal{S}\setminus\{A\}}}
\frac{\sigma_{PQ}(A)}{\sigma_{PQ}},
\label{eq:betweenness_def}
\end{equation}
where \(\sigma_{PQ}\) denotes the number of admissible shortest paths between families \(P\) and \(Q\), and \(\sigma_{PQ}(A)\) denotes the number of those paths that pass through \(A\) as an intermediate family.
If no admissible path exists between \(P\) and \(Q\), the corresponding contribution is taken to be zero.
A family with high betweenness lies frequently on admissible minimum-cost paths between other families and is interpreted as a \emph{relay} family.

These three measures are complementary:
strength resolves direct-access prominence, 
harmonic closeness resolves global accessibility  through admissible direct and multileg paths, and 
betweenness resolves intermediary routing. 
The three need not agree in a general weighted graph, and Section~\ref{subsec:budget_winners} shows that they nevertheless agree substantially in the present cislunar orbital network, which is a structural feature of the network itself rather than of the measures.


\subsection{Maximum-Budget Family Roles: (3,2)-Cycler as the Dominant Hub, Gateway, and Relay}\label{subsec:centralitiesrank}
(3,2)-cycler (C32) is the dominant family at the maximum-budget reference case, 
ranking first in strength, harmonic closeness, and betweenness.
It is therefore simultaneously the strongest direct-access family, the most efficient global-access family, and the dominant relay family in the cislunar orbital network.
Figure~\ref{fig:baseline_family_centralities} and Table~\ref{tab:baseline_family_centrality} report the underlying values.

These rankings are not determined by binary connectivity alone. 
Although \(\mathcal{G}_{\max}\) contains 75 of the 78 possible direct edges, the nonuniform 
edge-cost structure makes multileg paths competitive even where direct connections exist: 
44 of the 75 directly connected pairs admit a cheaper budget-feasible multileg route, 
and 47 of the 78 family pairs overall attain their minimum-cost admissible route through 
a multileg path rather than a single direct edge. Shortest-path centralities therefore 
remain informative even in this nearly complete direct network.

\par\medskip
\begin{center}
\captionsetup{type=figure}
    \centering
    \includegraphics[width=1\textwidth]{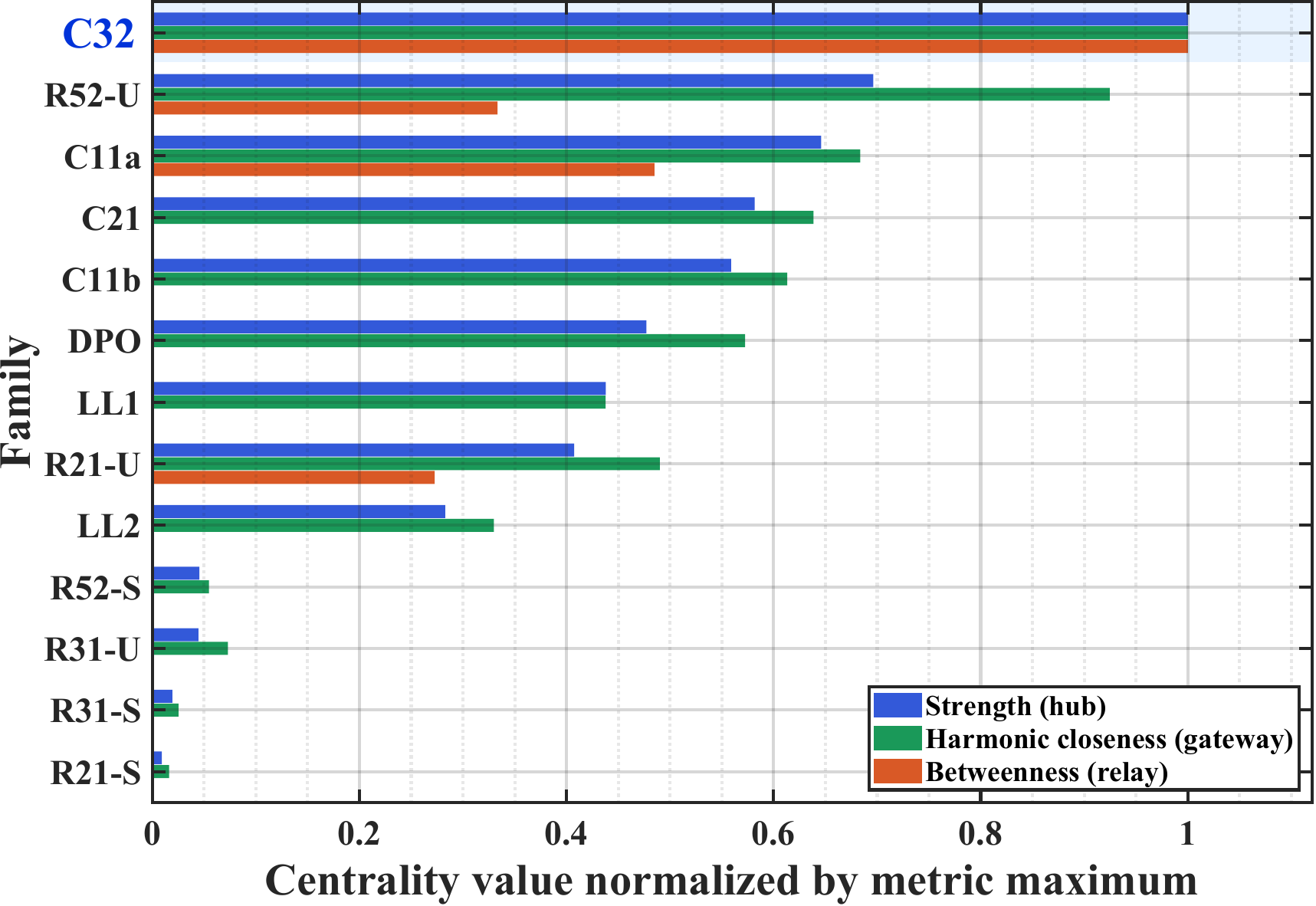}
    \caption{\footnotesize  Maximum-budget family role centralities in the cislunar orbital network.
    Bars show strength, harmonic closeness, and betweenness, each normalized to the maximum value of its metric, enabling direct comparison of hub, gateway, and relay signatures across families.
    (3,2)-cycler (C32) ranks first in all three measures, while other families show distinct role signatures; for example, 5:2 unstable resonant (R52-U) is strong in strength and harmonic closeness, whereas (1,1)a-cycler (C11a) is comparatively stronger in betweenness.
    The corresponding raw values and ranks are listed in Table~\ref{tab:baseline_family_centrality}.}
    \label{fig:baseline_family_centralities}
\end{center}
\medskip

The betweenness value \(B=0.500\) gives a more specific interpretation of the relay role: under the standard undirected normalization of Section~\ref{subsec:centrality}, half of the unordered family pairs not involving C32 have an admissible minimum-cost route that passes through it.
C32 is therefore not merely inexpensive to access directly; it also organizes a large fraction of the maximum-budget minimum-cost multileg structure.

The second tier exposes two opposite role asymmetries.
5:2 unstable resonant (R52-U) ranks second in strength and second in harmonic closeness, but only third in betweenness,
so its dominance is weighted toward direct-access and gateway prominence rather than relay prominence.
By contrast, (1,1)a-cycler (C11a) shows the inverse pattern:  third in strength and third in harmonic closeness, but second in betweenness, 
marking it as a more important  relay than its direct-access rank alone would suggest.

A sharper version of the relay-without-direct-access pattern appears in R21-U.
It ranks eighth in strength and seventh in harmonic closeness, yet fourth in betweenness,
so its network role derives almost entirely from repeated participation in admissible minimum-cost routes rather than from broad direct access.
The opposite pattern appears in (2,1)-cycler (C21) and (1,1)b-cycler (C11b), which rank fourth and fifth in both strength and harmonic closeness but both have zero betweenness;
these families are well-placed for direct-access and global staging but are not required as intermediaries in the maximum-budget minimum-cost path structure.

At the difficult end of the network, the stable resonant families remain clearly separated from the lower-cost cycler--Lyapunov--DPO portion of the family set.
R21-S is the limiting case:
last in both strength and harmonic closeness, with zero betweenness.
R31-S and R52-S also remain weak in both direct-access and global-access measures, although R52-S is less isolated than the \(3{:}1\) and \(2{:}1\) stable families.
These rankings  reinforce the hard R21-S edge already visible in the direct accessibility matrix of Section~\ref{subsec:matrix} and identify R21-S as the single most isolated family in the maximum-budget reference network.

The maximum-budget reference network was additionally checked against a 15-case coarsening, refinement, and stress envelope validation suite.
The direct topology, the persistent absence of the same three direct pairs, and the dominant role hierarchy remain unchanged across that suite; only the proxy values themselves vary quantitatively.
Full details are reported in Appendix~\ref{app:baseline_validation}.

\par\medskip
\begin{center}
\captionsetup{type=table}
\centering
\small
\caption{\footnotesize Maximum-budget family centrality values in the cislunar orbital network.
Values are reported at the maximum-budget reference case.
Numbers in parentheses give descending rank within each metric, with ties sharing the same rank.
Strength and harmonic closeness are reported in reciprocal-cost units \((\mathrm{m/s})^{-1}\), while betweenness is dimensionless.
Bold entries indicate the top-ranked value in each centrality metric.}
\label{tab:baseline_family_centrality}
\begin{tabular}{lccc}
\toprule
Family & Strength & Harmonic closeness & Betweenness \\
\midrule
\textbf{(3,2)-cycler} & \textbf{0.2850 (1)} & \textbf{0.2891 (1)} & \textbf{0.5000 (1)} \\
5:2 unstable resonant  & 0.1984 (2) & 0.2673 (2) & 0.1667 (3) \\
(1,1)a-cycler & 0.1841 (3) & 0.1976 (3) & 0.2424 (2) \\
(2,1)-cycler & 0.1657 (4) & 0.1846 (4) & 0.0000 (-) \\
(1,1)b-cycler & 0.1593 (5) & 0.1772 (5) & 0.0000 (-) \\
Distant prograde orbit & 0.1360 (6) & 0.1655 (6) & 0.0000 (-) \\
L1 Lyapunov  & 0.1248 (7) & 0.1265 (8) & 0.0000 (-) \\
2:1 unstable resonant  & 0.1161 (8) & 0.1417 (7) & 0.1364 (4) \\
L2 Lyapunov  & 0.0806 (9) & 0.0954 (9) & 0.0000 (-) \\
5:2 stable resonant  & 0.0130 (10) & 0.0158 (11) & 0.0000 (-) \\
3:1 unstable resonant  & 0.0127 (11) & 0.0211 (10) & 0.0000 (-) \\
3:1 stable resonant  & 0.0055 (12) & 0.0073 (12) & 0.0000 (-) \\
2:1 stable resonant  & 0.0026 (13) & 0.0047 (13) & 0.0000 (-) \\
\bottomrule
\end{tabular}
\end{center}
\medskip


\section{Budget-Dependent Accessibility Regimes in the Cislunar Orbital Network}\label{sec:budget}

The cislunar orbital network at the maximum-budget reference case is nearly complete and admits a clean family-role hierarchy.
Most operationally interesting regimes, however, sit below that limit, where tightening either the maneuver or the time budget forces structural choices that the maximum-budget network does not expose.
This section extends the construction of Section~\ref{sec:network} across the sampled two-parameter budget plane \((\Delta V_{\mathrm{cap}},T_{\mathrm{cap}})\) to determine how accessibility and family roles change as maneuver and time budgets are relaxed.

The central result is that direct accessibility, graph connectedness, and budget-feasible multileg closure separate into three distinct regimes across the budget plane:
the network can become connected before all family pairs are mutually reachable, and full multileg closure can occur before the direct graph becomes complete.
The dominant hub, gateway, and relay families also reorganize across the plane, and the transitions between them are controlled by different mechanisms for different centrality measures.

Throughout this section, the graph remains undirected, multileg admissibility is enforced through cumulative maneuver cost, and the hub/gateway/relay interpretations follow the strength, harmonic closeness, and betweenness measures of Section~\ref{subsec:centrality}.
The remainder of the section defines the sampled budget plane and per-budget networks, characterizes the three accessibility regimes, identifies the dominant families across the plane, and discusses mission-design implications.


\subsection{Definition of the Sampled Budget Plane and Budget-Dependent Networks}

The construction of Section~\ref{sec:network} extends to a two-parameter family of networks parameterized by the total maneuver and time budgets. For each sampled budget point
 \((\Delta V_{\mathrm{cap}},T_{\mathrm{cap}})\),
\begin{equation}
\mathcal{G}(\Delta V_{\mathrm{cap}},T_{\mathrm{cap}})
=
(\mathcal{V},\mathcal{E},w),
\end{equation}
is defined using the same family set \(\mathcal{S}\), common Jacobi manifold \(C_J=3.1294\), edge-retention rule \eqref{eq:network_edge_rule}, edge-cost definition \eqref{eq:edge_weight_def}, and path-budget condition \eqref{eq:path_budget_rule}
as in Section~\ref{sec:network}.
The sampled plane contains \(25\times 25=625\) such networks over,
\begin{equation}
51.16 \le \Delta V_{\mathrm{cap}} \le 409.26~\mathrm{m/s},
\qquad
6.83 \le T_{\mathrm{cap}} \le 27.32~\mathrm{days},
\end{equation}
with the maximum-budget reference case of Section~\ref{sec:network} sitting at the upper-right corner point of the plane.
All budget-plane figures in this section use these physical units.
An interactive version of the cislunar orbital network, including a budget explorer in which $\Delta V_{\mathrm{cap}}$ and $T_{\mathrm{cap}}$ can be varied continuously across the sampled plane, is available at the repository URL listed in the Data and Code Availability Statement.

Two sets of diagnostics are recorded at each sampled budget point. 
The first set characterizes the regime structure of the network: the number of direct family pairs, the size of the largest connected component, and the number of unordered family pairs admitting a budget-feasible minimum-cost path.
These three quantities distinguish direct accessibility, graph connectedness, and budget-feasible multileg closure, respectively.
The second set characterizes which family dominates each role:
the maximum-strength family, maximum-harmonic-closeness family, and maximum-betweenness family  at each budget point.
The resulting \emph{winner maps} track how the dominant hub, gateway, and relay families reorganize across the budget plane.
Appendix~\ref{app:winner_flips} provides additional diagnostics for the centrality winner transitions.


\subsection{Three Accessibility Regimes: Direct Access, Connectedness, and Multileg Closure}\label{subsec:regimes}


All three diagnostics shown in Figure~\ref{fig:section6_growth_metrics} increase monotonically as either \(\Delta V_{\mathrm{cap}}\) or \(T_{\mathrm{cap}}\) is relaxed, but they do not become equivalent.
\par\medskip
\begin{center}
\captionsetup{type=figure}
    \centering

    \begin{subfigure}[t]{0.32\textwidth}
        \centering
        \includegraphics[width=\textwidth]{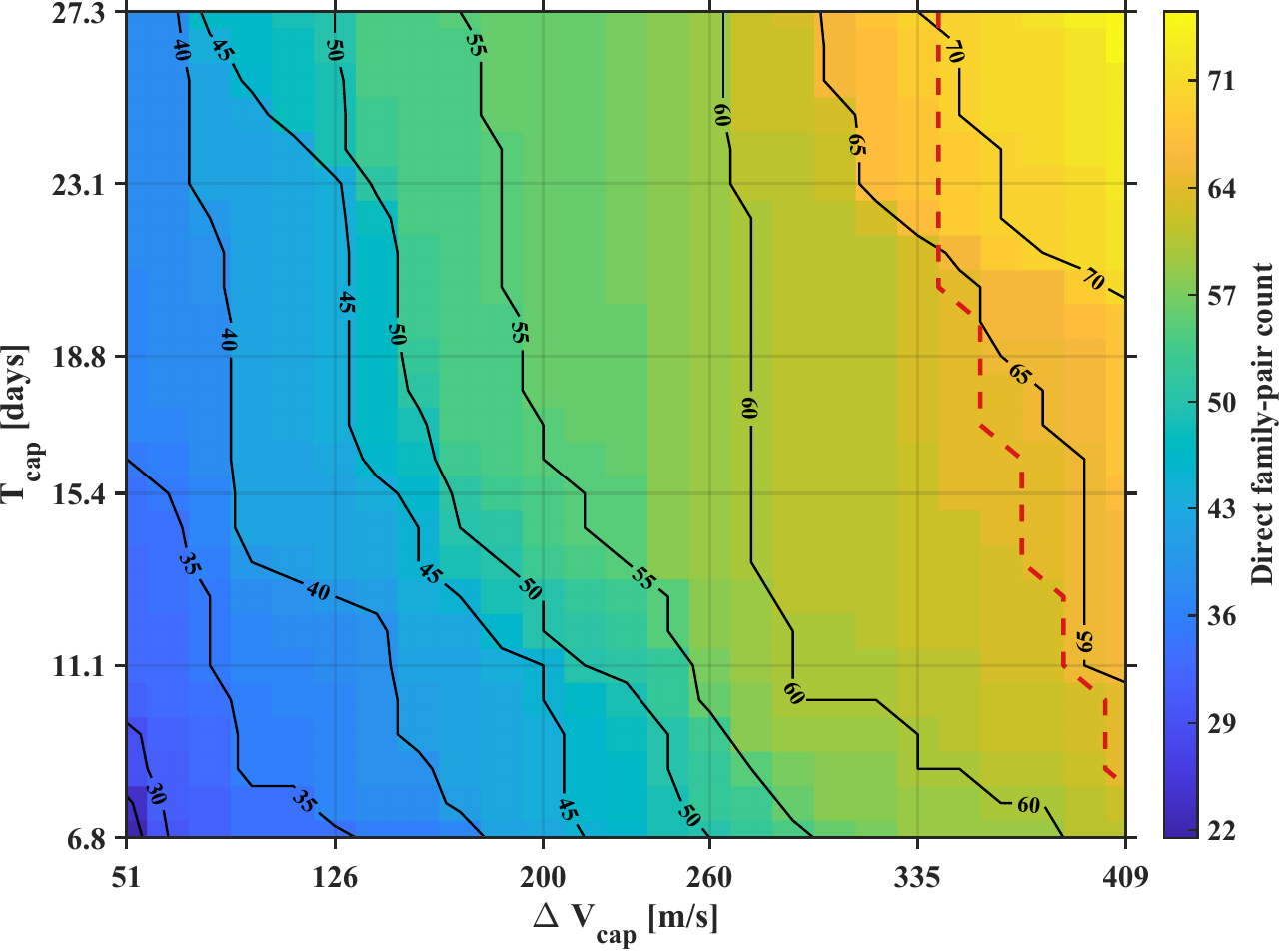}
        \caption{\footnotesize Direct family-pair count.}
    \end{subfigure}
    \begin{subfigure}[t]{0.32\textwidth}
        \centering
        \includegraphics[width=\textwidth]{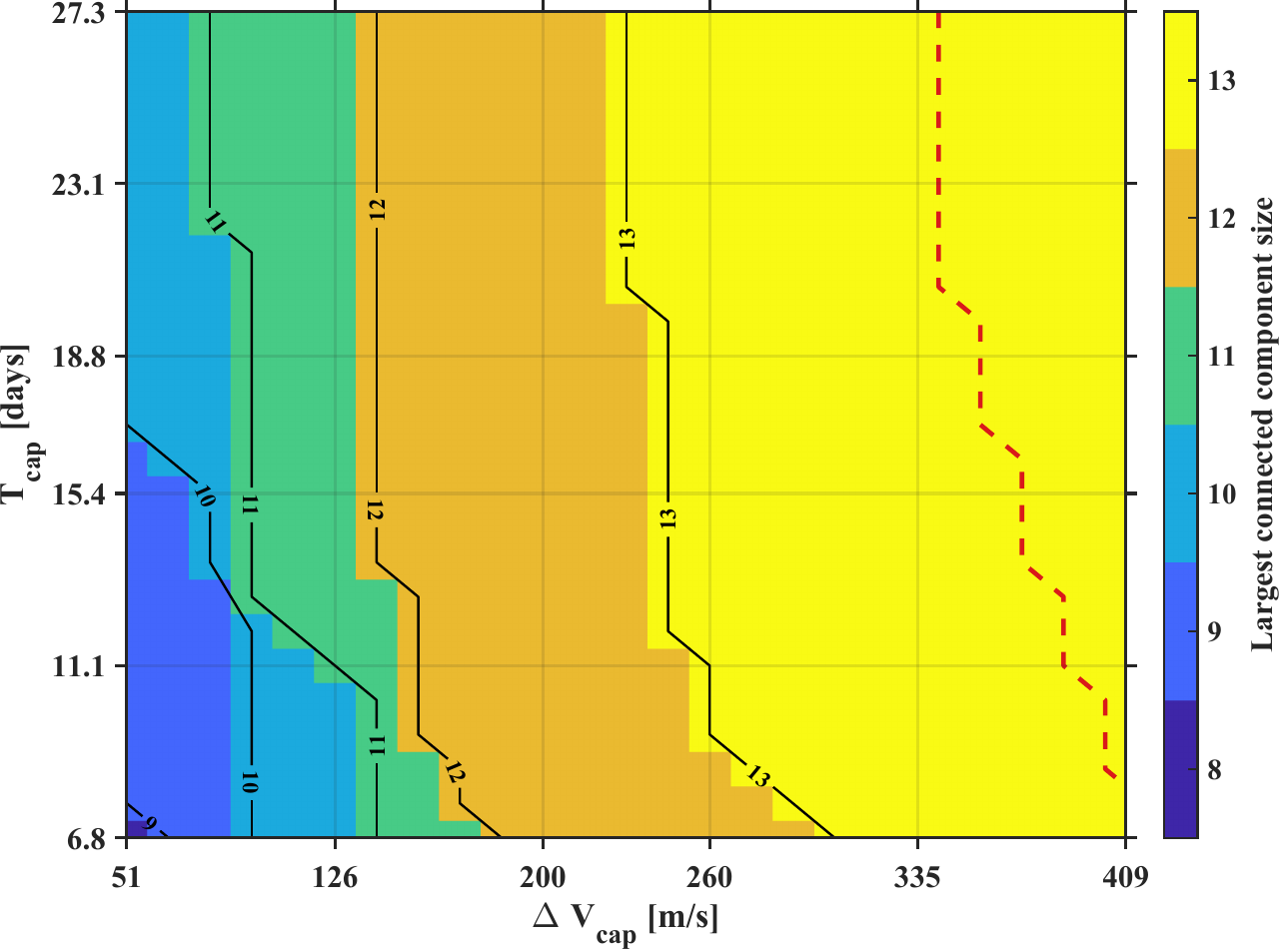}
        \caption{\footnotesize Largest connected component size.}
    \end{subfigure}
    \begin{subfigure}[t]{0.32\textwidth}
        \centering
        \includegraphics[width=\textwidth]{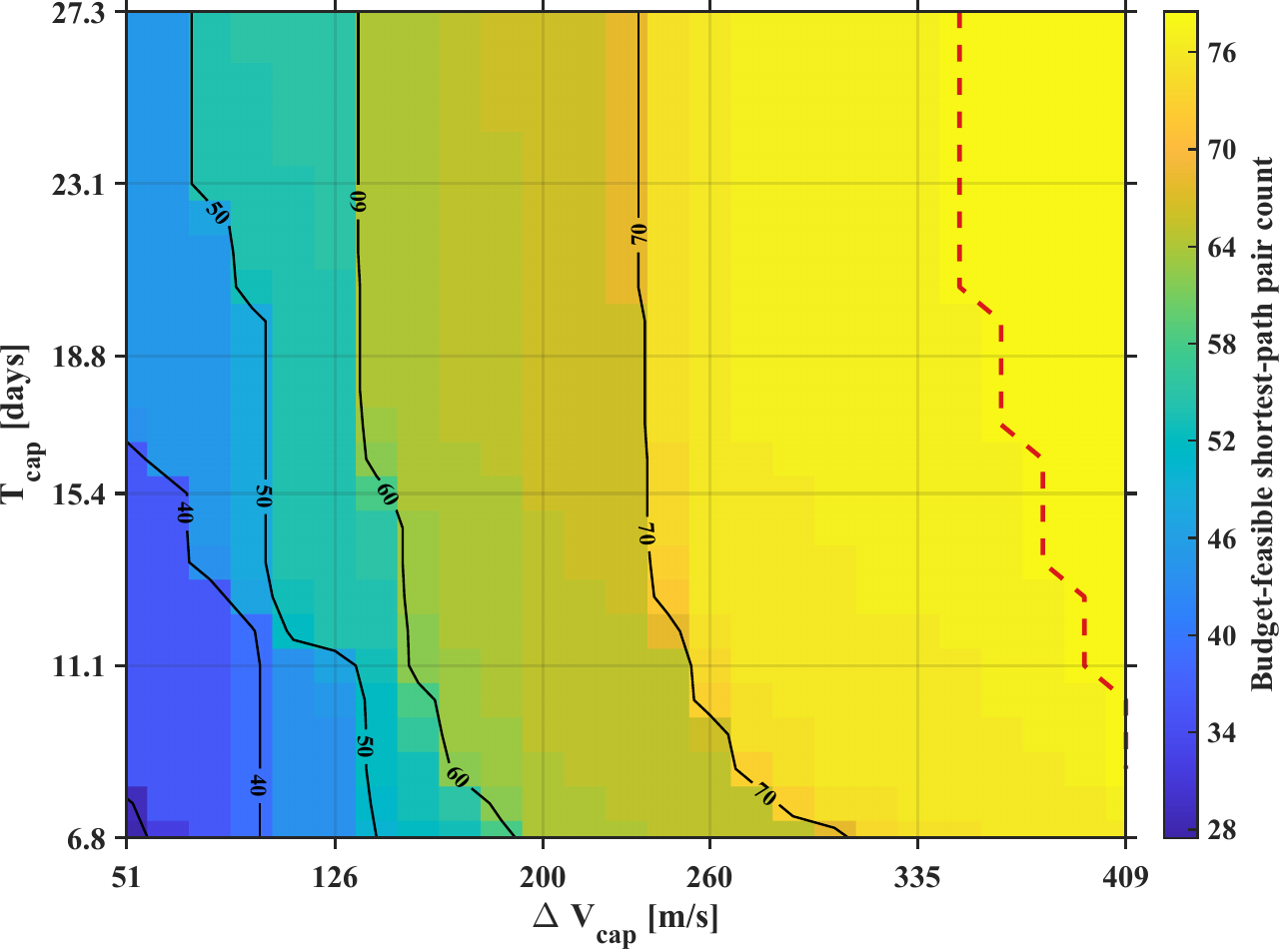}
        \caption{\footnotesize Budget-feasible shortest-path pair count.}
    \end{subfigure}

    \caption{\footnotesize Growth of the cislunar orbital network over the budget plane.
    (a) Direct family-pair count.
    (b) Largest connected component size.
    (c) Number of budget-feasible minimum-cost pairs.
    The red dotted contour marks the onset of full budget-feasible shortest-path connectivity.
    Axes are plotted in total budgets, with \(\Delta V_{\mathrm{cap}}\) on the horizontal axis and \(T_{\mathrm{cap}}\) on the vertical axis.}
    \label{fig:section6_growth_metrics}
\end{center}
\medskip
Across the sampled plane, the direct family-pair count increases from 22 to 75, the largest connected component increases from 8 to 13 families, and the budget-feasible minimum-cost pair count grows from 28 to the full 78.
These three diagnostics measure operationally distinct properties: 
one-leg accessibility, chains of individually admissible direct edges, and cumulative paths that remain within the same total maneuver budget.

The first key result is that full graph connectedness occurs before full all-pairs budget-feasible closure.
The earliest budget point with a 13-family largest connected component occurs at \(\Delta V_{\mathrm{cap}}=230.21~\mathrm{m/s}\) and \(T_{\mathrm{cap}}=11.95~\mathrm{days}\);
at that point the direct graph contains 53 direct family pairs, but only 65 of the 78 unordered family pairs are shortest-path feasible under the total maneuver budget.
The same separation appears when approached from the time-first direction, with the first fully connected component at \(T_{\mathrm{cap}}=6.83~\mathrm{days}\) and \(\Delta V_{\mathrm{cap}}=289.89~\mathrm{m/s}\), again with only 65 shortest-path-feasible pairs.
Graph connectedness and all-pairs budget-feasible closure are therefore distinct regimes.

The second key result is the reverse separation:
full all-pairs budget-feasible closure occurs before the direct graph saturates. 
Complete closure is first reached at \(\Delta V_{\mathrm{cap}}=349.58~\mathrm{m/s}\) with \(T_{\mathrm{cap}}=20.49~\mathrm{days}\), or  at \(T_{\mathrm{cap}}=8.54~\mathrm{days}\) with \(\Delta V_{\mathrm{cap}}=409.26~\mathrm{m/s}\) along the time-first side.
At these threshold points, all 78 unordered family pairs are shortest-path feasible while the direct graph contains only 65 and 66 direct family pairs, respectively.
The network captures genuine indirect accessibility rather than merely redescribing the direct matrix.
The budget-feasible shortest-path pair count exceeds the directly accessible pair count by up to 21 pairs across the sampled plane, so
multileg accessibility is  a persistent structural feature rather than a marginal correction near the high-budget limit.

Even at the most permissive sampled budgets, the direct accessibility graph does not become complete:
the maximum direct family-pair count attained is 75, reached only at the highest sampled maneuver budget together with the two largest sampled time budgets.
The same three missing direct pairs identified in Section~\ref{subsec:matrix} persist throughout the sampled range, and their mission-design significance is discussed in Section~\ref{subsec:mission_implications}.

Taken together, Fig.~\ref{fig:section6_growth_metrics} identifies three nested accessibility regimes:
at low budgets, the network is only partially connected and many family pairs remain inaccessible even through multileg paths;
at intermediate budgets, the direct graph becomes connected while all-pairs budget-feasible closure remains incomplete;
and at higher budgets, all families become mutually reachable through admissible minimum-cost paths even though the direct graph remains incomplete.
This separation between direct accessibility, graph connectedness, and budget-feasible multileg closure is one of the central structural features of the cislunar orbital network.


\subsection{Dominant Family Roles Across the Budget Plane}\label{subsec:budget_winners}

Figure~\ref{fig:section6_winner_regimes} shows the dominant family over the budget plane for each of the three centrality measures.
\par\medskip
\begin{center}
\captionsetup{type=figure}
    \centering

    \begin{subfigure}[t]{0.3\textwidth}
        \centering
        \includegraphics[width=\textwidth]{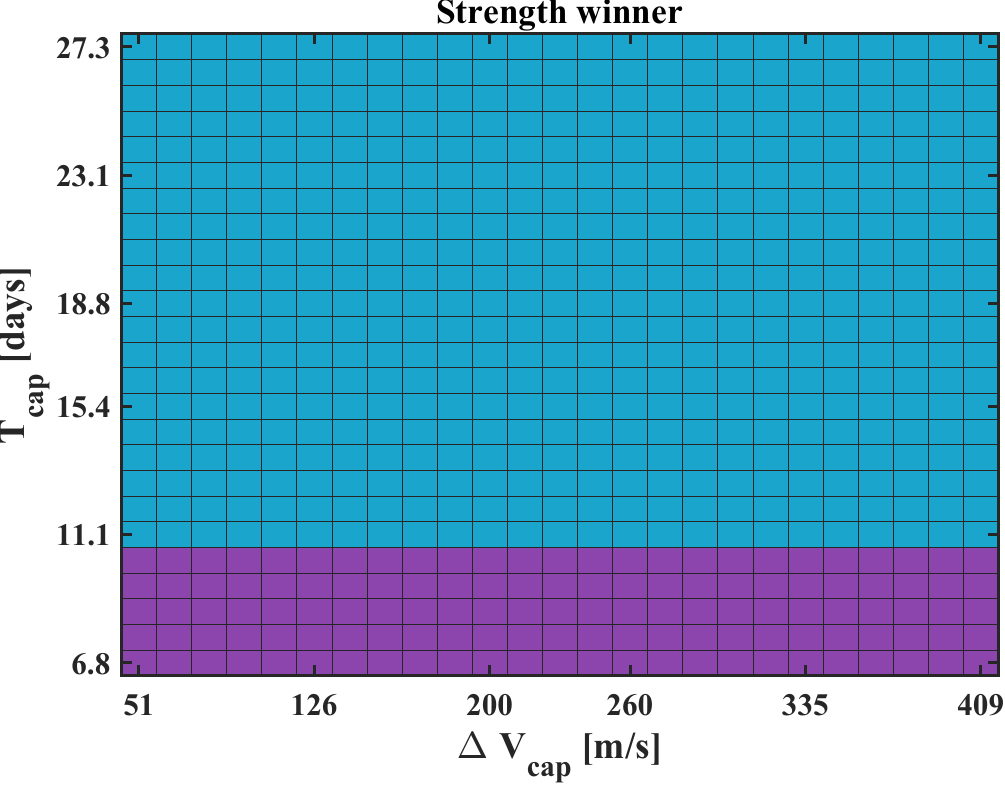}
        \caption{\footnotesize Strength winner.}
    \end{subfigure}
    \hfill
    \begin{subfigure}[t]{0.3\textwidth}
        \centering
        \includegraphics[width=\textwidth]{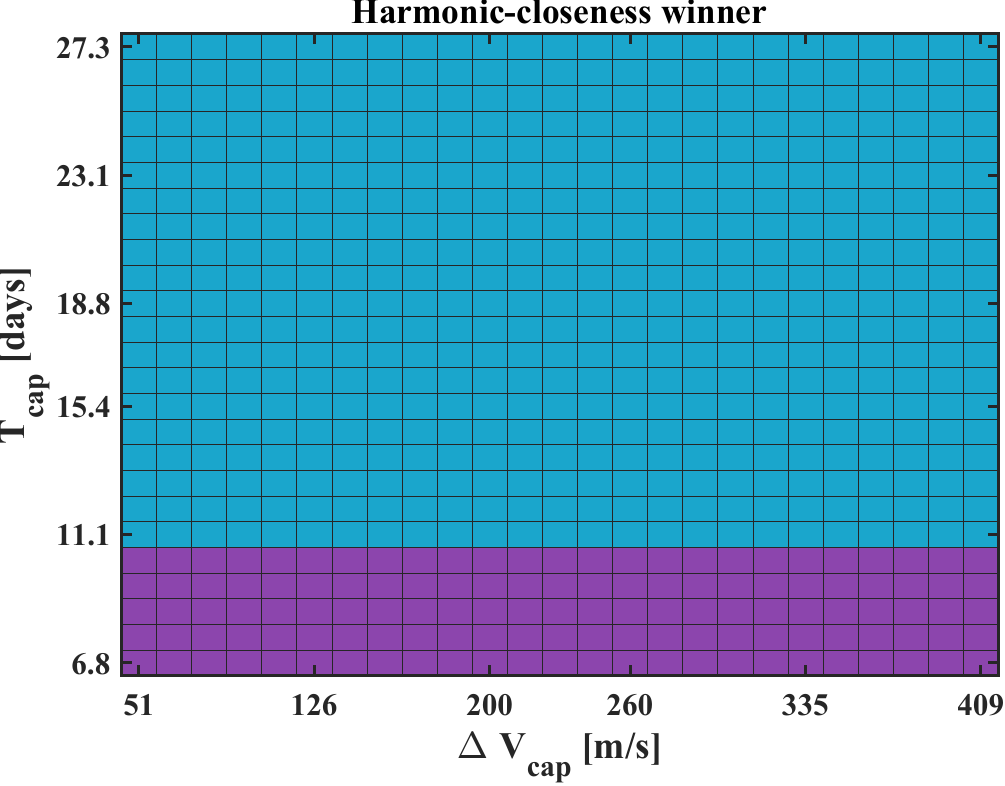}
        \caption{\footnotesize Harmonic-closeness winner.}
    \end{subfigure}
    \hfill
    \begin{subfigure}[t]{0.36\textwidth}
        \centering
        \includegraphics[width=\textwidth]{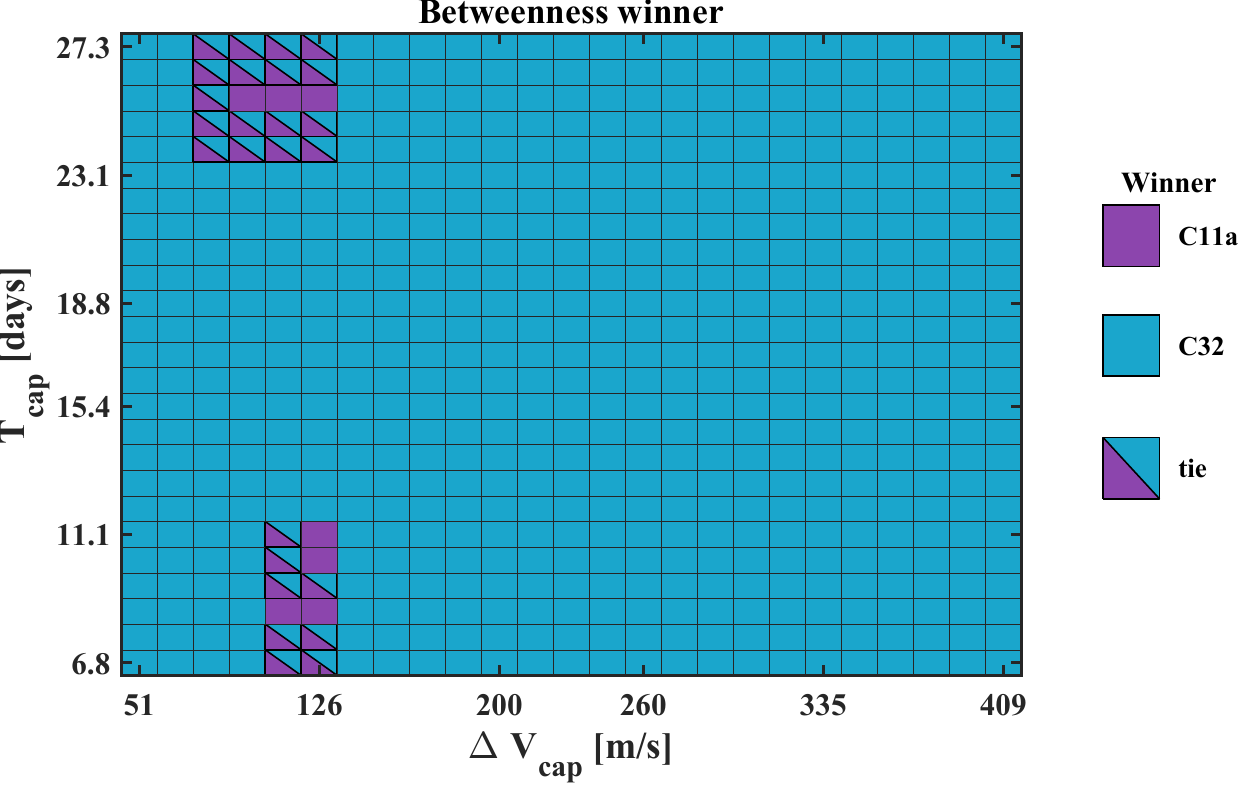}
        \caption{\footnotesize Betweenness winner.}
    \end{subfigure}

    \caption{\footnotesize Winner regimes of the three family-role centralities over the sampled budget plane.
    Panels (a) and (b) show the strength and harmonic-closeness winners, which coincide pointwise across the sampled plane.
    Panel (c) shows the betweenness winner and includes the color key.
    Diagonally split cells with black borders denote tied winner cells.
    Axes use total budgets, with \(\Delta V_{\mathrm{cap}}\) on the horizontal axis and \(T_{\mathrm{cap}}\) on the vertical axis.}
    \label{fig:section6_winner_regimes}
\end{center}
\medskip
At each sampled budget point, the displayed winner is the family attaining the maximum value of the corresponding centrality measure in \(\mathcal{G}(\Delta V_{\mathrm{cap}},T_{\mathrm{cap}})\).
\revone{
If multiple families attain the same maximum value, no tie-breaking rule is
applied; the corresponding budget point is retained and displayed explicitly
as a tied winner cell. In the present sweep, all such winner-level ties occur
between C11a and C32.
}
The clearest structural finding is that only two families, C11a and C32, appear as winners across all three maps and the full budget plane.
The strength and harmonic-closeness winner maps coincide pointwise across the entire sampled plane: in both maps, C11a occupies a low-\(T_{\mathrm{cap}}\) regime and C32
occupies the broader moderate- and high-\(T_{\mathrm{cap}}\) regime, with the transition between them controlled by $T_{\mathrm{cap}}$ rather than by $\Delta V_{\mathrm{cap}}$.
Strength is a direct-edge measure and harmonic closeness is a shortest-path measure, so their winners need not agree in a general weighted graph;
that they do agree pointwise here reflects a structural feature of the cislunar orbital network itself. 
Many of the dominant low-cost direct connections also define the admissible minimum-cost paths, so the family that dominates inexpensive direct access also dominates budget-feasible global accessibility at every sampled budget point.
The hub and gateway roles collapse onto the same two-regime organization.


The betweenness map shares the same two-family winner set but differs in geometry.
C32 dominates most of the sampled plane, while C11a appears only in small localized pockets at relatively low \(\Delta V_{\mathrm{cap}}\).
Two main C11a pockets are visible in Figure~\ref{fig:section6_winner_regimes}(c), 
both near $\Delta V_{\mathrm{cap}}$ about 126 m/s, 
one at low $T_{\mathrm{cap}}$
near 8 days, 
and one at high $T_{\mathrm{cap}}$ 
near 25 days, nearly bracketed by tie cells.
The diagonally split cells in Figure~\ref{fig:section6_winner_regimes}(c) mark exact C11a/C32 ties, identifying  the betweenness transition  not as a single clean regime line but as a localized tie-mediated region.
\revone{
The greater occurrence of ties in the betweenness map follows from the
discrete path-count structure of this metric. With 13 families, normalized
betweenness is evaluated over 66 unordered source--target pairs; in the
transition cases examined in Appendix~\ref{app:winner_flips}, newly admitted
routing contributions enter as \(1/66\) terms. Exact equality between
competing relay families can therefore persist across multiple neighboring
sampled budget points. Strength and harmonic closeness, by contrast, depend
on reciprocal edge-cost and shortest-path-cost terms whose magnitudes vary
with the admitted weighted connections. No exact winner ties occur for these
two metrics at the sampled budget points, producing the sharper boundaries
seen in panels (a) and (b).
}
This sensitivity is consistent with the meaning of betweenness:
it measures intermediary routing on admissible minimum-cost paths, which is more sensitive than the direct-edge and shortest-path-magnitude properties measured by strength and harmonic closeness to local changes in the edge set and shortest-path structure. 
The dominant relay family is still C32 over most of the budget plane, but the role is locally contested by C11a in restricted portions of the low-\(\Delta V_{\mathrm{cap}}\) regime.

Taken together, the winner maps reveal a two-family organization of dominant cislunar accessibility roles.
C11a governs the low-time-of-flight accessibility core revealed by strength and harmonic closeness; 
C32 governs the broader moderate- and high-time-of-flight accessibility structure and remains the dominant relay family over most of the budget plane.
The main difference among the centrality measures is therefore not \emph{which} families dominate but \emph{how} the transitions between their roles appear: clean and time-driven for strength and harmonic closeness, localized and tie-mediated for betweenness.
The mechanisms behind these transitions, e.g.,  reweighting of a few dominant reciprocal-access terms in the common strength/harmonic-closeness flip, and  admission of specific source-target pair contributions in the betweenness flips,
are analyzed in Appendix~\ref{app:winner_flips}.


\subsection{Mission Design Implications}\label{subsec:mission_implications}

The cislunar orbital network is best interpreted as a screening tool for accessibility on the common \(C_J=3.1294\) manifold,
playing a role analogous to a family-level $\Delta V$ map for mission-design trade studies in cislunar space:
it identifies which families are favorable starting points for detailed design, which families organize multileg transport, and which families carry persistent accessibility penalties.
It does not replace corrected transfer design, mission-specific targeting, or higher-fidelity modeling.
Table~\ref{tab:mission_design_implications} summarizes how the network's findings map onto representative mission classes.

\par\medskip
\begin{center}
\captionsetup{type=table}
\centering
\small
\caption{Mission-class interpretation of the cislunar orbital network.}
\label{tab:mission_design_implications}
\begin{tabular}{L{0.28\textwidth}L{0.22\textwidth}L{0.42\textwidth}}
\toprule
{\bf Mission class} & {\bf Relevant diagnostic} & {\bf Design implication } \\
\midrule
Rapid response, inspection, rescue, or short-notice retasking
& Direct accessibility and strength
& The (1,1)a-cycler  is the strongest candidate in the low-\(T_{\mathrm{cap}}\) regime; prioritize direct transfers, since multileg paths may require unmodeled phasing, dwell, or loiter time. \\
\addlinespace
Logistics, staging, or redistribution
& Harmonic closeness and multileg closure
& The (3,2)-cycler is the primary screening candidate for broad family-to-family redistribution in moderate- and high-time regimes. \\
\addlinespace
Relay or transfer-handling architectures
& Betweenness
& The (3,2)-cycler frequently lies on admissible minimum-cost routes between other families. \\
\addlinespace
Repeated access to hard-edge families
& Persistent missing direct pairs
& Direct access to the 2:1 stable resonance from L1 or L2 Lyapunov or DPO requires budgets beyond the sampled range; expanded \(\Delta V\) or multileg routes are required. \\
\bottomrule
\end{tabular}
\end{center}
\medskip

The central operational distinction is between direct accessibility and budget-feasible multileg accessibility.
For time-critical missions, a direct family-to-family connection is more important than a multileg route, because the family-level shortest-path model 
aggregates over orbital phase and therefore
does not account for additional phasing, dwell, or loiter time at intermediate families.
Indirect accessibility should not be treated as operationally equivalent to direct accessibility when elapsed time is the limiting resource.

For time-flexible missions, budget-feasible multileg closure becomes the more useful diagnostic, and C32 is the natural staging and redistribution candidate because it combines strong global accessibility with a dominant relay role.
This interpretation is most relevant to logistics redistribution, servicing campaigns, constellation maintenance, and transfer-handling architectures, where repeated movement across the network matters more than minimizing the elapsed time of any single transfer.

The persistent hard edge of the network supplies complementary  guidance of a different kind.
At the most permissive sampled budgets, the only direct pairs still missing are LL1--R21-S, LL2--R21-S, and DPO--R21-S.
This does not imply that R21-S is unreachable;
it indicates that direct access involving R21-S, particularly from LL1, LL2, or DPO, lies outside the maneuver--time budget envelope sampled here.
Mission classes requiring repeated direct access to R21-S would therefore need either expanded $\Delta V$ budgets or multileg routes, both with the time accounting caveats noted above.

\section{Validation through Trajectory Realization}\label{sec:validation}

The reduced family-to-family accessibility description developed in Sections~\ref{sec:overlap}--\ref{sec:budget} rests on overlap-based proxy costs. 
Whether those proxy values correspond to physically realizable trajectories has not yet been tested, and is the question this section addresses through selective trajectory realization rather than exhaustive correction of the full family-pair matrix. 
Selected proxy-supported direct and relay connections are converted into concrete patched trajectories through local differential correction, and the corrected maneuver costs are compared with the proxy values that identified them.
The corrected costs fall below the proxy values in every example examined, supporting the use of the overlap-based construction as a conservative screening measure rather than a final transfer estimate.


\subsection{From Proxy Overlaps to Corrected Trajectory Realizations}\label{subsec:from_proxy}

The pairwise overlap construction of Section~\ref{sec:overlap} assigns each family pair a representative voxel and a proxy cost; that proxy is a screening quantity derived from local family-to-family proximity, not a realized transfer trajectory.
To turn a selected proxy entry into a concrete trajectory, the minimizing overlap voxel is revisited and the lowest-turning-cost source-side and target-side arcs reaching that voxel are extracted. The two extracted arcs do not in general satisfy exact continuity at the voxel midpoint; a local differential correction is then applied to adjust the phase, heading, and time along each arc so that the two meet at a common patch point in the reduced state variables $(x, y, \theta)$ while locally minimizing the total turning effort. 
The full correction formulation (e.g., decision vector, normalized residual, and convergence criterion) is given in Appendix~\ref{app:correction}.

The procedure produces three distinct cost levels for each selected family pair: the \emph{matrix proxy value} from the overlap construction, the \emph{warm-start} cost of the two extracted arcs before correction, and the \emph{corrected cost} after local patching. The direct examples in Section~\ref{subsec:direct_examples} compare all three.


\subsection{Corrected Direct Trajectory Examples from the (1,1)a-Cycler}\label{subsec:direct_examples}

To illustrate the trajectory realization step, two corrected direct examples were generated from the representative source family, C11a.
The selected targets were LL1 and R31-S.
Together, these examples span a familiar libration-family target and a stable resonant target.

For each case, the minimum-proxy overlap voxel from Section~\ref{sec:overlap} was used to extract a concrete warm start trajectory pair.
A local differential correction was then applied to enforce continuity at the patch point in \((x,y,\theta)\) according to the normalized convergence criterion introduced in Section~\ref{subsec:from_proxy}.
Both selected cases satisfied this criterion and therefore provide realized direct trajectories associated with the proxy-identified family connections.

Figure~\ref{fig:section7_direct_examples} shows the two corrected direct trajectories together with the corresponding source and target family curves on the cislunar background.
The common source family C11a appears in dashed red in each panel, while the dashed blue curve identifies the target family.
The solid trajectory segments show the corrected source-side and target-side legs that meet at the patch point.
The LL1 example provides a low-cost realization toward a familiar libration family, whereas the R31-S example demonstrates that the same proxy-based realization procedure can reach a substantially more expensive stable resonant target with distinctly different trajectory geometry.

\par\medskip
\begin{center}
\captionsetup{type=figure}
    \centering

    \begin{subfigure}[t]{0.48\textwidth}
        \centering
        \includegraphics[width=\textwidth]{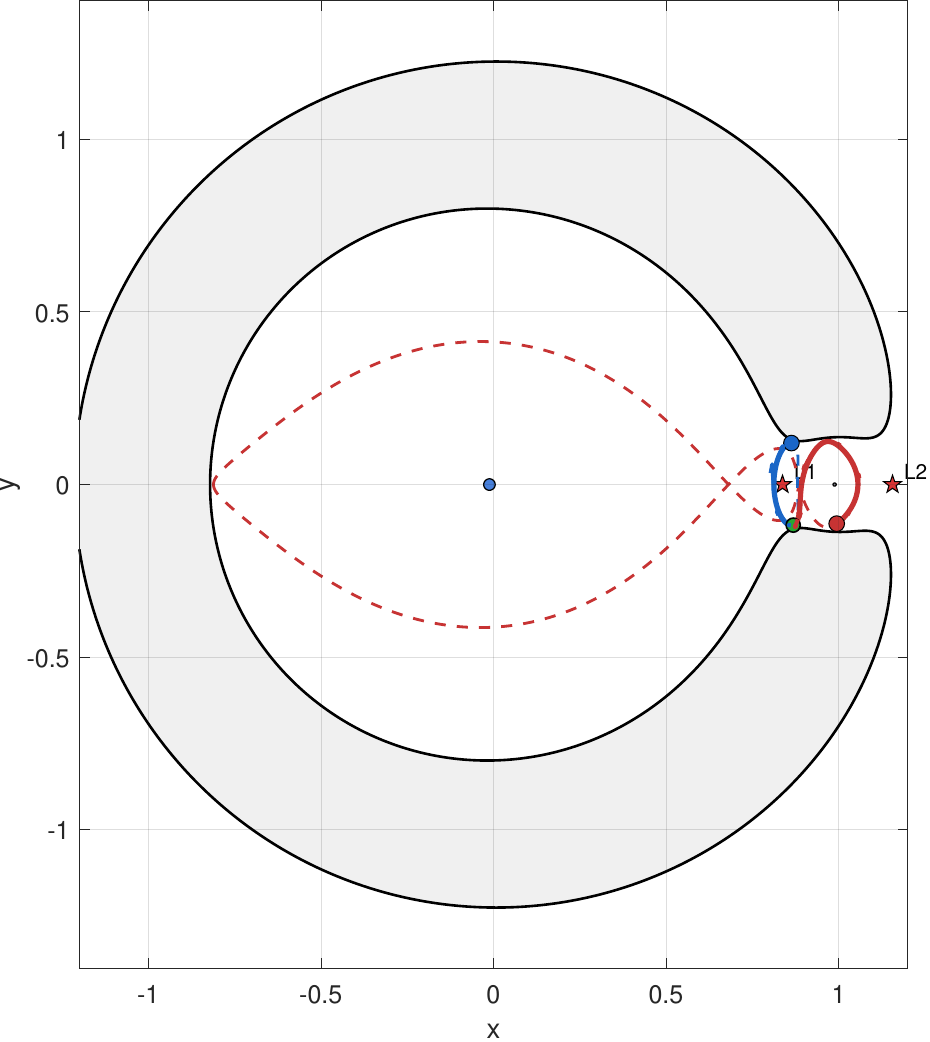}
        \caption{\footnotesize Corrected direct transfer from C11a to LL1.}
        \label{fig:section7_direct_L1}
    \end{subfigure}
    \hfill
    \begin{subfigure}[t]{0.48\textwidth}
        \centering
        \includegraphics[width=\textwidth]{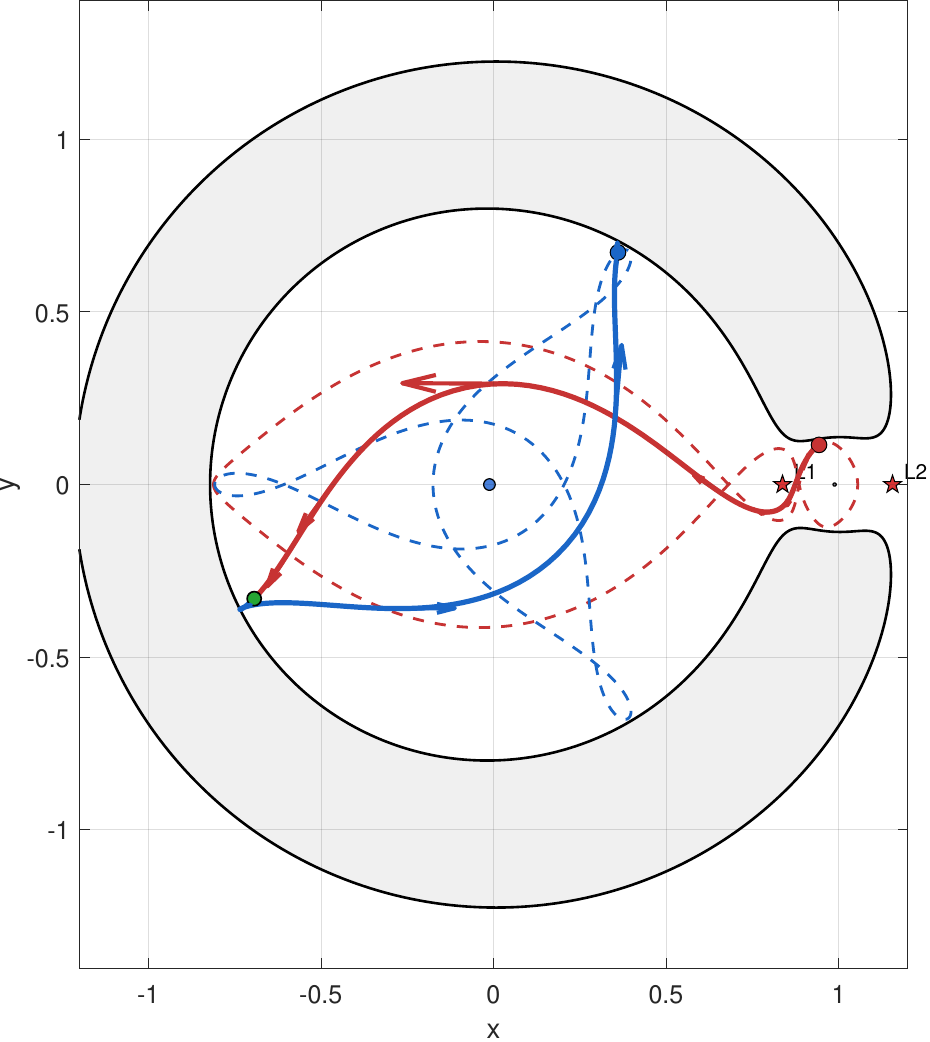}
        \caption{\footnotesize Corrected direct transfer from C11a to R31-S.}
        \label{fig:section7_direct_R31S}
    \end{subfigure}

    \caption{\footnotesize Selected corrected direct trajectories from (1,1)a-cycler (C11a) at the common Jacobi level \(C_J=3.1294\).
    The two examples span a libration-family target and a stable resonant target.
    In each panel, the dashed red curve denotes the source family C11a, the dashed blue curve denotes the target family, and the solid red and blue arcs denote the corrected source-side and target-side trajectory segments, respectively.
    Together, these examples show that proxy-supported family connections can be turned into exact patched direct realizations across geometrically distinct target-family types.}
    \label{fig:section7_direct_examples}
\end{center}
\medskip

The quantitative comparison for the examples is summarized in Table~\ref{tab:section7_direct_examples}.
For each case, the table reports the matrix proxy cost, the extracted warm-start cost, and the corrected trajectory cost.
This three-level comparison is important because the matrix entry is only a proxy screening quantity, whereas the extracted and corrected values are associated with a specific realized candidate.

\par\medskip
\begin{center}
\captionsetup{type=table}
\centering
\caption{\footnotesize Selected corrected direct trajectories from the (1,1)a-cycler.
All maneuver quantities are reported in \(\mathrm{m/s}\), and time of flight is reported in days.
For each target family, the table reports the matrix proxy cost, the extracted warm-start cost, the corrected trajectory cost, the absolute reduction from proxy to corrected cost, and the corrected total time of flight.}
\label{tab:section7_direct_examples}
\begin{tabular}{lccccc}
\hline
Target family & Proxy \(\Delta V\) & Extracted \(\Delta V\) & Corrected \(\Delta V\) & Reduction & TOF \\
 & \multicolumn{4}{c}{\(\mathrm{[m/s]}\)} & [days] \\
\hline
L1 Lyapunov& 4.02 & 3.44 & 2.25 & 1.77 & 16.99 \\
3:1 stable resonant & 178.31 & 178.01 & 174.85 & 3.46 & 24.69 \\
\hline
\end{tabular}
\end{center}
\medskip

In both cases, the corrected \(\Delta V\) is lower than the matrix proxy value.
This is the desired direction for a screening metric: the overlap-based proxy behaves conservatively for these examples, identifying candidate connections without underestimating the final corrected maneuver cost.
At the same time, the correction does not erase the family-level cost distinctions identified by the proxy matrix.
The LL1 transfer remains a low-cost connection, while the transfer to R31-S remains substantially more expensive.
Thus, the corrected trajectories preserve the broad proxy-based cost contrast between the selected examples.

These results support the intended role of the proxy construction in the present paper.
The proxy is not used as a rigorous bound or as a final transfer estimate.
Instead, it is used as a reduced accessibility measure derived from local family-to-family proximity in the common admissible domain.
The corrected direct examples show that selected proxy-supported connections can be turned into concrete trajectories with maneuver-cost scales consistent with the accessibility structure identified by the screening model.

The role of this subsection is therefore validation by representative realization rather than exhaustive correction of the direct family-pair matrix.
The examples demonstrate that the reduced accessibility model is capable of identifying physically meaningful direct candidates across both conventional and resonant target families.
This provides the direct-trajectory foundation for the relay-path interpretation developed in Section~\ref{subsec:relay_examples}.


\subsection{Relay-Path Realization and Interpretation}\label{subsec:relay_examples}

The same correction procedure extends naturally to admissible minimum-cost multileg routes identified by the network analysis of Sections~\ref{sec:network}--\ref{sec:budget}.
Two representative relay paths through C11a are realized here as concrete two-leg trajectories: R31-U \(\rightarrow\) LL2 and R31-S \(\rightarrow\) LL2.
In each case, the differential correction of Appendix~\ref{app:correction} is applied to both legs separately, producing two patched arcs that meet C11a at common bridge phases.
Figure~\ref{fig:section7_relay_examples} shows the two resulting relay trajectories, and Table~\ref{tab:section7_relay_examples} compares the relay-path cost with the direct-pair proxy cost for each endpoint pair.

\par\medskip
\begin{center}
\captionsetup{type=figure}
    \centering

    \begin{subfigure}[t]{0.48\textwidth}
        \centering
        \includegraphics[width=\textwidth]{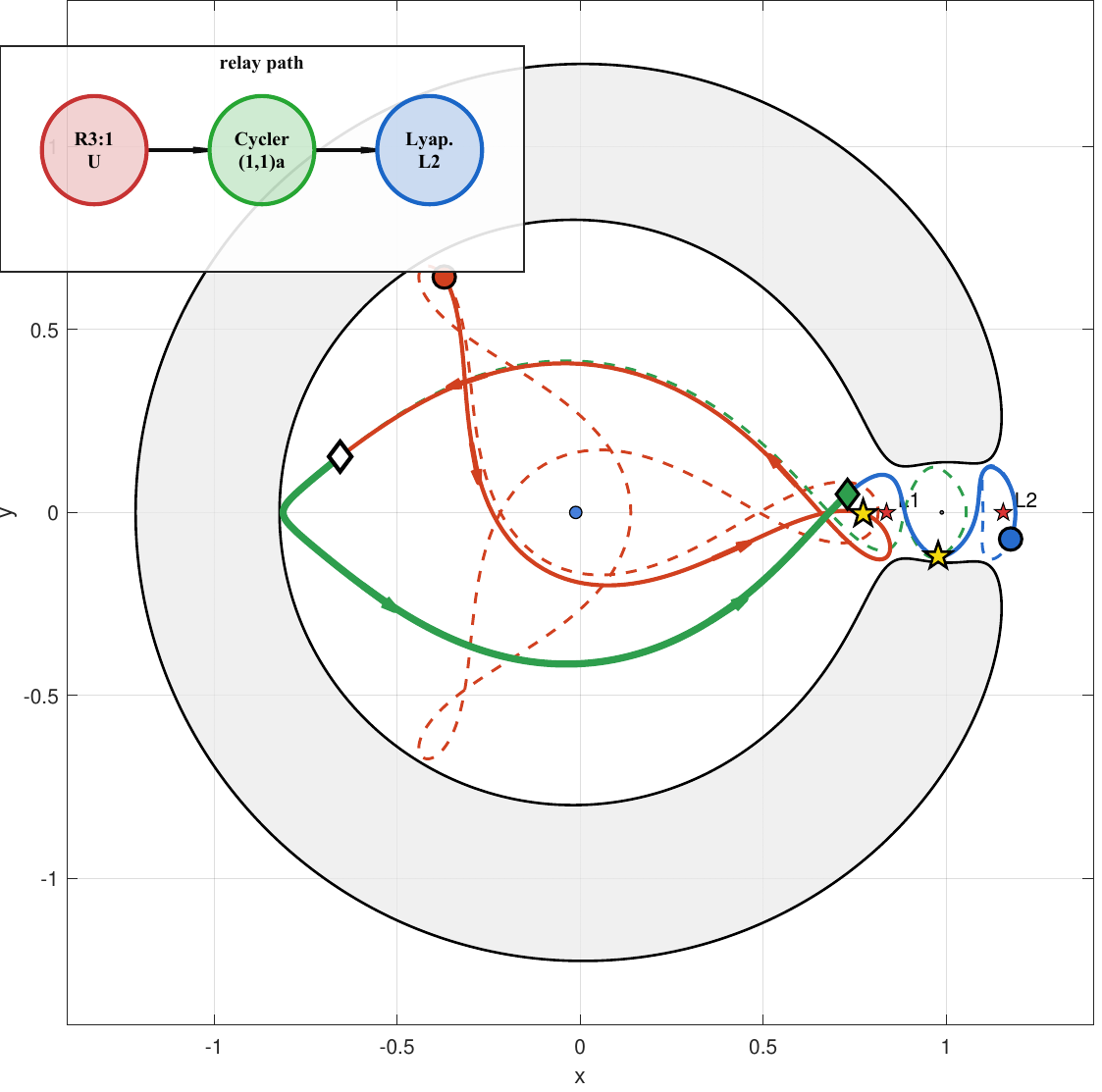}
        \caption{\footnotesize \(3{:}1\) unstable resonant to Lyapunov~\(L_2\) through (1,1)a-cycler. A video version is available at \url{https://youtu.be/_UUlLCnePGE}.}
        \label{fig:section7_relay_ex1}
    \end{subfigure}
    \hfill
    \begin{subfigure}[t]{0.48\textwidth}
        \centering
        \includegraphics[width=\textwidth]{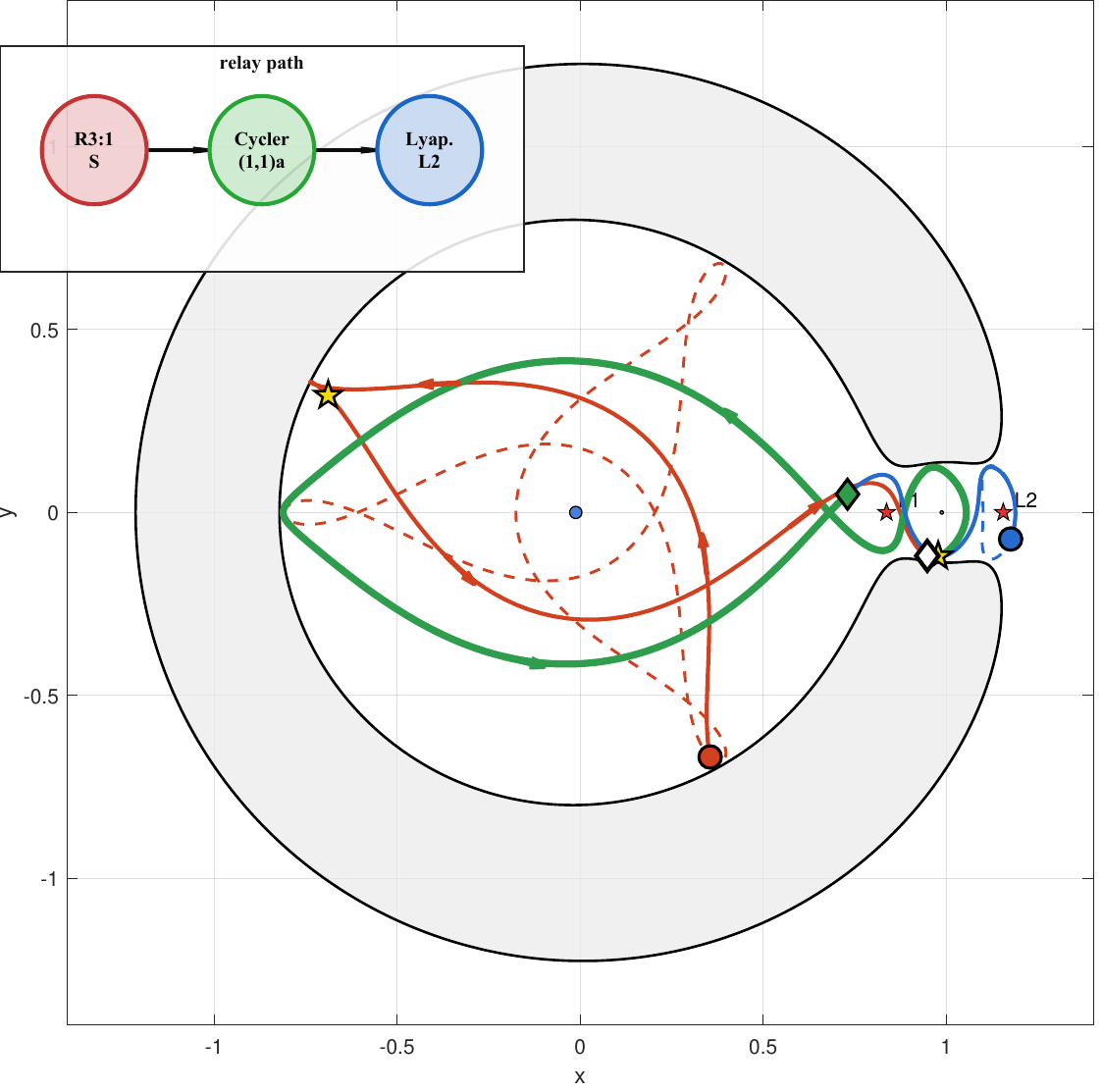}
        \caption{\footnotesize \(3{:}1\) stable resonant to Lyapunov~\(L_2\) through (1,1)a-cycler. A video version is available at \url{https://youtu.be/aJBk93Ts6DQ}.}
        \label{fig:section7_relay_ex3}
    \end{subfigure}

    \caption{\footnotesize Representative relay-path illustrations associated with admissible minimum-cost multileg routes in the reduced family network.
    Each panel combines the trajectory-space realization with an inset three-node schematic of the corresponding relay structure.
    In the trajectory plots, the dashed red curve denotes the origin family, the dashed green curve denotes the bridge family (1,1)a-cycler, and the dashed blue curve denotes the target family.
    The solid arcs show the first transfer leg, bridge-family coast segment, and second transfer leg.
    These examples physically realize the betweenness interpretation by showing how the (1,1)a-cycler acts as a dynamically accessible bridge between endpoint families.}
    \label{fig:section7_relay_examples}
\end{center}
\medskip


\par\medskip
\begin{center}
\captionsetup{type=table}
\centering
\caption{Representative relay-path examples through (1,1)a-cycler.
All maneuver quantities are reported in \(\mathrm{m/s}\).
For each case, the table compares the relay-path cost with the corresponding direct proxy cost and reports the implied savings.}
\label{tab:section7_relay_examples}
\begin{tabular}{lccc}
\hline
Endpoint pair & Relay \(\Delta V\) & Direct proxy \(\Delta V\) & Savings \\
 & \multicolumn{2}{c}{\(\mathrm{m/s}\)} & \(\mathrm{m/s}\) (\%) \\
\hline
\(3{:}1\) unstable resonant \(\rightarrow\) \(L_2\) Lyapunov & 62 & 181 & 119 (66\%) \\
\(3{:}1\) stable resonant \(\rightarrow\) \(L_2\) Lyapunov& 180 & 244 & 64 (26\%) \\
\hline
\end{tabular}
\end{center}
\medskip

These examples give direct geometric meaning to the intermediary role identified by the shortest-path and betweenness analyses of Sections~\ref{sec:network}--\ref{sec:budget}.
The C11a is not merely a high-betweenness graph node, but a dynamically accessible bridge through which one leg approaches the intermediary family and a second leg departs toward the final target.
The two examples show how an intermediate family can link endpoint families that would otherwise be connected only through a more expensive direct proxy route.

The savings are substantial in both cases.
For the R31-U--LL2 endpoint pair, the proxy cost decreases from 181~m/s to 62~m/s, a 66\% reduction.
For the harder-to-access R31-S--LL2 pair, it decreases from 244~m/s to 180~m/s, a 26\% reduction.
The relay-path cost is lower than the direct-pair proxy in both cases, supporting the network-level finding that multileg structure is a real source of accessibility savings rather than a graph-theoretic curiosity.

These relay examples should be interpreted more cautiously than the corrected direct trajectories of Section~\ref{subsec:direct_examples}.
The figures show patched two-leg trajectories, but the network analysis does not impose cumulative multileg time, so the bridge-family coast segment between the two legs is not phase-resolved.
Their purpose here is to illustrate the physical meaning of admissible multileg proxy routes identified by the reduced family network and to clarify why an intermediary family acquires routing importance, not to claim fully phase-corrected multileg transfers.

\section{Conclusion}\label{sec:conclusion}

This paper introduced a reachable-set-based framework for characterizing family-to-family accessibility in cislunar space on a common planar Earth--Moon CR3BP Jacobi manifold.
For each representative family, finite-\(\Delta V\), finite-time reachable-set atlases on a reduced \((x, y, \theta)\) phase space were constructed from energy-preserving heading-change maneuvers, with backward atlases obtained from forward ones through time-reversal symmetry.
Pairwise overlap on the common voxel grid produced proxy accessibility costs that were then assembled into the cislunar orbital network, a named weighted graph on the representative family set.

Three structural findings emerge.
First, at the maximum-budget reference case, the (3,2)-cycler is the dominant family in the network, ranking first in strength, harmonic closeness, and betweenness, and therefore simultaneously the strongest direct-access, most efficient global-access, and dominant relay family.
The (1,1)a-cycler dominates the same roles in the low-time-budget regime, and the 2:1 stable resonant orbit remains the persistent hard-access family across the entire sampled budget plane.
Second, direct accessibility, graph connectedness, and budget-feasible multileg closure separate into three distinct regimes across the budget plane, showing that network connectedness and operationally admissible multileg reachability are not interchangeable.
Third, selected proxy-supported family connections were converted into corrected patched trajectories whose maneuver costs lie below the proxy values in every tested case, supporting the use of the overlap-based construction as a conservative screening measure rather than a final transfer estimate.

Together, these findings have a natural geometric reading on the Jacobi manifold. The ten unstable representatives sit in what we hypothesize is a single connected chaotic region of the $C_J = 3.1294$ manifold, providing the low-cost connectivity medium that the network's dominant hub, gateway, and relay roles exploit. 
The reachable-set construction enters this medium through a restricted maneuver model: every admissible impulse is a pure rotation of the rotating-frame velocity vector at fixed rotating-frame speed, which in many configurations corresponds to a substantial normal component of the inertial impulse rather than the tangential burns that two-body intuition associates with efficient orbit transfer. The framework therefore shows that purely directional, fixed-rotating-speed control is sufficient to achieve low-$\Delta V$ family-to-family accessibility across the unstable representatives, with the chaotic-region geometry providing the medium that directional control alone can exploit. 
The three stable resonant representatives sit inside resonant tori that are dynamically isolated from the chaotic sea, and therefore occupy the persistent hard-access edge of the accessibility hierarchy. 
The reachable-set overlap construction does not independently verify connectivity of the chaotic region, but the dense network of low-cost admissible direct connections it identifies among the unstable representatives is consistent with this picture, and the network of admissible budget-feasible multileg paths sharpens it.

The present results should be interpreted within the scope of the model.
The analysis is planar, restricted to a single common Jacobi manifold, and based on representative periodic orbits rather than full continued families.
\revtwo{
The network structure obtained here should therefore be interpreted as
specific to the common fixed value \(C_J=3.1294\), rather than as an
energy-independent hierarchy. If the analysis were repeated at a different
common fixed \(C_J\), the representative periodic orbits would change in
geometry and stability, some branches could approach folds or bifurcations,
and changes in the zero-velocity topology could alter the accessible
transport regions. These effects can modify the reachable-set overlaps,
pairwise proxy costs, direct connectivity, and the resulting hub, gateway,
and relay roles. In particular, decreasing \(C_J\) can enlarge the
energetically accessible region, but the resulting changes in
family-to-family accessibility need not be monotonic because the orbit
families and their local transport structures vary simultaneously with
energy. The present thirteen-family network should therefore be regarded as
one common-energy slice of a broader energy-dependent orbital network.
}
The maneuver model uses bounded fixed-\(C_J\) heading changes, the proxy costs are screening quantities rather than optimized transfer costs, and the family-level network aggregates over orbital phase, so the bridge-family coast time on multileg routes is not phase-resolved.
The trajectory-realization study validates selected representative connections but does not constitute an exhaustive correction of the full family-pair matrix.

Natural extensions include carrying the framework into the spatial CR3BP and higher-fidelity ephemeris dynamics; verifying the single-connected-chaotic-region hypothesis; extending the maneuver model beyond energy-preserving heading changes to general impulsive maneuvers, which lifts the cislunar orbital network from a single Jacobi manifold to a multi-manifold network indexed by both family identity and energy; developing phase-resolved or time-varying family networks that close the multileg-time-accounting gap; and incorporating low-thrust or optimal-control refinement. In such settings, the reachable-set atlas and the cislunar orbital network could serve as a screening and initialization layer for higher-fidelity transfer design.

Within its present scope, the framework provides a principled and operationally interpretable family-level $\Delta V$ map of cislunar family-to-family accessibility.
Its value is not that it replaces detailed trajectory design, but that it exposes network-level structure that isolated pairwise transfers do not reveal: dominant access families, persistent hard-access families, distinct budget-dependent accessibility regimes, and relay pathways that organize indirect transport as budgets vary.

\section*{Appendices}
\appendix
\section{Maximum-Budget Network Validation}
\label{app:baseline_validation}

The maximum-budget cislunar orbital network was validated against a 15-case numerical suite spanning systematic coarsening, refinement, and stress-envelope configurations around the maximum-budget reference case.
For readability, each validation case is assigned an identifier \(C01\)--\(C15\).
Table~\ref{tab:validation_case_configs} defines these identifiers and lists the numerical settings used in each case.

\par\medskip
\begin{center}
\captionsetup{type=table}
\centering
\small
\caption{Validation cases used to assess the robustness of the maximum-budget cislunar orbital network.
The case identifiers \(C01\)--\(C15\) are used in Figs.~\ref{fig:validation_case_rank_heatmaps} and \ref{fig:validation_case_summary}.}
\label{tab:validation_case_configs}
\begin{tabular}{clcccc}
\toprule
Case ID & Validation case description & \(dx=dy\) & \(\Delta\theta\) & \(\Delta\delta_{\mathrm{fan}}\) & \(\Delta s_{\mathrm{seed}}\) \\
\midrule
C01 & maximum-budget reference case & 0.0010 & \(1.0^\circ\) & \(0.5^\circ\) & 0.010 \\
C02 & Spatial coarsening & 0.0015 & \(1.0^\circ\) & \(0.5^\circ\) & 0.010 \\
C03 & Heading-grid coarsening & 0.0010 & \(1.5^\circ\) & \(0.5^\circ\) & 0.010 \\
C04 & Heading-fan coarsening & 0.0010 & \(1.0^\circ\) & \(1.0^\circ\) & 0.010 \\
C05 & Seed-spacing coarsening & 0.0010 & \(1.0^\circ\) & \(0.5^\circ\) & 0.020 \\
C06 & Coupled spatial--heading coarsening & 0.0015 & \(1.5^\circ\) & \(0.5^\circ\) & 0.010 \\
C07 & Fine reference & 0.0005 & \(0.5^\circ\) & \(0.5^\circ\) & 0.010 \\
C08 & Fine reference with nominal spatial grid & 0.0010 & \(0.5^\circ\) & \(0.5^\circ\) & 0.010 \\
C09 & Fine reference with nominal heading grid & 0.0005 & \(1.0^\circ\) & \(0.5^\circ\) & 0.010 \\
C10 & Fine reference with coarse heading fan & 0.0005 & \(0.5^\circ\) & \(1.0^\circ\) & 0.010 \\
C11 & Fine reference with coarse seed spacing & 0.0005 & \(0.5^\circ\) & \(0.5^\circ\) & 0.020 \\
C12 & Fine reference with composite coarsening & 0.0010 & \(1.0^\circ\) & \(1.0^\circ\) & 0.020 \\
C13 & Stress-coarse envelope & 0.0020 & \(2.0^\circ\) & \(2.0^\circ\) & 0.030 \\
C14 & Very-fine reference & 0.00025 & \(0.25^\circ\) & \(0.25^\circ\) & 0.005 \\
C15 & Intermediate stress-coarse envelope & 0.0015 & \(1.5^\circ\) & \(1.5^\circ\) & 0.030 \\
\bottomrule
\end{tabular}
\end{center}
\medskip

The first validation result is topological.
As shown in Fig.~\ref{fig:validation_case_summary}(a,b), every validation case yields the same direct-network size, namely 75 finite family pairs out of the 78 possible undirected pairs, and all 78 family pairs remain connected by an admissible shortest path under the total maneuver budget.
The same three direct pairs are absent in every case, namely LL1--R21-S, LL2--R21-S, and DPO--R21-S.
Thus, the direct topology of the maximum-budget cislunar orbital network is not a fragile consequence of one particular numerical configuration.

The second validation result is quantitative.
As shown in Fig.~\ref{fig:validation_case_summary}(c,d), the mean admissible shortest-path proxy cost varies across the suite from \(62.1~\mathrm{m/s}\) to \(73.5~\mathrm{m/s}\), while the mean direct-edge proxy time varies from \(17.00~\mathrm{days}\) to \(20.02~\mathrm{days}\).
This variation is expected because the direct proxy costs are constructed from sampled reachable-set atlases and finite voxel overlap.
Changes in spatial resolution $(\Delta x,\Delta y)$, heading discretization \(\Delta\theta\), fan sampling increment \(\Delta\delta_{\mathrm{fan}}\), and seed spacing \(\Delta s_{\mathrm{seed}}\) modify the sampled candidate trajectories, the resolved overlap structure, and the voxel-scale patching contribution.
The numerical suite is therefore intended to test qualitative robustness of the network interpretation rather than invariance of every proxy value.

\par\medskip
\begin{center}
\captionsetup{type=figure}
    \centering
    \begin{subfigure}[t]{0.62\textwidth}
        \centering
        \includegraphics[width=\textwidth]{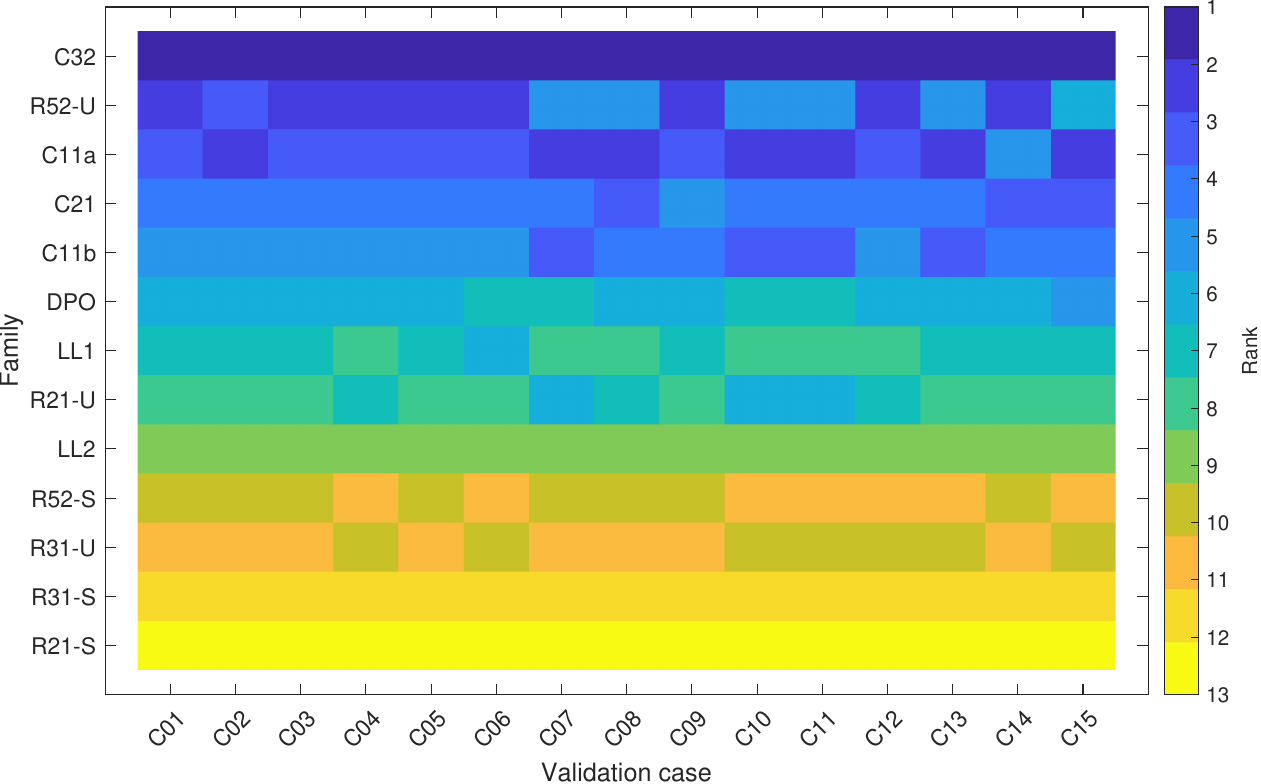}
        \caption{\footnotesize Strength rank.}
    \end{subfigure}

    \vspace{0.25em}

    \begin{subfigure}[t]{0.47\textwidth}
        \centering
        \includegraphics[width=\textwidth]{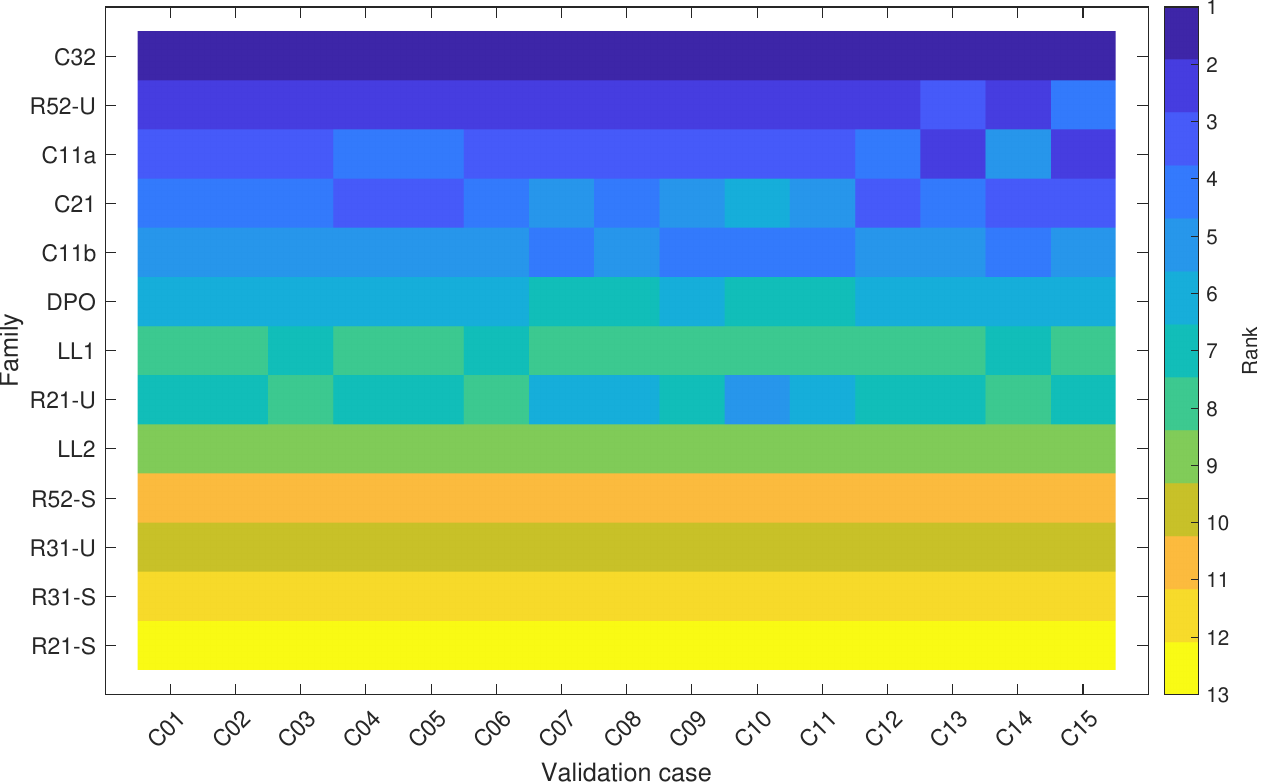}
        \caption{\footnotesize Harmonic-closeness rank.}
    \end{subfigure}
    \hfill
    \begin{subfigure}[t]{0.47\textwidth}
        \centering
        \includegraphics[width=\textwidth]{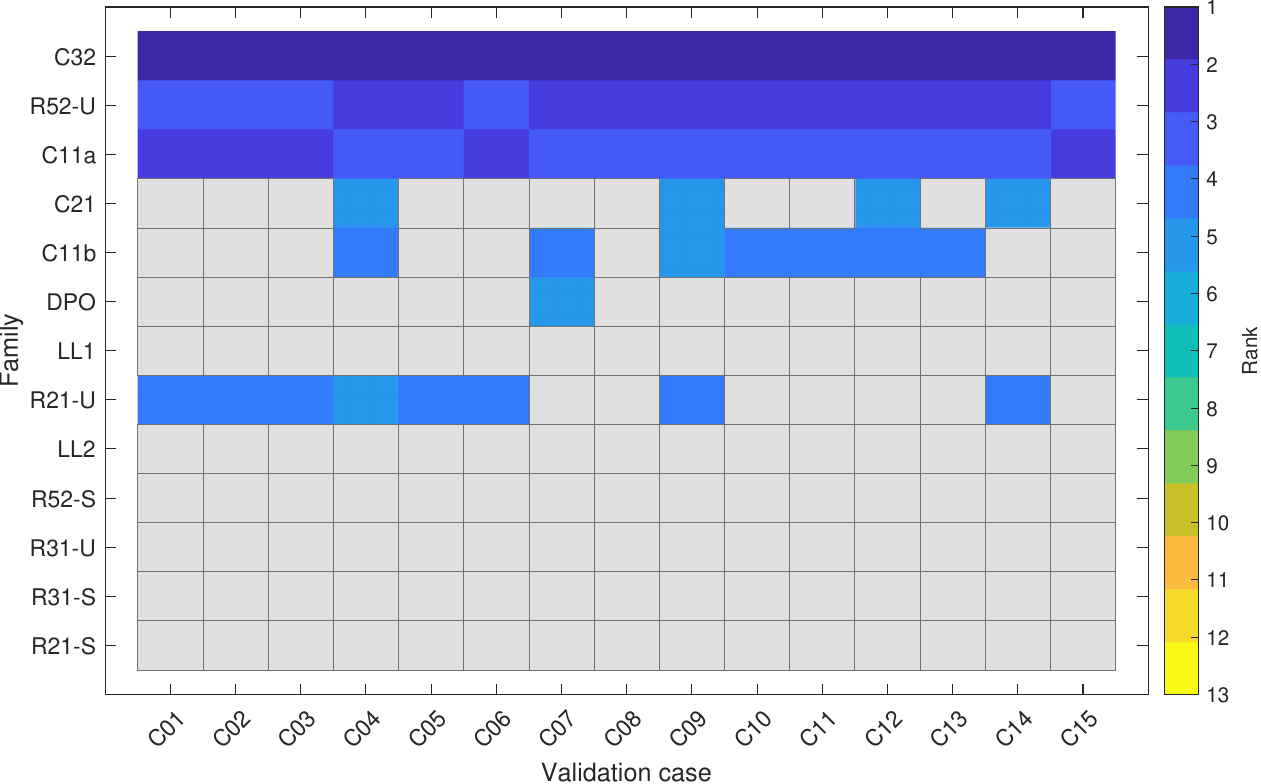}
        \caption{\footnotesize Betweenness rank.}
    \end{subfigure}
    \caption{\footnotesize Family-rank robustness across the validation suite.
    The horizontal axis uses the case identifiers \(C01\)--\(C15\) defined in Table~\ref{tab:validation_case_configs}.
    Families are ordered by the maximum-budget strength ranking from Section~\ref{subsec:centralitiesrank}, and lower rank is better.}
    \label{fig:validation_case_rank_heatmaps}
\end{center}
\medskip

The third validation result is the robustness of the family-role hierarchy.
Figure~\ref{fig:validation_case_rank_heatmaps}(a)--(c) shows the strength, harmonic-closeness, and betweenness ranks across the 15 validation cases, with family ordering fixed by the maximum-budget strength ranking from Section~\ref{subsec:centralitiesrank}.
The dominant feature is the complete persistence of C32, which remains rank~1 in all three measures for all 15 cases.
Accordingly, the maximum-budget interpretation of C32 as the dominant direct-access, gateway, and relay family is numerically stable across the entire validation suite.

Below C32, the same second tier also persists.
R52-U and C11a remain the two strongest competitors in the weighted network, although their relative ordering depends on the chosen measure.
R52-U remains especially strong in strength and harmonic closeness, whereas C11a remains consistently among the strongest relay families.
At the opposite extreme, R21-S remains last in strength and harmonic closeness in every validation case and retains zero betweenness throughout.
The cislunar orbital network therefore preserves the same hard-access family at the difficult end of the cislunar family set.

Taken together, Figs.~\ref{fig:validation_case_rank_heatmaps} and \ref{fig:validation_case_summary} show that the maximum-budget cislunar orbital network is quantitatively sensitive but qualitatively robust under the present parameter variations.

\par\medskip
\begin{center}
\captionsetup{type=figure}
    \centering
    \begin{subfigure}[t]{0.47\textwidth}
        \centering
        \includegraphics[width=\textwidth]{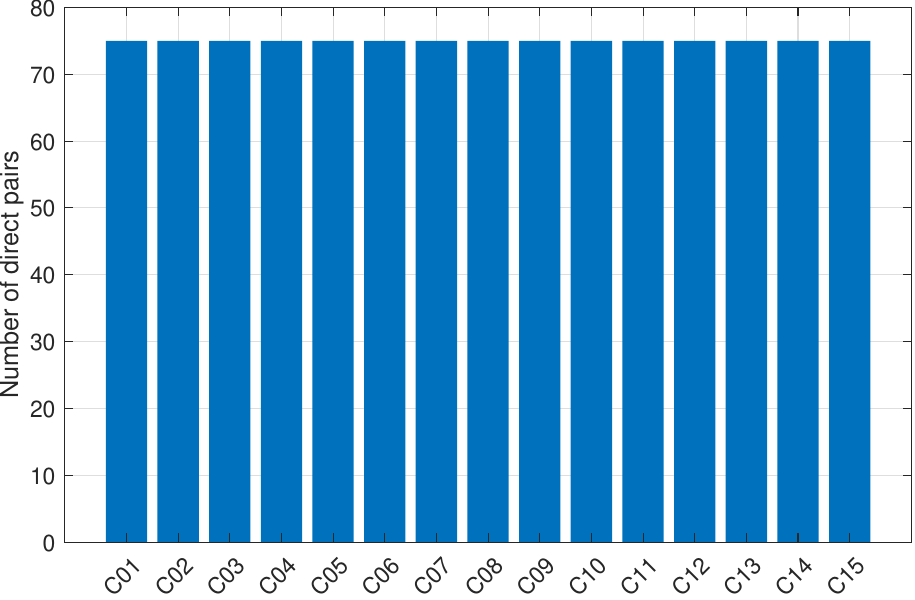}
        \caption{\footnotesize Finite direct pairs.}
    \end{subfigure}
    \hfill
    \begin{subfigure}[t]{0.47\textwidth}
        \centering
        \includegraphics[width=\textwidth]{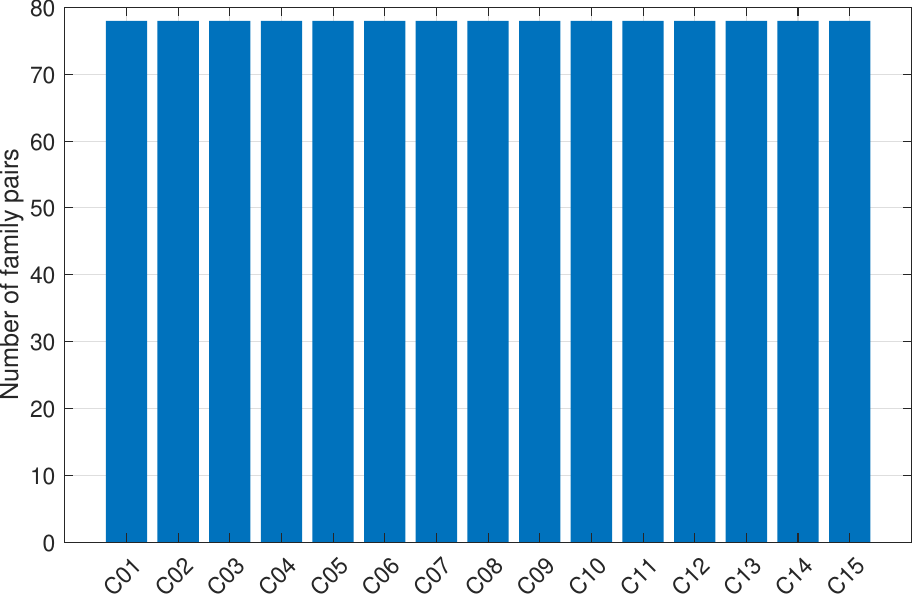}
        \caption{\footnotesize Budget-feasible shortest-path pairs.}
    \end{subfigure}

    \vspace{0.25em}

    \begin{subfigure}[t]{0.47\textwidth}
        \centering
        \includegraphics[width=\textwidth]{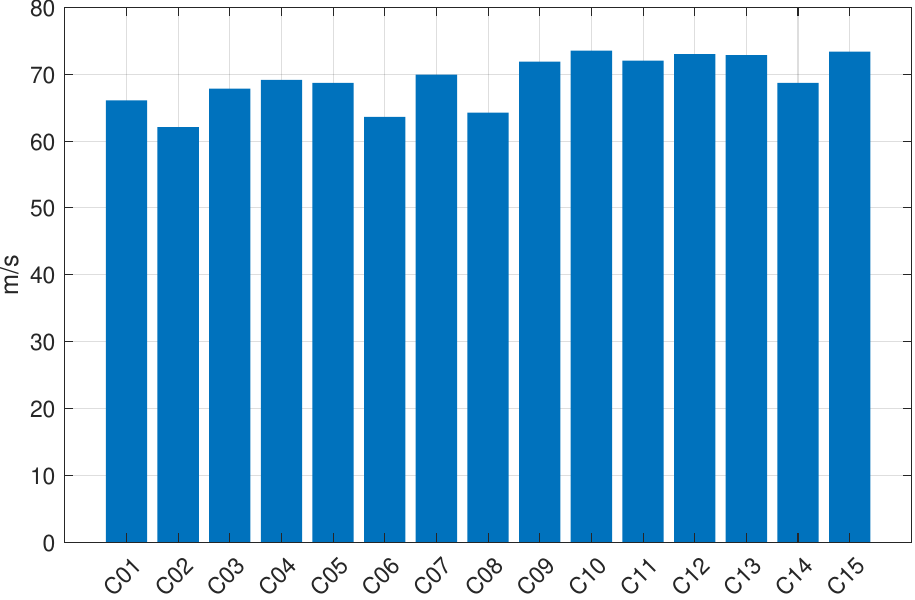}
        \caption{\footnotesize Mean admissible shortest-path proxy cost.}
    \end{subfigure}
    \hfill
    \begin{subfigure}[t]{0.47\textwidth}
        \centering
        \includegraphics[width=\textwidth]{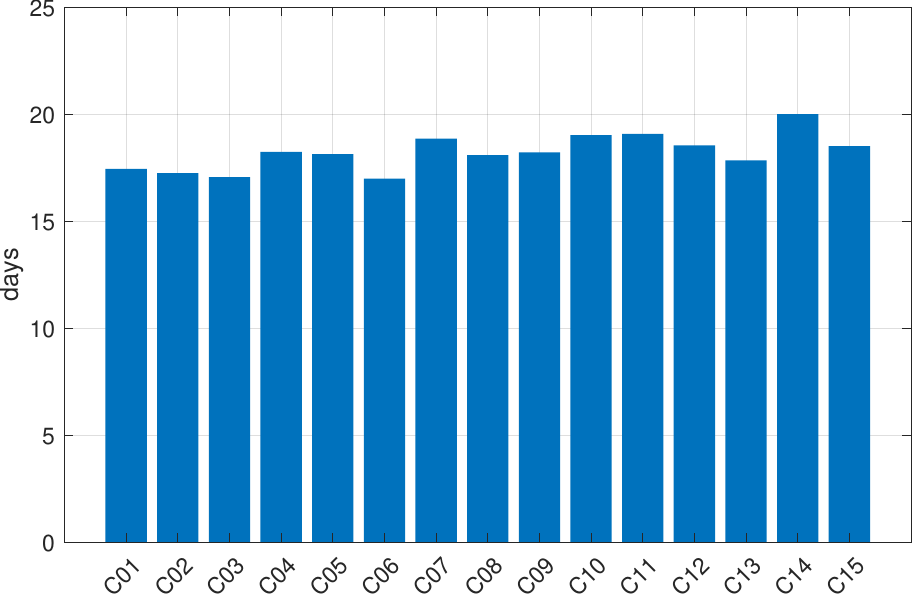}
        \caption{\footnotesize Mean direct-edge proxy time.}
    \end{subfigure}
    \caption{\footnotesize Case-level summary of the validation suite.
    The horizontal axis uses the case identifiers \(C01\)--\(C15\) defined in Table~\ref{tab:validation_case_configs}.
    The top row reports the invariant connectivity quantities, and the bottom row shows the remaining variation in weighted cost and time across the validation cases.}
    \label{fig:validation_case_summary}
\end{center}
\medskip


\section{Diagnostic Analysis of Winner Flips}
\label{app:winner_flips}

Section~\ref{subsec:budget_winners} identified where the dominant family changes over the budget plane.
This appendix diagnoses why those changes occur.
The analysis is carried out directly at the level of the metric terms.
For strength and harmonic closeness, the relevant quantities are reciprocal direct-cost and reciprocal shortest-path terms.
For betweenness, the relevant quantities are the individual source--target pair contributions.
Across the transition cases analyzed here, every winner-level tie occurs only between C11a and C32.

The common strength and harmonic-closeness flip is best illustrated at fixed \(\Delta V_{\mathrm{cap}}=230.21~\mathrm{m/s}\), where the winner changes from C11a to C32 as \(T_{\mathrm{cap}}\) increases from \(10.25\) to \(11.10~\mathrm{days}\).
This transition is not a new connectedness or closure regime.
The largest connected component remains unchanged, and the shortest-path-feasible pair count is also unchanged.
Instead, the flip is caused by a strong reweighting of a small number of reciprocal-cost and reciprocal-distance terms, as shown in Figure~\ref{fig:section6_strength_harmonic_term_deltas}.

\par\medskip
\begin{center}
\captionsetup{type=figure}
    \centering
    \includegraphics[width=0.85\textwidth]{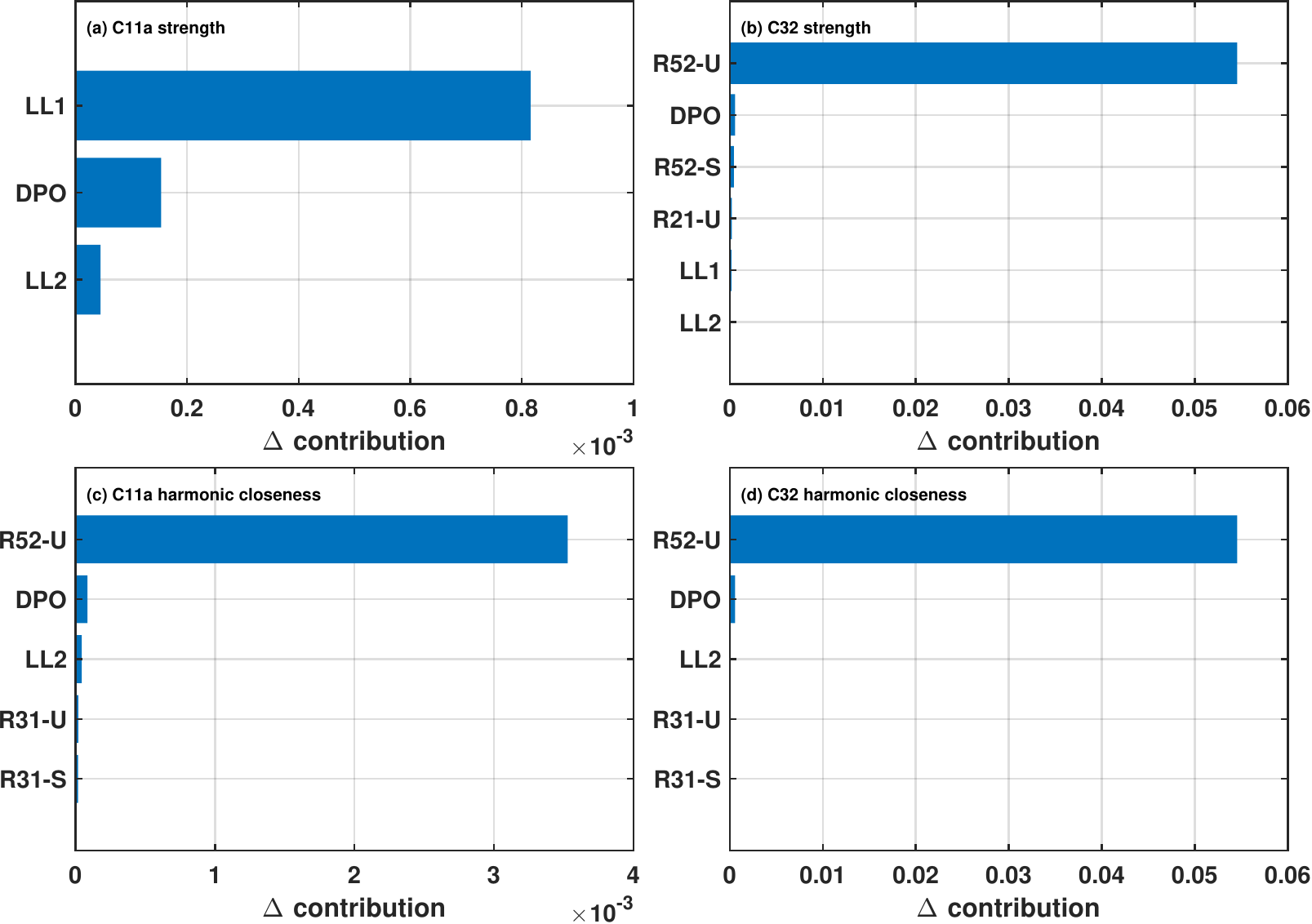}
    \caption{\footnotesize Term-level explanation of the common strength and harmonic-closeness winner flip at fixed \(\Delta V_{\mathrm{cap}}=230.21~\mathrm{m/s}\) as \(T_{\mathrm{cap}}\) increases from \(10.25\) to \(11.10~\mathrm{days}\).
    Panels (a) and (b) show the largest positive changes in the reciprocal direct-cost terms contributing to the strength of C11a and C32, respectively.
    Panels (c) and (d) show the largest positive changes in the reciprocal shortest-path terms contributing to the harmonic closeness of the same two families.
    In both metrics, the dominant gain is the R52-U term for C32, showing that the winner flip is caused primarily by the reweighting of existing access terms rather than by the mere appearance of a new direct edge.}
    \label{fig:section6_strength_harmonic_term_deltas}
\end{center}
\medskip

The dominant change in both metrics is associated with C32's access to R52-U.
That single term grows far more than any corresponding term for C11a, and it accounts for most of the increase in both the strength and harmonic-closeness scores of C32.
By contrast, the newly admitted direct transfer to R52-S contributes only weakly.
The strength and harmonic-closeness winner flip is therefore not driven by edge appearance alone.
It is driven by the disproportionate amplification of already important reciprocal-access terms for C32.

The betweenness transitions require a different interpretation.
With 13 families, normalized betweenness is evaluated over 66 unordered source--target pairs. In the transition cases analyzed here, newly admitted routing contributions enter as \(1/66\) terms for the intermediate family identified on the corresponding admissible minimum-cost route.
The winner changes can therefore be diagnosed by counting which new pair terms are gained by C11a and which are gained by C32, as summarized in Figure~\ref{fig:section6_betweenness_pair_term_deltas}.

\par\medskip
\begin{center}
\captionsetup{type=figure}
    \centering
    \includegraphics[width=0.72\textwidth]{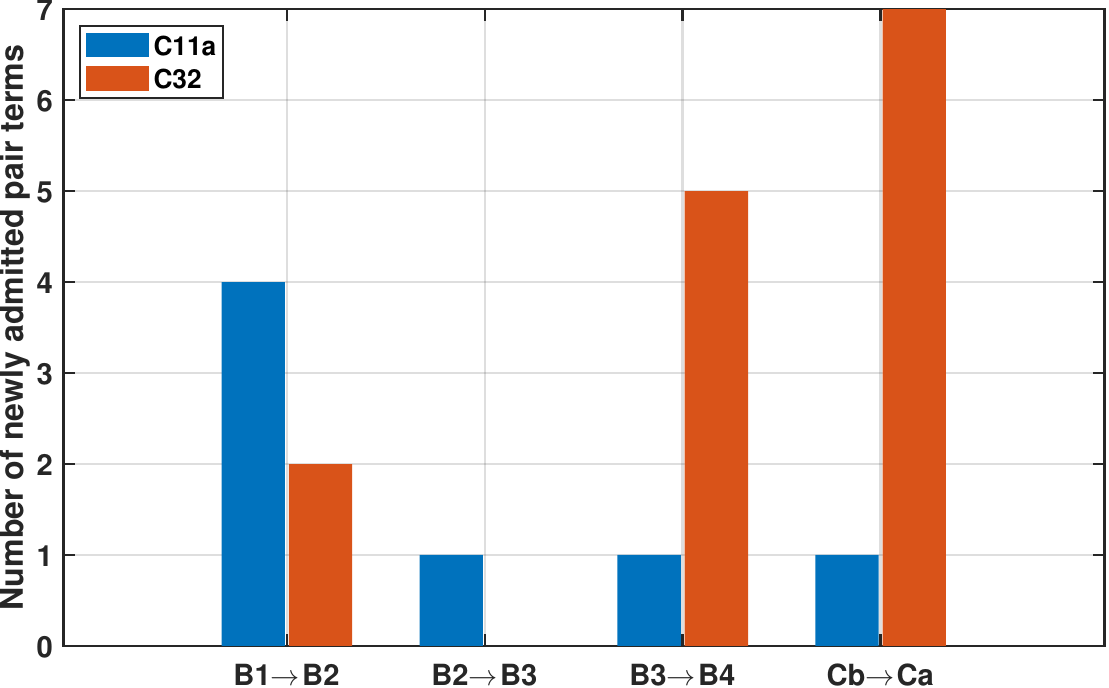}
    \caption{\footnotesize Betweenness winner changes explained by newly admitted \(1/66\) source--target pair terms.
    The four groups correspond to the diagnostic transitions discussed in the text:
    \(B_1\!\rightarrow\!B_2\), \(B_2\!\rightarrow\!B_3\), and \(B_3\!\rightarrow\!B_4\) at fixed \(T_{\mathrm{cap}}=10.25~\mathrm{days}\), and \(C_b\!\rightarrow\!C_a\) at fixed \(\Delta V_{\mathrm{cap}}=125.76~\mathrm{m/s}\).
    The bars show how many newly admitted pair terms are gained by C11a and C32 in each transition.
    The betweenness winner changes are therefore explained directly at the level of the source--target pair contributions, not merely by whether a new direct edge appears.}
    \label{fig:section6_betweenness_pair_term_deltas}
\end{center}
\medskip

At fixed \(T_{\mathrm{cap}}=10.25~\mathrm{days}\), increasing \(\Delta V_{\mathrm{cap}}\) through the sequence \(95.92 \rightarrow 110.84 \rightarrow 125.76 \rightarrow 140.68~\mathrm{m/s}\) produces the winner sequence
\[
C32
\;\rightarrow\;
\mathrm{tie}
\;\rightarrow\;
C11a
\;\rightarrow\;
C32\]
The first step, from \(95.92\) to \(110.84~\mathrm{m/s}\), creates an exact tie because the newly admitted routing structure generates four new pair terms mediated by C11a but only two new pair terms mediated by C32.
The second step, from \(110.84\) to \(125.76~\mathrm{m/s}\), breaks that tie in favor of C11a even though no new direct transfer appears.
Here the direct graph is unchanged, but one additional source--target pair mediated by C11a becomes admissible, and that single new \(1/66\) term is sufficient to reverse the winner.
The third step, from \(125.76\) to \(140.68~\mathrm{m/s}\), restores C32 as the relay winner because the expanded admissible routing structure contributes one new pair term mediated by C11a but five new pair terms mediated by C32.

A complementary time-driven betweenness transition occurs at fixed \(\Delta V_{\mathrm{cap}}=125.76~\mathrm{m/s}\) as \(T_{\mathrm{cap}}\) increases from \(11.10\) to \(11.95~\mathrm{days}\).
Again, the direct graph itself does not change.
Nevertheless, the set of admissible minimum-cost routes expands enough to add one new pair term mediated by C11a but seven new pair terms mediated by C32.
The relay winner therefore flips from C11a to C32 without any topology change in the network.

These diagnostics show that the winner transitions are controlled by different mechanisms for the three centrality measures used in the present study.
For strength and harmonic closeness, the transition is governed by the reweighting of a few dominant reciprocal-access terms, especially those tied to C32's access to the R52-U-side structure.
For betweenness, the transition is governed by the admission of individual source--target pair terms, and those terms may appear either because new direct transfers are admitted or because the same direct graph supports additional budget-feasible minimum-cost routes.
Betweenness is therefore more sensitive than strength or harmonic closeness to localized changes in admissible routing structure.

\section{Differential Correction Formulation}\label{app:correction}
This appendix details the local differential correction used in Section~\ref{sec:validation} to convert a proxy-identified overlap voxel into a continuous patched trajectory between two representative periodic orbits.

The correction operates on a six-component decision vector,
\begin{equation}
\mathbf{z} = (\alpha_A, \delta_A, t_A, \alpha_B, \delta_B, t_B),
\label{eq:decision_vector}
\end{equation}
where \(\alpha_A\) and \(\alpha_B\) are the phase locations on the source and target representative orbits, \(\delta_A\) and \(\delta_B\) are the heading changes applied at those phases, and \(t_A\) and \(t_B\) are the arc times from each orbit to the patch point.
The two arcs are integrated forward from the source side and backward from the target side under the reduced \((x,y,\theta)\) dynamics introduced in Section~\ref{subsec:reduced}, and the correction adjusts \(\mathbf{z}\) so that they meet at a common patch point.

Continuity is enforced through the patch residual,
\begin{equation}
\mathbf{r}(\mathbf{z}) =
\begin{bmatrix}
x_A(t_A) - x_B(t_B) \\
y_A(t_A) - y_B(t_B) \\
\operatorname{wrap}\!\left(\theta_A(t_A) - \theta_B(t_B)\right)
\end{bmatrix},
\label{eq:patch_residual}
\end{equation}
where \(\operatorname{wrap}(\cdot)\) maps a heading difference to \((-\pi, \pi]\) so that small differences across the angular discontinuity are not penalized as if they were large.
The three residual components are normalized by the atlas resolution scales used in the reachable-set construction (Table~\ref{tab:atlas_numerics}):
\begin{equation}
\mathbf{r}_{\mathrm{sc}}(\mathbf{z}) =
\begin{bmatrix}
\dfrac{x_A(t_A) - x_B(t_B)}{\Delta x} \\[0.8em]
\dfrac{y_A(t_A) - y_B(t_B)}{\Delta y} \\[0.8em]
\dfrac{\operatorname{wrap}\!\left(\theta_A(t_A) - \theta_B(t_B)\right)}{\Delta\theta}
\end{bmatrix}.
\label{eq:scaled_residual}
\end{equation}
A candidate is regarded as converged when
\begin{equation}
\left\|\mathbf{r}_{\mathrm{sc}}(\mathbf{z})\right\|_2 \le \varepsilon_{\mathrm{conv}},
\label{eq:convergence_criterion}
\end{equation}
so that continuity is enforced relative to the same spatial and angular resolution scales used to construct the reachable-set overlap atlases.
The convergence tolerance is set to $\varepsilon_{\mathrm{conv}} = 10^{-5}$ in the present study.

The objective minimized at convergence is the total boundary turning effort,
\begin{equation}
J(\mathbf{z}) = \Delta V_{\mathrm{turn},A}(\delta_A) + \Delta V_{\mathrm{turn},B}(\delta_B),
\label{eq:correction_objective}
\end{equation}
using the maneuver-cost expression of \eqref{eq:dv_turn}.
The corrected solution is obtained by Newton iteration the constrained nonlinear problem of minimizing \(J(\mathbf{z})\) subject to \(\mathbf{r}(\mathbf{z}) = \mathbf{0}\), with the warm start of Section~\ref{subsec:from_proxy} as initial guess.
The solver is terminated when \(\|\mathbf{r}_{\mathrm{sc}}(\mathbf{z})\|_2 \le \varepsilon_{\mathrm{conv}}\) and the gradient of \(J\) in the null space of the constraint Jacobian falls below a corresponding tolerance.

Both direct examples in Section~\ref{subsec:direct_examples} converged, and the resulting corrected \(\Delta V\) values are reported in Table~\ref{tab:section7_direct_examples}.

{\ifmarked\color{blue}\fi
\section{Large Time-Budget Limit of the Proxy Cost}\label{app:longtof}

The budget plane of Section~\ref{sec:budget} is sampled at finite
\((\Delta V_{\mathrm{cap}},T_{\mathrm{cap}})\). The behavior of the base turning cost as the
time budget is relaxed is informative because it separates the representative set into
two dynamically distinct classes and provides the sharpest available support for the
geometric reading given in the Conclusion.

Consider a pair of unstable representatives whose neighborhoods lie in a common
connected chaotic region of the \(C_J=3.1294\) manifold. As \(T_{\mathrm{a}}\) grows, the
forward atlas of the departure family and the backward atlas of the arrival family each
spread through that region along the invariant manifolds of the constituent unstable
orbits. Where an intersection of the unstable manifold of one family with the stable
manifold of the other is transverse, the resulting heteroclinic connection is an
asymptotic trajectory requiring no impulse at the patch point. The minimum base turning cost of
such a pair is therefore expected to decrease monotonically with time budget and to reach
zero at a finite \(T_{\mathrm{cap}}\), beyond which the pair is connected ballistically.

Stable resonant representatives behave differently. They lie on invariant tori inside
resonance islands that are dynamically isolated from the surrounding chaotic region, so no
asymptotic connection from the chaotic region reaches them. Their base turning cost instead
approaches a nonzero floor set by the impulse required to cross the island boundary, which
scales with the island width and therefore with the order of the resonance.

Sweeps of the one-sided propagation horizon confirm both behaviors for the
base turning cost \(\Delta V_{\mathrm{base}}\). This quantity isolates
the heading-change maneuvers associated with the underlying dynamics, whereas
the full proxy cost additionally contains the voxel-scale correction
\(\Delta V_{\mathrm{patch}}\) associated with the finite
\((\Delta x,\Delta y,\Delta\theta)\) discretization. Transfers from C32 to
unstable targets attain \(\Delta V_{\mathrm{base}}=0\) at finite horizon,
after which the base-cost Pareto front terminates. Transfers to the three
stable resonant representatives instead flatten toward the floors listed in
Table~\ref{tab:longtof_floors}. These behaviors are shown in
Figure~\ref{fig:longtof_sweeps}.

\par\medskip
\begin{center}
\captionsetup{type=figure}
\centering

\begin{subfigure}[t]{0.48\textwidth}
    \centering
    \includegraphics[width=\textwidth]{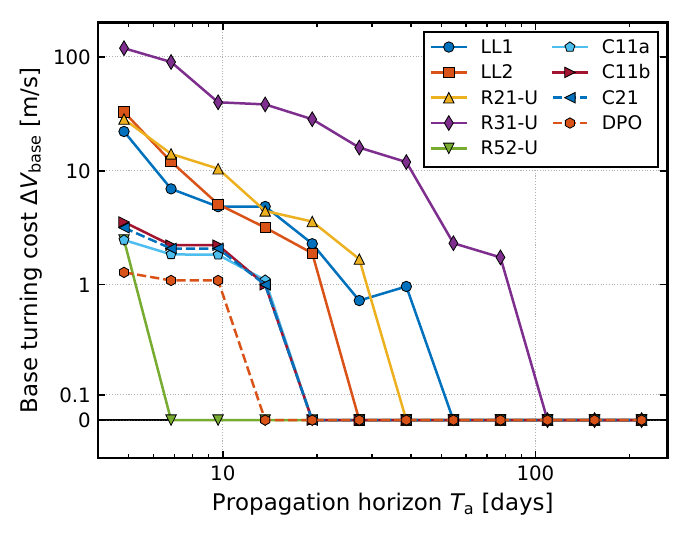}
    \caption{\footnotesize
    Unstable representatives.}
    \label{fig:longtof_unstable}
\end{subfigure}
\hfill
\begin{subfigure}[t]{0.48\textwidth}
    \centering
    \includegraphics[width=\textwidth]{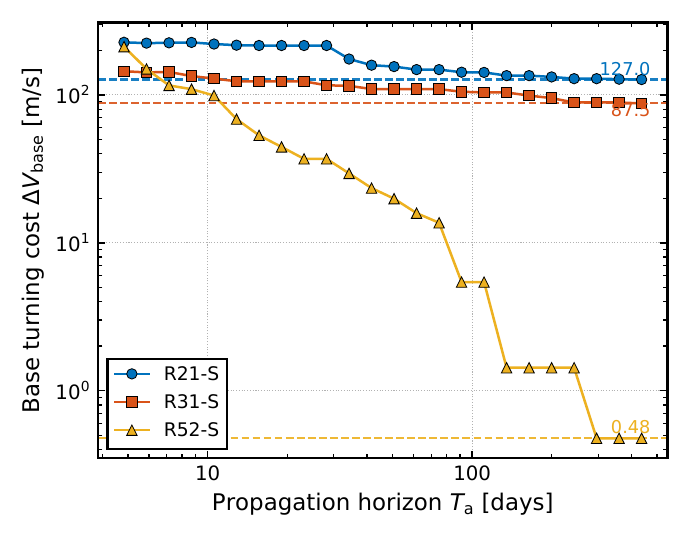}
    \caption{\footnotesize
    Stable resonant representatives.}
    \label{fig:longtof_stable}
\end{subfigure}

\caption{\footnotesize
Base turning contribution \(\Delta V_{\mathrm{base}}\) as a function of the
one-sided propagation horizon \(T_{\mathrm{a}}\) for transfers from C32.
(a) For the nine unstable representatives,
\(\Delta V_{\mathrm{base}}\) reaches zero at a finite propagation horizon.
(b) For the three stable resonant representatives,
\(\Delta V_{\mathrm{base}}\) approaches nonzero large-time floors.
}
\label{fig:longtof_sweeps}
\end{center}
\medskip

\par\medskip
\begin{center}
\captionsetup{type=table}
\centering
\caption{\footnotesize Approximate large-time-budget floors of the base turning cost \(\Delta V_{\mathrm{base}}\) for the
stable resonant representatives. Unstable representatives reach zero at finite horizon.}
\label{tab:longtof_floors}
\begin{tabular}{lcc}
\toprule
Target family & Abbreviation & Asymptotic floor [m/s] \\
\midrule
2:1 stable resonant & R21-S & 127.0 \\
3:1 stable resonant & R31-S & 87.5 \\
5:2 stable resonant & R52-S & 0.48 \\
\bottomrule
\end{tabular}
\end{center}
\medskip

Two consequences follow. First, the low direct cost recovered for the LL1--LL2 pair in
Section~\ref{subsec:matrix} is the finite-budget signature of a connection that becomes exactly
ballistic in this limit, which is why that pair is recovered at a cost of order
\(10~\mathrm{m/s}\) rather than at a cost comparable to the resonant targets. Second, the
dichotomy between vanishing and floored asymptotic cost is a direct expression of the
partition of the representative set into a connected chaotic component and a set of isolated
resonant tori, and it is independent of the particular voxel resolution used here.

Two limitations should be noted. The floors in Table~\ref{tab:longtof_floors} are obtained
numerically and are not established analytically; in particular, the R52-S floor is
substantially smaller than a first-order island-width estimate would suggest, and the
discrepancy is the subject of ongoing work. The argument above also assumes transversality of
the relevant manifold intersections, which is generic but not verified pair by pair here.
}

\section*{Acknowledgments}

S.D.R. gratefully acknowledges support from
the Air Force Office of Scientific Research (AFOSR) under
Grant No.\ FA9550-24–1-0194.
The authors acknowledge Advanced Research Computing at Virginia Tech for
providing computational resources that contributed to the results reported
in this paper.

\section*{Data and Code Availability Statement}

The MATLAB source code used in this work, together with the family-to-family accessibility matrices, budget-plane diagnostics, and corrected-transfer summaries reported in this paper, is archived on Zenodo (release v1.0.0, \url{https://doi.org/10.5281/zenodo.24022012}) and maintained at \url{https://github.com/binbraik/cislunar_orbital_network}. 
The full reachable-set atlases (several hundred GB) are not archived; they are regenerated deterministically by the archived code. 
The interactive $\Delta V$-map visualization is hosted at \url{https://binbraik.github.io/cislunar_orbital_network}. 
Video supplements associated with Figure~\ref{fig:section7_relay_examples} are referenced in the figure captions and hosted on YouTube.

\bibliographystyle{cas-model2-names}
\bibliography{references}

\begin{thebibliography}{48}
\expandafter\ifx\csname natexlab\endcsname\relax\def\natexlab#1{#1}\fi
\providecommand{\url}[1]{\texttt{#1}}
\providecommand{\href}[2]{#2}
\providecommand{\path}[1]{#1}
\providecommand{\DOIprefix}{doi:}
\providecommand{\ArXivprefix}{arXiv:}
\providecommand{\URLprefix}{URL: }
\providecommand{\Pubmedprefix}{pmid:}
\providecommand{\doi}[1]{\href{http://dx.doi.org/#1}{\path{#1}}}
\providecommand{\Pubmed}[1]{\href{pmid:#1}{\path{#1}}}
\providecommand{\bibinfo}[2]{#2}
\ifx\xfnm\relax \def\xfnm[#1]{\unskip,\space#1}\fi
\bibitem[{Barrat et~al.(2004)Barrat, Barth\'{e}lemy, Pastor-Satorras and Vespignani}]{barrat2004}
\bibinfo{author}{Barrat, A.}, \bibinfo{author}{Barth\'{e}lemy, M.}, \bibinfo{author}{Pastor-Satorras, R.}, \bibinfo{author}{Vespignani, A.}, \bibinfo{year}{2004}.
\newblock \bibinfo{title}{The architecture of complex weighted networks}.
\newblock \bibinfo{journal}{Proceedings of the National Academy of Sciences} \bibinfo{volume}{101}, \bibinfo{pages}{3747--3752}.
\newblock \DOIprefix\doi{10.1073/pnas.0400087101}.
\bibitem[{Bate et~al.(1971)Bate, Mueller and White}]{BaMuWh1971}
\bibinfo{author}{Bate, R.R.}, \bibinfo{author}{Mueller, D.D.}, \bibinfo{author}{White, J.E.}, \bibinfo{year}{1971}.
\newblock \bibinfo{title}{Fundamentals of Astrodynamics}.
\newblock \bibinfo{publisher}{Dover}, \bibinfo{address}{New York}.
\bibitem[{Boldi and Vigna(2014)}]{boldi2014}
\bibinfo{author}{Boldi, P.}, \bibinfo{author}{Vigna, S.}, \bibinfo{year}{2014}.
\newblock \bibinfo{title}{Axioms for centrality}.
\newblock \bibinfo{journal}{Internet Mathematics} \bibinfo{volume}{10}, \bibinfo{pages}{222--262}.
\newblock \DOIprefix\doi{10.1080/15427951.2013.865686}.
\bibitem[{Bowerfind and Taheri(2024)}]{bowerfind2024}
\bibinfo{author}{Bowerfind, S.R.}, \bibinfo{author}{Taheri, E.}, \bibinfo{year}{2024}.
\newblock \bibinfo{title}{Rapid approximation of low-thrust spacecraft reachable sets within complex two-body and cislunar dynamics}.
\newblock \bibinfo{journal}{Aerospace} \bibinfo{volume}{11}, \bibinfo{pages}{380}.
\newblock \DOIprefix\doi{10.3390/aerospace11050380}.
\bibitem[{Braik and Ross(2025)}]{braik2025}
\bibinfo{author}{Braik, A.}, \bibinfo{author}{Ross, S.D.}, \bibinfo{year}{2025}.
\newblock \bibinfo{title}{Heteroclinic transfer between {L1} and {L3} in the {Earth-Moon} system}, in: \bibinfo{booktitle}{AAS/AIAA Astrodynamics Specialist Conference}, \bibinfo{address}{Boston, MA}.
\bibitem[{Broucke(1968)}]{broucke1968}
\bibinfo{author}{Broucke, R.A.}, \bibinfo{year}{1968}.
\newblock \bibinfo{title}{Periodic Orbits in the Restricted Three-Body Problem with {Earth--Moon} Masses}.
\newblock \bibinfo{type}{Technical Report} \bibinfo{number}{JPL Technical Report 32-1168}. Jet Propulsion Laboratory.
\bibitem[{Bruchko and Bosanac(2025)}]{bruchko2025}
\bibinfo{author}{Bruchko, K.L.}, \bibinfo{author}{Bosanac, N.}, \bibinfo{year}{2025}.
\newblock \bibinfo{title}{Rapid trajectory design in multi-body systems using sampling-based kinodynamic planning}.
\newblock \bibinfo{journal}{The Journal of the Astronautical Sciences} \bibinfo{volume}{72}, \bibinfo{pages}{33}.
\newblock \DOIprefix\doi{10.1007/s40295-025-00506-6}.
\bibitem[{Capdevila and Howell(2018)}]{capdevila2018}
\bibinfo{author}{Capdevila, L.R.}, \bibinfo{author}{Howell, K.C.}, \bibinfo{year}{2018}.
\newblock \bibinfo{title}{A transfer network linking {Earth}, {Moon}, and the triangular libration point regions in the {Earth--Moon} system}.
\newblock \bibinfo{journal}{Advances in Space Research} \bibinfo{volume}{62}, \bibinfo{pages}{1826--1852}.
\newblock \DOIprefix\doi{10.1016/j.asr.2018.06.045}.
\bibitem[{Cheetham et~al.(2022)Cheetham, Gardner, Kayser, Parrish et~al.}]{cheetham2022}
\bibinfo{author}{Cheetham, B.W.}, \bibinfo{author}{Gardner, T.}, \bibinfo{author}{Kayser, E.}, \bibinfo{author}{Parrish, N.L.}, et~al., \bibinfo{year}{2022}.
\newblock \bibinfo{title}{{CAPSTONE}: A {CubeSat} pathfinder for the lunar {Gateway} ecosystem}, in: \bibinfo{booktitle}{36th Annual Small Satellite Conference}.
\bibitem[{Davis et~al.(2017)Davis, Phillips and Howell}]{davis2017}
\bibinfo{author}{Davis, D.C.}, \bibinfo{author}{Phillips, S.M.}, \bibinfo{author}{Howell, K.C.}, \bibinfo{year}{2017}.
\newblock \bibinfo{title}{Stationkeeping and transfer trajectory design for spacecraft in cislunar space}, in: \bibinfo{booktitle}{AAS/AIAA Astrodynamics Specialist Conference}.
\bibitem[{Dellnitz et~al.(2005a)Dellnitz, Junge, Koon, Lekien, Lo, Marsden, Padberg, Preis, Ross and Thiere}]{dellnitz2005}
\bibinfo{author}{Dellnitz, M.}, \bibinfo{author}{Junge, O.}, \bibinfo{author}{Koon, W.S.}, \bibinfo{author}{Lekien, F.}, \bibinfo{author}{Lo, M.W.}, \bibinfo{author}{Marsden, J.E.}, \bibinfo{author}{Padberg, K.}, \bibinfo{author}{Preis, R.}, \bibinfo{author}{Ross, S.D.}, \bibinfo{author}{Thiere, B.}, \bibinfo{year}{2005}a.
\newblock \bibinfo{title}{Transport in dynamical astronomy and multibody problems}.
\newblock \bibinfo{journal}{International Journal of Bifurcation and Chaos} \bibinfo{volume}{15}, \bibinfo{pages}{699--727}.
\newblock \DOIprefix\doi{10.1142/S0218127405012545}.
\bibitem[{Dellnitz et~al.(2005b)Dellnitz, Junge, Lo, Marsden, Padberg, Preis, Ross and Thiere}]{DeJuLoMaPaPrRoTh2005}
\bibinfo{author}{Dellnitz, M.}, \bibinfo{author}{Junge, O.}, \bibinfo{author}{Lo, M.W.}, \bibinfo{author}{Marsden, J.E.}, \bibinfo{author}{Padberg, K.}, \bibinfo{author}{Preis, R.}, \bibinfo{author}{Ross, S.D.}, \bibinfo{author}{Thiere, B.}, \bibinfo{year}{2005}b.
\newblock \bibinfo{title}{Transport of {M}ars-crossing asteroids from the quasi-{H}ilda region}.
\newblock \bibinfo{journal}{Physical Review Letters} \bibinfo{volume}{94}, \bibinfo{pages}{231102}.
\bibitem[{Dellnitz et~al.(2006)Dellnitz, Junge, Post and Thiere}]{dellnitz2006target}
\bibinfo{author}{Dellnitz, M.}, \bibinfo{author}{Junge, O.}, \bibinfo{author}{Post, M.}, \bibinfo{author}{Thiere, B.}, \bibinfo{year}{2006}.
\newblock \bibinfo{title}{On target for venus--set oriented computation of energy efficient low thrust trajectories}.
\newblock \bibinfo{journal}{Celestial Mechanics and Dynamical Astronomy} \bibinfo{volume}{95}, \bibinfo{pages}{357--370}.
\bibitem[{Dijkstra(1959)}]{dijkstra1959}
\bibinfo{author}{Dijkstra, E.W.}, \bibinfo{year}{1959}.
\newblock \bibinfo{title}{A note on two problems in connexion with graphs}.
\newblock \bibinfo{journal}{Numerische Mathematik} \bibinfo{volume}{1}, \bibinfo{pages}{269--271}.
\newblock \DOIprefix\doi{10.1007/BF01386390}.
\bibitem[{Doedel et~al.(2007)Doedel, Romanov, Paffenroth, Keller, Dichmann, Gal\'{a}n-Vioque and Vanderbauwhede}]{doedel2007}
\bibinfo{author}{Doedel, E.J.}, \bibinfo{author}{Romanov, V.A.}, \bibinfo{author}{Paffenroth, R.C.}, \bibinfo{author}{Keller, H.B.}, \bibinfo{author}{Dichmann, D.J.}, \bibinfo{author}{Gal\'{a}n-Vioque, J.}, \bibinfo{author}{Vanderbauwhede, A.}, \bibinfo{year}{2007}.
\newblock \bibinfo{title}{Elemental periodic orbits associated with the libration points in the circular restricted 3-body problem}.
\newblock \bibinfo{journal}{International Journal of Bifurcation and Chaos} \bibinfo{volume}{17}, \bibinfo{pages}{2625--2677}.
\newblock \DOIprefix\doi{10.1142/S0218127407018671}.
\bibitem[{Freeman(1977)}]{freeman1977}
\bibinfo{author}{Freeman, L.C.}, \bibinfo{year}{1977}.
\newblock \bibinfo{title}{A set of measures of centrality based on betweenness}.
\newblock \bibinfo{journal}{Sociometry} \bibinfo{volume}{40}, \bibinfo{pages}{35--41}.
\newblock \DOIprefix\doi{10.2307/3033543}.
\bibitem[{G\'{o}mez et~al.(2004)G\'{o}mez, Koon, Lo, Marsden, Masdemont and Ross}]{gomez2004}
\bibinfo{author}{G\'{o}mez, G.}, \bibinfo{author}{Koon, W.S.}, \bibinfo{author}{Lo, M.W.}, \bibinfo{author}{Marsden, J.E.}, \bibinfo{author}{Masdemont, J.}, \bibinfo{author}{Ross, S.D.}, \bibinfo{year}{2004}.
\newblock \bibinfo{title}{Connecting orbits and invariant manifolds in the spatial restricted three-body problem}.
\newblock \bibinfo{journal}{Nonlinearity} \bibinfo{volume}{17}, \bibinfo{pages}{1571--1606}.
\newblock \DOIprefix\doi{10.1088/0951-7715/17/5/002}.
\bibitem[{Guzzetti et~al.(2016)Guzzetti, Bosanac, Haapala, Howell and Folta}]{guzzetti2016}
\bibinfo{author}{Guzzetti, D.}, \bibinfo{author}{Bosanac, N.}, \bibinfo{author}{Haapala, A.}, \bibinfo{author}{Howell, K.C.}, \bibinfo{author}{Folta, D.C.}, \bibinfo{year}{2016}.
\newblock \bibinfo{title}{Rapid trajectory design in the {Earth--Moon} ephemeris system via an interactive catalog of periodic and quasi-periodic orbits}.
\newblock \bibinfo{journal}{Acta Astronautica} \bibinfo{volume}{126}, \bibinfo{pages}{439--455}.
\newblock \DOIprefix\doi{10.1016/j.actaastro.2016.06.009}.
\bibitem[{Haapala and Howell(2016)}]{haapala2016}
\bibinfo{author}{Haapala, A.F.}, \bibinfo{author}{Howell, K.C.}, \bibinfo{year}{2016}.
\newblock \bibinfo{title}{A framework for constructing transfers linking periodic libration point orbits in the spatial circular restricted three-body problem}.
\newblock \bibinfo{journal}{International Journal of Bifurcation and Chaos} \bibinfo{volume}{26}, \bibinfo{pages}{1630013}.
\newblock \DOIprefix\doi{10.1142/S0218127416300135}.
\bibitem[{H\'{e}non(1969)}]{henon1969}
\bibinfo{author}{H\'{e}non, M.}, \bibinfo{year}{1969}.
\newblock \bibinfo{title}{Numerical exploration of the restricted problem. {V}. {Hill's} case: Periodic orbits and their stability}.
\newblock \bibinfo{journal}{Astronomy \& Astrophysics} \bibinfo{volume}{1}, \bibinfo{pages}{223--238}.
\bibitem[{H\'{e}non(1997)}]{henon1997}
\bibinfo{author}{H\'{e}non, M.}, \bibinfo{year}{1997}.
\newblock \bibinfo{title}{Generating Families in the Restricted Three-Body Problem}. volume~\bibinfo{volume}{52} of \textit{\bibinfo{series}{Lecture Notes in Physics Monographs}}.
\newblock \bibinfo{publisher}{Springer}, \bibinfo{address}{Berlin}.
\bibitem[{Hiraiwa et~al.(2026)Hiraiwa, Bando, Sato and Hokamoto}]{hiraiwa2026}
\bibinfo{author}{Hiraiwa, N.}, \bibinfo{author}{Bando, M.}, \bibinfo{author}{Sato, Y.}, \bibinfo{author}{Hokamoto, S.}, \bibinfo{year}{2026}.
\newblock \bibinfo{title}{Design of low-energy transfers in cislunar space using sequences of lobe dynamics}.
\newblock \bibinfo{journal}{arXiv preprint arXiv:2602.17444} \URLprefix \url{https://arxiv.org/abs/2602.17444}, \DOIprefix\doi{10.48550/arXiv.2602.17444}, \href{http://arxiv.org/abs/2602.17444}{\tt arXiv:2602.17444}.
\bibitem[{Howell(1984)}]{howell1984}
\bibinfo{author}{Howell, K.C.}, \bibinfo{year}{1984}.
\newblock \bibinfo{title}{Three-dimensional, periodic, {`Halo'} orbits}.
\newblock \bibinfo{journal}{Celestial Mechanics} \bibinfo{volume}{32}, \bibinfo{pages}{53--71}.
\newblock \DOIprefix\doi{10.1007/BF01358403}.
\bibitem[{Jerg et~al.(2009)Jerg, Junge and Ross}]{JeJuRo2009}
\bibinfo{author}{Jerg, S.}, \bibinfo{author}{Junge, O.}, \bibinfo{author}{Ross, S.D.}, \bibinfo{year}{2009}.
\newblock \bibinfo{title}{Optimal capture trajectories using multiple gravity assists}.
\newblock \bibinfo{journal}{Communications in Nonlinear Science and Numerical Simulation} \bibinfo{volume}{14}, \bibinfo{pages}{4168 -- 4175}.
\newblock \DOIprefix\doi{DOI: 10.1016/j.cnsns.2008.12.009}.
\bibitem[{{Jet Propulsion Laboratory}(2025)}]{JPL_PeriodicOrbits}
\bibinfo{author}{{Jet Propulsion Laboratory}}, \bibinfo{year}{2025}.
\newblock \bibinfo{title}{{Three-Body Periodic Orbits}}.
\newblock \bibinfo{howpublished}{JPL Solar System Dynamics, \url{https://ssd.jpl.nasa.gov/tools/periodic_orbits.html}}.
\bibitem[{Koon et~al.(2000)Koon, Lo, Marsden and Ross}]{koon2000}
\bibinfo{author}{Koon, W.S.}, \bibinfo{author}{Lo, M.W.}, \bibinfo{author}{Marsden, J.E.}, \bibinfo{author}{Ross, S.D.}, \bibinfo{year}{2000}.
\newblock \bibinfo{title}{Heteroclinic connections between periodic orbits and resonance transitions in celestial mechanics}.
\newblock \bibinfo{journal}{Chaos} \bibinfo{volume}{10}, \bibinfo{pages}{427--469}.
\newblock \DOIprefix\doi{10.1063/1.166509}.
\bibitem[{Koon et~al.(2011)Koon, Lo, Marsden and Ross}]{koon2011}
\bibinfo{author}{Koon, W.S.}, \bibinfo{author}{Lo, M.W.}, \bibinfo{author}{Marsden, J.E.}, \bibinfo{author}{Ross, S.D.}, \bibinfo{year}{2011}.
\newblock \bibinfo{title}{Dynamical Systems, the Three-Body Problem and Space Mission Design}.
\newblock \bibinfo{publisher}{Marsden Books}.
\bibitem[{Kumar et~al.(2026)Kumar, Rawat, Rosengren and Ross}]{kumar2006}
\bibinfo{author}{Kumar, B.}, \bibinfo{author}{Rawat, A.}, \bibinfo{author}{Rosengren, A.J.}, \bibinfo{author}{Ross, S.D.}, \bibinfo{year}{2026}.
\newblock \bibinfo{title}{Cislunar resonant transport and heteroclinic pathways: From 3:1 to 2:1 to l1}.
\newblock \bibinfo{journal}{Advances in Space Research} \bibinfo{volume}{77}, \bibinfo{pages}{3815--3843}.
\newblock \DOIprefix\doi{https://doi.org/10.1016/j.asr.2025.12.005}.
\bibitem[{Latora and Marchiori(2001)}]{latora2001}
\bibinfo{author}{Latora, V.}, \bibinfo{author}{Marchiori, M.}, \bibinfo{year}{2001}.
\newblock \bibinfo{title}{Efficient behavior of small-world networks}.
\newblock \bibinfo{journal}{Physical Review Letters} \bibinfo{volume}{87}, \bibinfo{pages}{198701}.
\newblock \DOIprefix\doi{10.1103/PhysRevLett.87.198701}.
\bibitem[{Leiva and Briozzo(2006)}]{leiva2006control}
\bibinfo{author}{Leiva, A.}, \bibinfo{author}{Briozzo, C.}, \bibinfo{year}{2006}.
\newblock \bibinfo{title}{Control of chaos and fast periodic transfer orbits in the {E}arth-{M}oon {CR3BP}}.
\newblock \bibinfo{journal}{Acta Astronautica} \bibinfo{volume}{58}, \bibinfo{pages}{379--386}.
\bibitem[{Marchiori and Latora(2000)}]{marchiori2000}
\bibinfo{author}{Marchiori, M.}, \bibinfo{author}{Latora, V.}, \bibinfo{year}{2000}.
\newblock \bibinfo{title}{Harmony in the small-world}.
\newblock \bibinfo{journal}{Physica A} \bibinfo{volume}{285}, \bibinfo{pages}{539--546}.
\newblock \DOIprefix\doi{10.1016/S0378-4371(00)00311-3}.
\bibitem[{McCarthy and Howell(2023)}]{mccarthy2023}
\bibinfo{author}{McCarthy, B.P.}, \bibinfo{author}{Howell, K.C.}, \bibinfo{year}{2023}.
\newblock \bibinfo{title}{Construction of heteroclinic connections between quasi-periodic orbits in the three-body problem}.
\newblock \bibinfo{journal}{The Journal of the Astronautical Sciences} \bibinfo{volume}{70}, \bibinfo{pages}{24}.
\newblock \DOIprefix\doi{10.1007/s40295-023-00389-5}.
\bibitem[{Newman(2010)}]{newman2010}
\bibinfo{author}{Newman, M.E.J.}, \bibinfo{year}{2010}.
\newblock \bibinfo{title}{Networks: An Introduction}.
\newblock \bibinfo{publisher}{Oxford University Press}, \bibinfo{address}{Oxford}.
\bibitem[{Parker and Anderson(2014)}]{parker2014}
\bibinfo{author}{Parker, J.S.}, \bibinfo{author}{Anderson, R.L.}, \bibinfo{year}{2014}.
\newblock \bibinfo{title}{Low-Energy Lunar Trajectory Design}. volume~\bibinfo{volume}{12} of \textit{\bibinfo{series}{JPL Deep-Space Communications and Navigation Series}}.
\newblock \bibinfo{publisher}{Wiley}, \bibinfo{address}{Hoboken, NJ}.
\bibitem[{Rawat et~al.(2026)Rawat, Kumar, Rosengren and Ross}]{rawat2026}
\bibinfo{author}{Rawat, A.}, \bibinfo{author}{Kumar, B.}, \bibinfo{author}{Rosengren, A.J.}, \bibinfo{author}{Ross, S.D.}, \bibinfo{year}{2026}.
\newblock \bibinfo{title}{Cislunar mean-motion resonances: Definitions, widths, and comparisons with resonant satellites}.
\newblock \bibinfo{journal}{Journal of Guidance, Control, and Dynamics} \bibinfo{volume}{49}, \bibinfo{pages}{1068--1082}.
\newblock \DOIprefix\doi{10.2514/1.G009336}.
\bibitem[{Ross and Roberts-Tsoukkas(2025)}]{ross2025cyclers}
\bibinfo{author}{Ross, S.D.}, \bibinfo{author}{Roberts-Tsoukkas, M.}, \bibinfo{year}{2025}.
\newblock \bibinfo{title}{Stable, low-energy prograde {Earth-Moon} cycler orbits}, in: \bibinfo{booktitle}{AAS/AIAA Astrodynamics Specialist Conference}, \bibinfo{address}{Boston, MA}.
\bibitem[{Ross and Roberts-Tsoukkas(2026)}]{ross2026cyclers}
\bibinfo{author}{Ross, S.D.}, \bibinfo{author}{Roberts-Tsoukkas, M.}, \bibinfo{year}{2026}.
\newblock \bibinfo{title}{Stable ballistic prograde cyclers in the three-body problem}.
\newblock \bibinfo{journal}{Phys. Rev. E} \URLprefix \url{https://link.aps.org/doi/10.1103/6jy8-bls1}, \DOIprefix\doi{10.1103/6jy8-bls1}, \href{http://arxiv.org/abs/2606.29189}{\tt arXiv:2606.29189}. \bibinfo{note}{accepted}.
\bibitem[{Saxena and Iyengar(2020)}]{SaxenaIyengar2020Survey}
\bibinfo{author}{Saxena, A.}, \bibinfo{author}{Iyengar, S.}, \bibinfo{year}{2020}.
\newblock \bibinfo{title}{Centrality measures in complex networks: A survey}.
\newblock \bibinfo{journal}{arXiv preprint arXiv:2011.07190} \DOIprefix\doi{10.48550/arXiv.2011.07190}.
\bibitem[{Schroer and Ott(1997)}]{ScOt1997}
\bibinfo{author}{Schroer, C.G.}, \bibinfo{author}{Ott, E.}, \bibinfo{year}{1997}.
\newblock \bibinfo{title}{Targeting in {H}amiltonian systems that have mixed regular/chaotic phase spaces}.
\newblock \bibinfo{journal}{Chaos} \bibinfo{volume}{7}, \bibinfo{pages}{512--519}.
\bibitem[{Smith and Bosanac(2022)}]{smith2022}
\bibinfo{author}{Smith, T.R.}, \bibinfo{author}{Bosanac, N.}, \bibinfo{year}{2022}.
\newblock \bibinfo{title}{Constructing motion primitive sets to summarize periodic orbit families and hyperbolic invariant manifolds in a multi-body system}.
\newblock \bibinfo{journal}{Celestial Mechanics and Dynamical Astronomy} \bibinfo{volume}{134}.
\newblock \DOIprefix\doi{10.1007/s10569-022-10063-x}.
\bibitem[{Smith and Bosanac(2023)}]{smith2023}
\bibinfo{author}{Smith, T.R.}, \bibinfo{author}{Bosanac, N.}, \bibinfo{year}{2023}.
\newblock \bibinfo{title}{Motion primitive approach to spacecraft trajectory design in a multi-body system}.
\newblock \bibinfo{journal}{The Journal of the Astronautical Sciences} \bibinfo{volume}{70}, \bibinfo{pages}{34}.
\newblock \DOIprefix\doi{10.1007/s40295-023-00395-7}.
\bibitem[{Spear et~al.(2026)Spear, Miceli and Bosanac}]{spear2026}
\bibinfo{author}{Spear, R.L.}, \bibinfo{author}{Miceli, G.E.}, \bibinfo{author}{Bosanac, N.}, \bibinfo{year}{2026}.
\newblock \bibinfo{title}{Tree-based approach to spacecraft trajectory design in multi-body systems}.
\newblock \bibinfo{journal}{The Journal of the Astronautical Sciences} \bibinfo{volume}{73}, \bibinfo{pages}{15}.
\newblock \DOIprefix\doi{10.1007/s40295-025-00563-x}.
\bibitem[{Szebehely(1967)}]{szebehely1967}
\bibinfo{author}{Szebehely, V.G.}, \bibinfo{year}{1967}.
\newblock \bibinfo{title}{Theory of Orbits: The Restricted Problem of Three Bodies}.
\newblock \bibinfo{publisher}{Academic Press}, \bibinfo{address}{New York}.
\bibitem[{Tsirogiannis(2012)}]{tsirogiannis2012}
\bibinfo{author}{Tsirogiannis, G.A.}, \bibinfo{year}{2012}.
\newblock \bibinfo{title}{A graph based methodology for mission design}.
\newblock \bibinfo{journal}{Celestial Mechanics and Dynamical Astronomy} \bibinfo{volume}{114}, \bibinfo{pages}{353--363}.
\newblock \DOIprefix\doi{10.1007/s10569-012-9444-9}.
\bibitem[{Vaquero and Howell(2014)}]{vaquero2014}
\bibinfo{author}{Vaquero, M.}, \bibinfo{author}{Howell, K.C.}, \bibinfo{year}{2014}.
\newblock \bibinfo{title}{Leveraging resonant-orbit manifolds to design transfers between libration-point orbits}.
\newblock \bibinfo{journal}{Journal of Guidance, Control, and Dynamics} \bibinfo{volume}{37}, \bibinfo{pages}{1143--1157}.
\newblock \DOIprefix\doi{10.2514/1.62230}.
\bibitem[{Whitley and Martinez(2016)}]{whitley2016}
\bibinfo{author}{Whitley, R.J.}, \bibinfo{author}{Martinez, R.}, \bibinfo{year}{2016}.
\newblock \bibinfo{title}{Options for staging orbits in cislunar space}, in: \bibinfo{booktitle}{2016 IEEE Aerospace Conference}, pp. \bibinfo{pages}{1--9}.
\newblock \DOIprefix\doi{10.1109/AERO.2016.7500635}.
\bibitem[{Zimovan(2017)}]{zimovan2017}
\bibinfo{author}{Zimovan, E.M.}, \bibinfo{year}{2017}.
\newblock \bibinfo{title}{Characteristics and Design Strategies for Near Rectilinear Halo Orbits Within the {Earth--Moon} System}.
\newblock Master's thesis. Purdue University.
\bibitem[{Zimovan-Spreen et~al.(2022)Zimovan-Spreen, Howell and Davis}]{zimovanspreen2022}
\bibinfo{author}{Zimovan-Spreen, E.M.}, \bibinfo{author}{Howell, K.C.}, \bibinfo{author}{Davis, D.C.}, \bibinfo{year}{2022}.
\newblock \bibinfo{title}{Dynamical structures nearby {NRHOs} with applications in cislunar space}.
\newblock \bibinfo{journal}{The Journal of the Astronautical Sciences} \bibinfo{volume}{69}, \bibinfo{pages}{718--744}.
\newblock \DOIprefix\doi{10.1007/s40295-021-00279-4}.

\end{thebibliography}

\section*{Declaration of competing interest}

The authors declare that they have no known competing
financial interests or personal relationships that could have
appeared to influence the work reported in this paper.

\section*{Declaration of generative AI and AI-assisted technologies in the manuscript preparation process}

During the preparation of this work, the authors used ChatGPT and Claude to assist with reviewing and refining code implementation, improving clarity and coherence in the manuscript structure, and editing language for readability. After using these tools, the authors reviewed, verified, and edited the content as needed and take full responsibility for the content of the published article.

\end{document}